\documentclass[12pt,reqno]{amsart}
\usepackage{amsmath,amsthm,amscd,amsfonts, amssymb,eepic}

\swapnumbers

\newtheorem{thm}[subsection]{Theorem}

\newtheorem{lem}[subsection]{Lemma}
\newtheorem{prop}[subsection]{Proposition}

\theoremstyle{definition}
\newtheorem{defn}[subsection]{Definition}

\numberwithin{equation}{subsection}
\newcommand{\parens}[1]{{\rm{(}}#1{\rm{)}}}
\newcommand{\quash}[1]{}
\newcommand{\qdr}[1]{_{\left\lceil\nu_P,{#1}\right\rfloor}}
\newcommand{\q}{\qbezier}

\newcommand{\alarge}{large }

\newcommand{\ao}{0 }
\newcommand{\ainfty}{\infty }

\newcommand{\age}{\ge }
\newcommand{\alt}{< }
\newcommand{\agt}{> }
\newcommand{\reals}{\mathbb R}
\newcommand{\Sum}{\sum}
\newcommand{\intersect}{\cap}

\begin{document}
\title{ The Topological Trace Formula}
\author{Mark Goresky\thanks{School of Mathematics, Institute for Advanced Study,
Princeton N.J.  Research partially supported by NSF grants \# DMS 9626616 and DMS
9900324.}\and
Robert MacPherson \thanks{School of Mathematics, Institute for Advanced Study, Princeton
N.J.} }

\begin{abstract}  The topological trace formula is a computation of  the
Lefschetz number of a Hecke correspondence $C$ acting on the weighted cohomology groups,
defined in \cite{GHM}, of a locally symmetric space $X$.  It expresses this Lefschetz number
as a sum of contributions from fixed point components of $C$ on the reductive Borel Serre
compactification of $X$.  The proof uses the Lefschetz fixed point formula of \cite{LFP}.
  
When the $L^2$ cohomology of $X$ is finite dimensional (the equal rank case),   the Lefschetz
number of a Hecke correspondence acting on the $L^2$ cohomology was computed by J. Arthur in
\cite{Ar1, Ar3} using the trace formula.  In this cased the $L^2$ cohomology coincides with
the ``middle'' weighted cohomology (see \cite{GHM}), so this paper gives an independent
computation of this Lefschetz number.  In \cite{GKM}, it was shown that these two computations
agree.  Consequently, this paper completes an independent proof of Arthur's formula.

We also describe some of the rich geometry of Hecke correspondences and their fixed points.
While \cite{Ar1} and \cite{GKM}, because of their adelic setting, treated congruence
subgroups, we treat here the slightly more general case of an arbitrary neat arithmetic
subgroup.  Some of the results of this paper were announced in \cite{Montreal}.
\end{abstract}
\maketitle

\section{Introduction.}
\subsection{Geometric overview.}  

\subsubsection{Locally symmetric spaces.}\label{Locally symmetric spaces}  For the purposes of
this introduction, a {\it locally symmetric
space} $X$ is a complete connected Riemannian manifold with finite volume and non-positive
curvature, such that every point $p\in X$ has a neighborhood $U_p$ with a {\em Cartan
symmetry:}  an isometry $U_p \to U_p$ that takes $p$ to itself and induces minus the identity
on the tangent space to $X$ at $p$.  As for any manifold, we have $X=\Gamma \backslash D$
where $D$ is the simply connected covering space of $X$ and $\Gamma$ is the fundamental group
of $X$.  Because $X$ has nonpositive curvature, it follows that $D$ is a Riemannian symmetric
space of noncompact type, that is,  
the metric product of a negatively curved symmetric space from Cartan's classification, and a
Euclidean factor $\reals^n$).  The discrete group $\Gamma$ acts by Riemannian automorphisms of
$D$.   We assume that this action is ``arithmetic'' (see \S \ref{arithmetic}).  

\subsubsection{Correspondences.}  We are interested in the automorphisms of $X$.  A {\em
morphism} $f: X\to Y$ of locally symmetric spaces is a local isometry; i.e. a map $f$ that
restricts to an isometry $U_p \to U_{f(p)}$ for appropriate choices of neighborhoods.
Topologically, a morphism is a covering map of finite degree.  There are finitely many
morphisms $X \to X.$  Instead, we consider a correspondence on $X$, i.e.  a locally symmetric
space $C$ together with two morphisms $c_1$ and $c_2$ from $C$ to $X$.  We can think of
$(c_1,c_2):C\rightrightarrows X$ as a multi-valued function, whose values at $p\in X$ are the
points in the set $c_2(c_1^{-1}(p))$. There is a rich supply of correspondences.  They include
the  {\em Hecke correspondence}, see \S\ref{Hecke correspondences}.

\subsubsection{Lefschetz numbers.} Consider a cohomology theory of $X$, such as the $L^2$
cohmology $H^i_{(2)}(X)$.  A correspondence $(c_1,c_2):C\rightrightarrows X$ acts on
$H^i_{(2)}(X)$ by sending a differential form $\omega$ to $C^* \omega = (c_1)_* (c_2)^*
\omega$. (The map $(c_1)_*$ adds the differential form over the sheets of the
finite covering map $c_1$).  It is thought that the induced maps $C^*: H^i_{(2)}(X)\to
H^i_{(2)}(X)$ carry deep number theoretic significance. One would like to compute them.
Unfortunately, this is too difficult.  Even the rank of $H^i_{(2)}(X)$ is too hard to compute
in most cases.   As often happens, however, computing the alternating sum over the index $i$
is a tractable problem.  Define the Lefschetz number of $C$ by

\[L(C)=\Sum_i (-1)^i \text{Tr } \left( C^*: H^i_{(2)}(X) \to H^i_{(2)}(X)\right) \]
Our goal is to compute this, using the Lefschetz fixed point theorem to find an expression
\[L(C)=\Sum_F L(F)\]
as a sum over fixed point components $F$ of some local contribution $L(F)$.

\subsubsection{Compactifying $X$.}  The first obstacle is the fact that $X$ is, in general,
non-compact.  (When $X$ is compact, the Lefschetz formula was described by M. Kuga and J. H.
Sampson \cite{Kuga}.)  There is no hope for a Lefschetz fixed point formula on a noncompact
space.  To see this, consider the example where $X$ is the real line with the self map that
sends $x\in X$ to $x+1$.  The Lefschetz number for ordinary cohomology is $1$.  But there are
no fixed points, so the right hand side is $0$ no matter how $L(F)$ is interpreted.  There are
similar examples for $L^2$ cohomology and locally symmetric spaces $X$.  

% reference to section on examples? Put in the Klein icosahedral sphere?

The solution is to pass to a compactification $\overline{X}$ of $X$.  We need a
compactification
$X\subseteq \overline{X}$ satisfying

\begin{enumerate}
\item The $L^2$ cohomology of $X$ can be expressed locally on $\overline{ X}.$
\item The correspondence $(c_1,c_2):C\rightrightarrows X$ extends to a compactified
correspondence $(c_1,c_2):\overline{ C}\rightrightarrows \overline{ X}$.
\item The singularities of $\overline{ X}$ are simple enough to allow a calculation of $L(F)$.
\end{enumerate}

\subsubsection*{Remarks on these properties} In (1),  ``Expressed locally'' means that the
$L^2$ cohomology is the cohomology of a complex of sheaves on $\overline{X}$.   For (2), we
want a functorially constructed compactification $\overline{C}$ of $C$.  The map $C^*$ on
$L^2$ cohomology should be induced by a lift of the complex of sheaves to $\overline{C}$.
Properties (1) and (2) together will imply that there exists an expression $L(C)=\Sum L(F)$
for $L(C)$ as a sum of contributions $L(F)$ over fixed point components $F$ of $\overline{C}$,
by applying the fixed point formula of Grothendick and Illusie \cite{GI}. 

A lot of effort has gone in to constructing various compactifications of $X$.  Most
of these fail some of the criteria, however.  For example, the toroidal compactification
\cite{AMRT} for Hermitian $X$ satisfies (3) but neither (1) nor (2).  The Borel-Serre
compactification \cite{BS} satisfies (2) and (3) but not (1).  (It does satisfy (1) in the
case of ordinary cohomology rather than $L^2$ cohomology.  In this case, U. Weselmann,
following suggestions of G. Harder, has carried out the Lefschetz computations (\cite{We}).) 

It is likely that (for sufficiently high rank) any compactification satisfying (1) and (2)
must be singular.  A well known example is the Baily-Borel compactification for Hermitian $X$.
This satisfies (1) because of the Zucker conjecture (i.e. the Looijenga, Saper-Stern theorem)
which expresses the $L^2$ cohomology of $X$ on the Baily-Borel compactification as the
intersection cohomology, which is the cohmology of a complex of sheaves (see \cite{IHII}.)  It
satisfies (2) because it is functorial.  However, the singularities of the Baily-Borel
compactification are as complex as a locally symmetric space only slightly smaller than $X$
and are too complicated to allow a direct computation of $L(F)$.  

The first miracle is that there is a compactification satisfying all three properties:  the
reductive Borel-Serre compactification $\overline{X}$ (defined for all $X$, Hermitian or not).
In the Hermitian case, it may be thought of as a (non algebraic) partial resolution of
singularities of the Baily-Borel compactification. The reductive Borel-Serre compactification
satisfies property (2) because it is functorial with respect to morphisms of locally symmetric
spaces.  So in the diagram $(c_1,c_2):\overline{C}\rightrightarrows\overline{X}$, the space
$\overline{C}$ is the reductive Borel-Serre compactification of $C$.  It satisfies property
(1) because of the existence of weighted cohomology described below, and it satisfies property
(3) because its singularities may be explicitly constructed from certain nilmanifolds (see
\S\ref{nilmanifold}). 

\subsubsection{Weighted Cohomology.}  The Lefschetz fixed point formula of this paper is for
the {\em weighted cohomology groups} $W^\nu H^i(\overline{X}, \mathbf E)$ where $X$ is any
locally symmetric space and $\mathbf E$ is a local system over it.  These were introduced in
\cite{GHM}. The weighted cohomology is the cohomology of a complex of sheaves ${\mathbf
W}^\nu{\mathbf C}^\bullet (\mathbf E)$.  Therefore it satisfies property (1) mentioned above.

The weighted cohomology groups $W^\nu H^i(\overline{X},\mathbf E)$  depend on an auxiliary
parameter $\nu$ called a {\em weight profile}.  When $X$ has finite dimensional $L^2$
cohomology, then $W^\nu H^i(\overline{X},\mathbf E)=H^i_{(2)}(X,\mathbf E)$ provided the
weight profile $\nu$ is chosen to
be the ``middle weight'' (\cite{GHM},\cite{Nair}), so our formula includes the $L^2$ case.
Another weight profile gives the ordinary cohomology of $X$.  A. Nair \cite{Nair} has shown
that for any weight profile, the weighted cohomology $W^\nu H^i(X,\mathbf E)$ is equal to J.
Franke's ``weighted $L^2$ cohomology'' \cite{Franke} for a particular weight function.  For a
leisurely account of the properties of weighted cohomology, see the introduction to
\cite{GHM}.

\subsubsection{The Lefschetz formula.} Even on a compact space with tractable singularities,
the fixed point contribution $L(F)$ is usually too difficult to compute.  The second miracle
is that on the reductive Borel-Serre compactification $\overline{X}$, each correspondence is
{\em hyperbolic}.  A formula for the contributions $L(F)$ for hyperbolic correspondences was
developed in \cite{LFP} expressly for the application in this paper to Hecke correspondences.  

The rest of the introduction proceeds as follows.  The next section enumerates the fixed point
components and determines their topology.  Section \ref{local contribution} describes the
local contribution from each fixed point component, and \S\ref{thm-intro} adds them up to
give the Lefschetz number $L(C)$.

\subsection{The structure of a correspondence on $\overline{X}$.}  The theory of
correspondences on $\overline{X}$ is very self-referential. The reductive Borel-Serre
compactification $\overline{X}$ is a stratified space whose strata are themselves locally
symmetric spaces.  The closure of such a stratum is its reductive Borel-Serre
compactification.   The fixed point components of a correspondence on $X$ are (almost) locally
symmetric.  A correspondence restricted to a stratum of $\overline{X}$ is itself a
correspondence.

\subsubsection{Notations concerning algebraic groups.}\label{arithmetic}  In order to be
precise, we need the language of algebraic groups.  We use boldface symbols for linear
algebraic groups, and Roman symbols for their Lie group of real points, for example,
$G=\mathbf G(\reals)$. 

\subsubsection*{ Obtaining $X$ from $\mathbf G$} Throughout this paper we
fix a reductive linear algebraic group $\mathbf G$ defined over the rational numbers $\mathbb
Q.$   The {\em symmetric space $D$ for $\mathbf G$} is defined to be $G/S_G K$.  Here $K$ is a
maximal compact subgroup of $G$ and $\mathbf S_{\mathbf G}$ is the greatest $\mathbb Q$-split
torus in the center of $\mathbf G$.  The group $G$ acts on $D$ by Riemannian
automorphisms.  Let $X=\Gamma\backslash D = \Gamma \backslash G /S_G K$, where
$\Gamma\subset\mathbf G(\mathbb Q)\subset G$ is assumed to be a neat arithmetic subgroup.
  (This is the arithmeticity assumption on the locally symmetric space $X$
of \S\ref{Locally symmetric spaces}. By results of Margulis, in most cases the arithmeticity
assumption is automatically satisfied.) (The space $D$ may have Euclidean factors because $G$
may have a part of its center that is split over $\reals$ but not over $\mathbb Q$.  After
dividing by $\Gamma$, these Euclidean factors will become wound into tori.)

\subsubsection*{ Rational parabolic subgroups $\mathbf P$ of $\mathbf G$}  If $\mathbf P$ is a
rational parabolic subgroup, let $\mathbf {L_P}$ be its Levi quotient; let $\nu_\mathbf P :
\mathbf P \to \mathbf{L_P}$ be the projection; let $\mathcal U_{\mathbf P}$ be the unipotent
radical; let $\mathbf S_{\mathbf P}\subseteq \mathbf L_\mathbf P$ be the maximal $\mathbb
Q$-split torus in the center of $\mathbf{L_P};$ and let $\Delta_P$ be the set of
simple roots occuring in $\mathfrak N_P = \text{Lie }(\mathcal U_P).$  Let $K_P = K \cap L_P$
be the maximal compact subgroup of $L_P$ which corresponds to $K;$ set $\Gamma_P = \Gamma
\intersect P$; and $\Gamma_L = \nu_P(\Gamma_P)$.  

\subsubsection{The reductive Borel-Serre compactification $\overline{X}$.} (\S
\ref{subsec-RBS})
The strata of $\overline{X}$ are indexed by $\Gamma$-congugacy classes of rational parabolic
subgroups ${\mathbf P}$ of $\mathbf G$. The stratum $X_P\subseteq \overline{X}$ corresponding
to the conjugacy class containing a parabolic ${\mathbf P}$ is the locally symmetric space
$\Gamma_L\backslash D_P$, where $D_P = L_P/S_P K_P$ is the symmetric space of the
Levi $\mathbf L_\mathbf P$.  If $\mathbf P \subseteq {\mathbf P}'$ then the stratum $X_P$ is
contained in the closure of $X_{P'}.$

\subsubsection{Hecke correspondences.}\label{Hecke correspondences} Let $g \in \mathbf
G(\mathbb Q).$  Let $\Gamma' \subset \Gamma \cap g^{-1} \Gamma g$ be a subgroup of finite
index.  This data determines a correspondence $(c_1,c_2):C[g, \Gamma'] \rightrightarrows X$
as follows.  Let $C[g,\Gamma'] = \Gamma' \backslash D.$  The mapping $c_1$ is obtained by
factoring the projection $d_1:D \to \Gamma \backslash D = X$ through $C[g,\Gamma']$ which may
be done since $\Gamma' \subset \Gamma.$  The mapping $c_2$ is obtained by factoring the
projection $d_2:D \to \Gamma \backslash D$ through $C[g,\Gamma']$, where $d_2(x) = d_1(gx).$
Such a factoring exists because $\Gamma' \subset g^{-1} \Gamma g.$
 
It is a fact (cf. Proposition \ref{prop-covering correspondence}) that every correspondence
may be obtained in this way.  For the maximal choice $\Gamma' = \Gamma \cap g^{-1} \Gamma g$,
the resulting correspondence is called a {\it Hecke correspondence} and is denoted
$C[g]\rightrightarrows X.$  Up to isomorphism, this correspondence depends only on the double
coset $\Gamma g \Gamma \in \Gamma \backslash \mathbf G(\mathbb Q) /
\Gamma$ (cf. Lemma \ref{lem-isomorphism}).

The correspondence $C[g,\Gamma']\rightrightarrows X$ is a covering of the correspondence
$C[g]\rightrightarrows X$ of degree $d=[\Gamma:  \Gamma'].$  The action of $C[g,\Gamma']$ on
weighted cohomology is simply $d$ times the action of $C[g].$ So, without loss of generality,
we will concentrate on computing the Lefschetz number of the Hecke correspondence $C[g]$ for
a fixed double coset $\Gamma g\Gamma\in \Gamma \backslash \mathbf G(\mathbb Q) / \Gamma$.  

\subsubsection{The correspondence on a stratum of $\overline{X}$.} \label{The correspondence
on a stratum} Each Hecke correspondence $C[g] \rightrightarrows X$ has a unique continuous
extension to the reductive Borel-Serre compactification $\overline{C}[g] \rightrightarrows
\overline{X}$.  Every boundary stratum $C_Q$ of $\overline{C}[g]$ will be a correspondence
taking a boundary stratum of $\overline{X}$ to another one.  Since we are interested in fixed
points, we want to classify those $C_Q$ which take a stratum $X_P \subseteq \overline{X}$ to
itself.  There is one of these for every  double coset $\Gamma_P g_i \Gamma_P$ in the
intersection $P \cap \Gamma g \Gamma$  (cf. Proposition
\ref{prop-restriction}).  It is isomorphic to a correspondence of the form $C[\bar g_i ,
\Gamma'_L]$ as described above, but with $G$ replaced by $L_P.$  (Here  $\bar g_i =
\nu_P(g_i)\in L_P$  and $\Gamma'_L = \nu_P(\Gamma_P \intersect g_i^{-1} \Gamma g_i)\subset
L_P$.)

\subsubsection{Fixed point components} \label{Fixed point components} The fixed point set
(sometimes called the {\it coincidence set}) of a correspondence  $(c_1,c_2):C
\rightrightarrows X$ is by definition the set of points in $C$ on which the two maps $c_1$ and
$c_2$ agree.

The fixed point set of the correspondence $C[g, \Gamma']\rightrightarrows X$ (before
compactification) breaks up into constituent pieces $F(e)$ indexed by $\Gamma$ conjugacy
classes of elliptic (modulo $S_G$) elements $e\in \Gamma g \Gamma$  (\S
\ref{subsec-char-elt}).  The piece $F(e)$ corresponding to the conjugacy class containing $e$
is the space $\Gamma'_e\backslash G_e/K'_e$, where $G_e$ is the centralizer of $e$ in $G$,
$\Gamma'_e=\Gamma'\intersect G_e$, and $K'_e=G_e \intersect (z(S_G K)z^{-1})$ where $zA_G K
\in D = G/A_G K$ is a fixed point of the action of $e$ on $D.$  (Such a point exists since $e$
is elliptic) (\S\ref{prop-char}). The constituent $F(e)$ is a finite union of spaces, each of
which is almost a locally symmetric spaces for the group $G_e.$ (It have infinite volume
because it may have Euclidean factors that are not wound into tori.)

Applying this result to the boundary stratum $X_P$, in conjunction with the calculation
(\S\ref{The correspondence on a stratum}) of the part of the correspondence $\overline{C}[g]$
which preserves $X_P$,  we get a group theoretic enumeration of
all the fixed points of $\overline{C}[g]$ which lie over points in $X_P:$ For each choice
of a double coset $\Gamma_P  g_i \Gamma_P\subset \Gamma g \Gamma \cap P$, and for each
conjugacy class of elliptic (modulo $S_P$) elements $e$ in $\Gamma_L \bar g_i \Gamma_L$, there
is a fixed point constituent $F_P(e).$  (It is a smooth submanifold of the stratum $C_Q$ of
$\overline{C}[g]$ which is determined by the double coset $\Gamma_P g_i \Gamma_P$ as in \S
\ref{The correspondence on a stratum}.)  Summing over $\Gamma$ conjugacy classes of rational
parabolics $P$ gives the complete enumeration of fixed points $\overline{C}[g]$.

\subsubsection{The topology of the fixed point set.}\label{The topology of the fixed point
set} There are finitely many constituents $F_P(e)$ of the fixed point set and they are
disjoint.  Unfortunately however they may not be topologically
isolated from each other.  If $X_{P'} \subseteq \overline{X}_P$, then the closure of $F_P(e)$
can contain points in some $F_{P'}(e')$.  So a single connected component of the fixed point
set may have a very complicated structure.  This phenomenon is the main source of technical
difficulty of this paper.  (The only real limit we have found on the possible complexity of a
connected component of the fixed point set is  Proposition \ref{prop-structure}.)  We get
around this problem by composing the correspondence with a mapping, very close to the
identity, which shrinks a neighborhood of the singularity set $\overline{X} - X$ into the
singularity set, and which does something similar on the closure of each stratum of
$\overline{X}.$  This has the effect of ``truncating'' each connected component of the fixed
point set into pieces each of which is contained in a single stratum of $\overline{X}$ (a
process which may be considered as a sort of topological analog to Arthur's truncation
procedure).  The Lefschetz number of this ``modified'' Hecke correspondence is equal to that
of the original one.  We prove that the modified Hecke correspondence is hyperbolic.  
The resulting formula (Theorem \ref{thm-intro}) would be the same if no truncation were used,
however the proof would be even more technical.   

\subsection{Calculating the local contribution $L(F)$.}\label{local contribution} 

\subsubsection{Weighted cohomology.} Let $E$ be a
finite dimensional representation of $\mathbf G$, and let $\mathbf E$ be the associated local
system over $X$.  Denote by $\mathbf{P_0}$ a fixed minimal (``standard'') rational parabolic
subgroup of $G$ and by $\mathbf{S_0}$ a maximal $\mathbb Q$-split torus in the center of its
Levi factor.  Then $\mathbf S_{\mathbf G}\subseteq \mathbf{S_0}.$
A weight profile $\nu\in \chi^*_\mathbb Q (\mathbf{S_0})$ (\S\ref{subsec-WC}) is a
(quasi-)character of $\mathbf{S_0}$ whose restriction to $\mathbf S_{\mathbf G}$ coincides
with the character by which $\mathbf {S_G}$ acts on $\mathbf E$.   The Hecke correspondence
$\overline{C}[g] \rightrightarrows \overline{X}$ has a canonical lift (\S \ref{subsec-Lef1})
to
the weighted cohomology sheaf $\mathbf W^\nu \mathbf C^\bullet(\mathbf E)$, so it induces a
homomorphism on weighted cohomology whose Lefschetz number 
\begin{equation}\label{eqn-Lefschetznumber}
 L(C[g])=\Sum_{i \geq 0}(-1)^i {\rm Tr}(C[g]; W^\nu H^i(\overline{X};\mathbf E))
\end{equation}
is what we want to compute.  

\subsubsection{Hyperbolicity of the correspondence.}  Let us assume for the moment that the
fixed point set is topologically the disjoint union of the constituent pieces $F_P(e)$.  This
is not always the case, but the formula we obtain is nevertheless always valid, as explained
in \S\ref{The topology of the fixed point set}.

We focus on a single stratum $X_P$ which is preserved by the correspondence, and on a single
stratum of the correspondence $\overline{C}$ corresponding to a single double coset $\Gamma_P
g_i \Gamma_P \subset P \cap \Gamma g \Gamma.$  Within this stratum, we focus on a single
constituent $F_P(e)$ of the fixed point set. Each stratum $X_Q$ which contains $X_P$ in its
closure correspondence to a rational parabolic subgroup  $\mathbf Q$
containing $\mathbf P$, and therefore to a unique subset $I\subset
\Delta_P$.  The empty subset corresponds to $X_P$ itself and the largest subset $\Delta_P
\subseteq \Delta_P$ corresponds to $X$.  Let $a_e$ be the projection of $e$ to the identity
component $A_P$ of $\mathbf{S_P}(\mathbb R).$  The elements of $\Delta_P$ are rational
characters of $\mathbf {S_P}$ so we may define
\[\Delta^+_P(e) = \{ \alpha \in \Delta_P | \alpha (a_e) < 1\}.\]
Let $X_Q$ be the stratum containing $X_P$ which corresponds to the subset
$\Delta^+_P(e)\subseteq \Delta_P$.  The correspondence $C[g]$ is hyperbolic near $F_P(e)$
\ref{subsec-hypercorres}, with ``expanding'' (or ``unstable'') set $X_Q$ (Theorem
\ref{thm-hyperbolic}, \S\ref{subsec-outline}).  In other words, near $F_P(e)$ the Hecke
correspondence is ``expanding'' in those directions normal to $X_P$ which point into $X_Q.$

Let $F' = c_1(F_P(e)) = c_2(F_P(e)) \subset X_P$ and let $L_e$ be the centralizer of $e$ in
$L_P.$  There are diffeomorphisms (Proposition \ref{prop-char}),
\begin{equation*}
F_P(e) \cong \Gamma'_e \backslash L_e /K'_e \text{  and  } F' \cong \Gamma_e \backslash L_e
/K'_e  \end{equation*}
where $K'_e = L_e \cap z(K_PA_P)z^{-1}$ (for appropriately chosen $z \in L_P$), $\Gamma_e =
L_e \cap \Gamma_L$, and $\Gamma'_e = L_e \cap \Gamma'_L.$  The projection $F_P(e) \to F'$ is a
covering of degree $d = [\Gamma_e : \Gamma'_e]$ (cf \S \ref{subsec-addendum}).
 
  It follows from the Lefschetz fixed point theorem of \cite{LFP} that the
local contribution is given by 
\[L(F_P(e)) = \chi_c(F') \Sum_{i\geq 0} (-1)^i Tr\left(C[g]^*: H^i_x(\mathbf A^\bullet) \to
H^i_x(\mathbf A^\bullet)\right)\]
(See Theorem \ref{thm-topological}).  Here $\mathbf A^\bullet = h^! j^* \mathbf W^\nu \mathbf
C^\bullet (E)$ where $h$ is the inclusion $F'\hookrightarrow {X}_P
\hookrightarrow \overline{X}_Q$ and $j$ is the inclusion $\overline{X}_Q \hookrightarrow
\overline{X}$;  $H^i_x(\mathbf A^\bullet)$ denotes its stalk cohomology at $x\in F'$; and
$\chi_c$ denotes the Euler characteristic with compact supports.     See the introduction to
\cite{LFP} for a geometric account of hyperbolic correspondences.

\subsubsection{The stalk cohomology} Let $c$ denote the codimension of $F'$ in $X_P$ and let
$\mathcal O$ be the top exterior power of the normal bundle of $F'$ in $X_P.$  Let $r$ denote
the index $[\Gamma_P \cap \mathcal U_P : \Gamma'_P \cap \mathcal U_P]$ where $\Gamma'_P =
\Gamma_P \cap g_i^{-1} \Gamma_P g_i.$  

The stalk cohomology of the sheaf $\mathbf A^\bullet$ at
a fixed point $x \in F'\subset X_P$ is given by (\ref{eqn-stalk cohomology}) and Proposition
\ref{prop-Kostant}:
\begin{equation*}
{\textstyle  H^\bullet_x(\mathbf A^\bullet) =  \underset{\substack{w\in W^1_P \\
I_{\nu}(w) = \Delta^+_P(e)}}{\bigoplus} V^L_{w(\lambda_B+\rho_B)-\rho_B}[-\ell(w) -
|\Delta_P^+(e)|-c] \otimes \mathcal O_x. }
\end{equation*}
The Hecke correspondence $\overline{C}[g]$ acts on the first factor by $rd$ times the action
of $e^{-1}$ and it acts on the second factor by $(-1)^c,$ cf. (\ref{eqn-dr}).   We now
describe the other symbols in this formula.

Let $\mathbf T$ be a maximal torus (over $\mathbb C$) in $\mathbf G$, and let $\mathbf B$ be a
Borel subgroup (over $\mathbb C$) of $\mathbf G$ containing $\mathbf T.$ These may be chosen
as in \S \ref{subsec-Kostant} so that $\mathbf{S_0}(\mathbb C)\subset \mathbf T(\mathbb C)$
and so that  $\mathbf B \subset \mathbf{P_0}.$  Let $W_G= W(\mathbf G(\mathbb C),\mathbf
T(\mathbb C))$ denote the Weyl group of $\mathbf G.$  The choice of $\mathbf B$ determines
positive roots $\Phi^+_G=\Phi^+(\mathbf G(\mathbb C), \mathbf T(\mathbb C))$, and a length
function $\ell$ on $W_G.$  Let $W^1_P\subset W_G$ denote the set of Kostant representatives:
the unique elements of minimal length from each of the cosets $W_Px \in W_P \backslash W_G,$
where $W_P = W(\mathbf{L_P}(\mathbb C), \mathbf T(\mathbb C))$, (\S \ref{subsec-Kostant}).
The last sum is over those $v\in W^1_P$ such that the set 
\begin{equation*}
I_{\nu}(w) = \{\alpha \in \Delta_P|\ 
\langle (w(\lambda_B+\rho_B)-\rho_B -\nu, t_{\alpha}\rangle < 0
\}
\end{equation*}
coincides with the set $\Delta_P^+(e)$ defined above (after conjugating $\mathbf P$ so as to
contain $\mathbf B$). Here, as in \S \ref{subsec-Kostant}, $\lambda_B$ denotes the highest
weight of the representation $E$, and $\{ t_{\alpha} \}$ form the basis of the cocharacter
group $\chi_*^{\mathbb Q}(\mathbf{S_P}/\mathbf{S_G})$ which is dual to the basis $\Delta_P$ of
the simple roots.  Also, $\rho_B$ denotes the half-sum of the positive roots $\Phi^+_G.$  The
product $\langle (v(\lambda_B+\rho_B)-\rho_B -\nu, t_{\alpha}\rangle$ makes sense:  the
restriction $(v(\lambda_B + \rho_B)-\rho_B-\nu) | \mathbf{S_P}$ is trivial on $\mathbf{S_G}$
and hence defines an element of $\chi^*(\mathbf{S'_P}) \otimes \mathbb Q$ which can then be
paired with $t_{\alpha}.$ For any $\mathbf B$-dominant weight $\beta$, the symbol
$V^L_{\beta}$ denotes the irreducible $\mathbf{L_P}$-module with highest weight $\beta \in
\chi^*(\mathbf T(\mathbb C))$ and $V^L_{\beta}[- m]$ means that the module
$V^L_{\beta}$ is placed in degree $m$.  

The geometry behind this formula is roughly this: Consider the intersection of a small
neighborhood of $x$ in $\overline{X}$ with the largest stratum $X$.  This intersection will
deformation retract to the nilmanifold \label{nilmanifold} $(\Gamma \intersect \mathcal
U_P)\backslash \mathcal U_P$.  The cohomology of this intersection with coefficients in
$\mathbf E$ coincides with the $\mathfrak N_P$ cohomology by Van Est's theorm, which is
computed by Kostant's theorem to be $\bigoplus_{w\in W^1_P} V^L_{w(\lambda_B+\rho_B)-\rho_B}
[-\ell(w)]$.  The cut-off ($I_{\nu}(w) = \Delta_P^+(e)$) and the degree shift (by $\ell(w) + |
\Delta_P^+(e)|$) come from the computation of $h^! j^* \mathbf W^\nu \mathbf C^\bullet (E)$ in
\S \ref{sec-localweighted}.  The integer $r$ is the ramification index:  the degree of the
mapping $c_1$ when it is restricted to this nilmanifold ( \S \ref{pf-vanEst}). 

By adding the contributions $L(F_P(e))$ over all the fixed point constituents
$F_P(e)$ we arrive at the main result in this paper.

\begin{thm}\label{thm-intro} Let $g \in \mathbf G(\mathbb Q).$  Let $\overline{C}[g]
\rightrightarrows \overline{X}$ be the resulting Hecke correspondence.  Fix a weight profile
$\nu \in \chi^*_{\mathbb Q}(\mathbf{S_0}).$  The Lefschetz number $L(C[g])$
(\ref{eqn-Lefschetznumber}) is given by
\begin{equation*}
\textstyle
\underset{\{\mathbf P\}}{\sum}\ \underset{i}{\sum} \ \underset{\{e\}}{\sum}\ r
\chi_c(\Gamma'_e\backslash L_e/K'_e)
(-1)^{|\Delta_P^+(e)|}  \underset{\substack{w\in W^1_P \\
I_{\nu}(w) = \Delta^+_P(e)}}{\sum} (-1)^{\ell(w)}
\text{Tr}(e^{-1};V^L_{w(\lambda_B+\rho_B)-\rho_B})
\end{equation*}
\end{thm}
\noindent
The first sum is over a choice of representative $\mathbf P$, one from each $\Gamma$-conjugacy
class of rational parabolic subgroups of $\mathbf G$.  For such a $\mathbf P$, set $\Gamma g
\Gamma \cap \mathbf P(\mathbb Q) = \coprod_i \Gamma_P g_i \Gamma_P$ (where $\Gamma_P = \Gamma
\cap \mathbf P(\mathbb Q)$ and where $g_i \in \mathbf P(\mathbb Q).$  The second sum is over
these finitely many double cosets.  Set $\bar g_i = \nu_P(g_i) \in L_P$ and $\Gamma_L =
\nu_P(\Gamma_P).$  The third sum is over a choice of representatives $e$, one
from each $\Gamma_L$-conjugacy class of elliptic (modulo $S_P$) elements $e \in \Gamma_L \bar
g_i \Gamma_L$.  The rest of the notations are explained above.

\subsection{Adelic formulation}  One of the main goals of the series of papers
\cite{Montreal}, \cite{LFP}, \cite{GHM}, \cite{GKM}, and the present paper is Theorem 7.14.B
(p.~535) of \cite{GKM}, an expression for the Lefschetz number $L(C[g])$ in the adelic
setting.  If the weight profile $\nu$ is the ``middle'' weight (and if the rank of $\mathbf G$
equals the rank of $\mathbf K$) then the weighted cohomology coincides with the $L^2$
cohomology, and this formula coincides with Arthur's formula \cite{Ar1} (Theorem 6.1).  If the
weight profile $\nu = -\infty$ then the weighted cohomology coincides with the ordinary ``full
cohomology'' $H^*(X,E)$ and this formula coincides with Franke's formula \cite{Franke}
(thm.~21 p.~273).  The paper \cite{GKM} uses the above Theorem \ref{thm-intro} as its starting
point, then modifies it using three main steps.

(1)  The quantity $r \chi_c(\Gamma'_e \backslash L_e/K'_e)$ which appears in Theorem
\ref{thm-intro}, and the sum $\sum_i$ over double cosets $\Gamma_P g_i \Gamma_P \subset
\Gamma g \Gamma \cap P$ (which preceeds it) are replaced by an orbital integral.

(2)  If $L_P/A_P$ does not contain a compact maximal torus, then the stratum $C_P$ makes no
contribution to Arthur's formula or to Franke's formula.  The same holds for the general
formula in Theorem 7.14.B of \cite{GKM}.  However fixed points in such a stratum may make a
nonzero contribution to the formula in Theorem \ref{thm-intro} above.  In \cite{GKM} \S 7.14
the method of descent is used to re-attribute such a contribution to smaller strata $C_Q$ for
which $L_Q/A_Q$ does admit a compact maximal torus.  See also \S \ref{subsec-Euler} of this
paper.

(3)  Theorem \ref{thm-intro} above involves a sum over parabolic subgroups, while Theorem
7.14.B of \cite{GKM} involves a sum over Levi subgroups.  This is achieved in \cite{GKM} 
(p.~529) by grouping together the contributions from those parabolic subgroups with a given
Levi factor.  (This has the remarkable effect of grouping together fixed points with different
contracting-expanding behavior.)  In \cite{GKM} it is shown that the resulting contribution
from a single Levi subgroup may be interpreted in terms of the (Harish-Chandra) character of a
certain admissible representation. In the case of the
middle weight, this fact gives rise to a combinatorial formula for the stable discrete series
characters, which is the second main result of \cite{GKM}.  (Although this discrete series
character formula was discovered by comparing Arthur's formula to Theorem \ref{thm-intro}, the
statement and proof of the character formula in \cite{GKM} is independent of the part of the
paper dealing with Lefschetz numbers.)

\subsection{Related literature.}  Besides the articles listed above, and an extensive
literature on the co-compact case, we mention several other closely related papers.  The
Lefschetz formula in the rank one case was studied by Moscovici \cite{Moscovici} and
Barbasch-Moscovici \cite{Barbasch}, also by Bewersdorff \cite{Bewersdorff} and Rapoport
\cite{Rapoport}.  In \cite{Stern}, M. Stern gave a general Lefschetz formula for Hecke
correspondences.  We do not easily see how to compare his formula with ours.  In
\cite{Shokranian} S. Shokranian, following the outline in \cite{GKM}, describes a formula for
the Lefschetz numbers of Hecke operators on twisted groups. We wish to draw attention to
Langlands' article \cite{Langlands}, in which the expanding and contracting nature of the
fixed points on the boundary was first isolated (see especially Proposition 7.12 p.~485).

\subsection{Acknowledgements}
We would like to thank J. Arthur, W. Casselman and R. Langlands for encouraging us to 
work on this question.  We would especially like to thank R. Kottwitz for patiently 
explaining Arthur's formula to us and for helping to interpret our early results in this
direction.  Some of the results in this paper appeared earlier in the adelic setting in our
joint paper \cite{GKM} with R. Kottwitz.  We have profited from useful conversations with A.
Nair, A. Borel, W. Casselman, P. Deligne, G. Harder, E. Looijenga, A. Nicas, M. Rapoport, L.
Saper, J. Steenbrink, M. Stern, and S. Zucker.
The first author is grateful to the Institute for Advanced Study for its support while much of
this paper was written.  This research was begun and partially completed when the authors
were at Northeastern University and the Massachusetts Institute of Technology, respectively.
We are also grateful to the following institutions for their
hospitality and support during various phases of this project:  the Centre de Recherches
Math\'ematiques at the Universit\'e de Montr\'eal, the MaxPlanck Institut f\"ur Mathematik in
Bonn, the Department of Mathematics at the University of Chicago, the Universita di Roma la
Sapienza.  This research was partially supported by the National Science Foundation under
grants number DMS-8802638, DMS-9001941, DMS-9303550, DMS-9626616, DMS-9900324 (Goresky) and
DMS-8803083, DMS-9106522,(MacPherson). 

\subsection{List of symbols}
\begin{itemize}
\item[\S\ref{sec-notation}] {\bf(\ref{subsec-LSS}):}
 $\mathbf{G},$ $\mathbf{S_G},$ $A_G,$ $\mathbf{{}^0G},$ $K,$ $D,$
$\mathbf{G^{(1)}},$ $K^{(1)},$ $A_0,$ $T_g,$ $K(x),$ $\psi_x,$ $K',$ elliptic,
$\theta,$ $\Gamma,$ $X,$ $\tau,$
{\bf(\ref{subsec-parabolics}):} $\mathbf{P},$ $\mathcal U_P,$ $R_dP,$ $L_P,$ $\nu_P,$ $M_P,$
$\mathbf{S_P},$ $A_P,$ $\Gamma_P,$ $K_P,$ $\mathbf{S'_P},$ $A'_P,$ $i_{x_0},$
{Langlands' decomposition}, $a_g,$ {geodesic action}, {torus factor},
{\bf (\ref{subsec-roots}):} $\mathbf{P_0},$ $\mathbf{S_0},$ $\Phi,$ $\mathfrak N_0,$ $\Delta,$
$\mathbf{P}_0(I),$ $\chi_Q(\mathbf{S'_P}),$
{\bf (\ref{subsec-twoparab}):} $\mathbf{\overline{P}},$ $i(\Delta_Q),$ 
{complementary decomposition}, {orthogonal decomposition},
{\bf(\ref{subsec-boundary components}):} {boundary component}, {boundary
stratum}, $e_P,$ $Y_P,$ $D_P,$ $X_P,$ $\mu,$ $F_P,$ {canonical cross section},
{\bf(\ref{subsec-connection}):} $\mathcal{M}_x,$
{\bf(\ref{subsec-Borel-Serre}):} $\overline{A}'_P,$ $A'_P(>1),$ $\overline{A}'_P(\ge 1),$
$\widetilde{D},$ $D(P),$ $\widetilde{X},$
{\bf(\ref{subsec-RBS}):} $\pi_P,$ {geodesic projection}, $\overline{D},$ $\mu,$
$D[P],$ $\overline{X},$ $\tau$
\item[\S\ref{sec-parabolicneighborhood}]
{\bf(\ref{subsec-parabnbhd}):} $\alpha,$ $\beta,$ $\Gamma-${parabolic},
{\bf(\ref{subsec-rootfunction}):} {root function}, $f_{\alpha}^P,$ $\pi_P$
\item[\S\ref{sec-tilings}]
{\bf(\ref{subsec-tilings}):} $\mathcal P_1,$ $\mathbf{b},$ {parameter}, $\mathcal B,$ 
{tiling}, $D^P,$ $\partial^PD^0,$ $D_Q^0,$
{\bf(\ref{thm-Saper}):} $T(\overline{D}_Q),$ $\partial T(\overline{D}_Q),$ $r_{\alpha}^Q,$
partial distance function,
{\bf(\ref{subsec-tilingofXbar}):} $X^P,$ $T(\overline{X}_P),$ $\partial T(\overline{X}_P),$
{\bf(\ref{subsec-retraction}):} $R,$ {retraction}, $W,$ {exhaustion}, $R_Q,$
$W_Q$
\item[\S\ref{sec-shrink}]
{\bf(\ref{subsec-5.2}):} $D\{Q\},$ $\rho,$ $Sh(Q,t),$ $Sh_{\mathbf{Q}}(t),$
{\bf(\ref{subsec-5.4}):} $Sh(\mathbf{t})$
\item[\S\ref{sec-Hecke}]
{\bf(\ref{defn-morphism}):}  {morphism}, $\text{Mor}(X',X),$
{\bf(\ref{lem-extension}):} $\widetilde{f},$ $\overline{f},$
{\bf(\ref{def-Hecke}):} {correspondence}, $\Gamma[g],$ $C[g],$
{\bf(\ref{subsec-narrow}):} {narrow}
\item[\S\ref{sec-restriction}]
{\bf(\ref{subsec-restriction}):} {parabolic correspondence}, $\Gamma_P[y],$
{modeled},
{\bf(\ref{prop-restriction}):} $\Xi$
\item[\S\ref{sec-counting}]
{\bf(\ref{subsec-char-elt}):} {fixed point}, {characteristic element}, $e,$
$F_P(e),$ {elliptic}, $L_e$
\item[\S\ref{sec-hyperbolic}]
{\bf(\ref{subsec-hyperbolic}):} $\Delta_P^{+},$ $\Delta_P^{-},$ $\Delta_P^0,$
{\bf(\ref{prop-neutral}):} {neutral}, $\mathbf{P} \prec \mathbf{Q},$
{\bf(\ref{subsec-dagger}):} $P^{\dagger}$
\item[\S\ref{sec-modifiedHecke}] {\bf(\ref{subsec-tangential}):} $E,$ $d_E,$
{\bf(\ref{subsec-hypercorres}):} hyperbolic
\item[\S\ref{sec-localweighted}] {\bf(\ref{subsec-quadrants}):} $t_{\alpha},$ 
$\chi^*_{\mathbb Q}(\mathbf{S_P})\qdr{J},$  $\chi^*_{\mathbb Q}(\mathbf{S_P})_{\ge
\nu_P(J)},$
{\bf(\ref{subsec-WC}):} $\mathbf E,$ $\nu,$ weight profile, $\mathbf{W^{\nu}C^{\bullet}}
(\mathbf E),$
{\bf(\ref{subsec-sheafremarks}):} $\mathcal O_{X/Y},$
{\bf(\ref{subsec-proof}):} $N_y,$ $\mathcal L_y,$ $\delta,$ $\triangleright^{s-1},$
$\triangleright_J,$
{\bf(\ref{subsec-Kostant}):} $\Phi_G^{+},$ $\Phi_L^{+},$ $\rho_B,$ $W_G,$
$W_P,$ $W^1_P,$ $V^L_{\beta},$ $\lambda_B,$ $I_{\nu}(w)$
\item[\S\ref{sec-Lefschetz}]  {\bf(\ref{subsec-Lef1}):} $\mathbf{A}^{\bullet},$
{\bf(\ref{thm-topological}):} $\chi_c,$
{\bf(\ref{subsec-computation}):} $r,$
{\bf(\ref{subsec-nilmanifold}):} $N_P,$
{\bf(\ref{pf-vanEst}):} $\mathbf{C}^{\bullet}(N_P,\mathbf E),$ $C^{\bullet}(\mathfrak N_P,
E),$ $\Omega_{\text{inv}}(\mathcal U_P, E)$
\end{itemize}

\newpage

\section{Notation and Terminology}\label{sec-notation}
\subsection{Locally symmetric spaces}\label{subsec-LSS}  %sect. 2.1

Linear algebraic groups will be represented by boldface symbols (e.g., $\mathbf G$, $\mathbf
S$) and their real points will be in Roman type (e.g., $G = \mathbf G(\mathbb R)$, $S =
\mathbf S(\mathbb R)$).  Throughout this paper we fix a connected reductive linear algebraic
group $\mathbf G$ defined over $\mathbb Q$.  Denote by $\mathbf{S_G}$ the greatest $\mathbb
Q$-split torus in the center of $G$, and let $A_G = \mathbf{S_G}(\mathbb R)^0$ denote the
identity component of the group of real points of $\mathbf{S_G}$.  Following \cite{BS} \S 1.1
let 
\begin{equation*}
\mathbf{{ }^0G} = \bigcap_{\chi} \ker (\chi^2) \end{equation*}
be the intersection of the kernels of the squares of all the rationally defined characters
$\chi: \mathbf G \to \mathbf{GL_1}.$  Then ${}^0G$ is normal in $G$;  it contains every
compact subgroup and every arithmetic subgroup of $G$, and $G = A_G \times {}^0G.$  Let $K
\subset \mathbf G(\mathbb R)$ be a maximal compact subgroup and define $D = G/KA_G.$  We refer
to $D$ as the
 ``symmetric space'' associated to $\mathbf G.$  The derived group $\mathbf{G^{(1)}}$ is
semisimple and $K^{(1)} = G^{(1)} \cap K$ is a maximal compact subgroup.  The space $D$ is
diffeomorphic to the Cartesian product of the Riemannian symmetric space $D^{(1)} =
G^{(1)}/K^{(1)}$ with $A_0/A_G$ where $A_0$ is the identity component of the
greatest $\mathbb R$-split torus in the center of $G$.  Both $\mathbf G(\mathbb R)$ and
${}^0G$ act transitively on $D,$ an action which we usually denote by $(g,x) \mapsto gx$ or
$g. x$ but occasionally it will be necessary to refer to this action as a mapping, in
which case we write
\begin{equation}\label{eqn-action}
T_g:D \to D \end{equation}
for $g \in G.$
(For most geometric questions involving the symmetric space $D$, one could replace $G$ by
${}^0G$, however there are Hecke correspondences for $G$ which do not necessarily come from
${}^0G.$) For each $x\in D$ the stabilizer $K(x)$ of $x$ in ${}^0G$ is a maximal compact
subgroup of ${}^0G$ so we obtain a $G$-equivariant diffeomorphism 
\begin{equation*} \psi_x:G/A_GK(x) \to D. \end{equation*}  

The choice of $K \subset G$ corresponds to a ``standard'' basepoint $x_0 \in D.$  We write $K
= K(x_0)$ and $K' = A_GK(x_0).$  An element $x \in G$ is {\it elliptic mod $A_G$} (often
shortened to ``elliptic'') if it is $\mathbf G(\mathbb R)$-conjuate to an element of $K'.$
There is a unique ``algebraic'' Cartan involution $\theta = \theta_{x_0}:G \to G$ whose fixed
point set is $K$.  If $x_1\in D$ is another basepoint with $x_1 = gx_0$ then the Cartan
involution for the new basepoint is given by 
\begin{equation}\label{eqn-basepoint}
\theta_{x_1}(y)=g\theta_{x_0} (g^{-1}yg)g^{-1}\end{equation} %2.1.1
and the composition $\psi_{x_1}^{-1}\psi_{x_0}:G/A_GK(x_0) \to G/A_GK(x_1)$ is given by
\begin{equation}\label{eqn-basepoint2}
yA_GK(x_0) \mapsto yg^{-1}A_GK(x_1). \end{equation}
Let $\mathfrak g = \mathfrak k \oplus \mathfrak p$ be the $\pm 1$ eigenspace decomposition of
$\theta$ in $\text{Lie}(G).$  The Cartan involution $\theta$ preserves ${}^0G$ and determines
a decomposition of its Lie algebra, ${}^0\mathfrak g = \mathfrak k \oplus {}^0\mathfrak p.$
Then ${}^0\mathfrak p$ may be canonically identified with the tangent space $T_{x_0}D.$  Any
choice of $K$-invariant inner product on ${}^0\mathfrak p$ induces a $G$-invariant Riemannian
metric on $X$.

Throughout this paper we also fix an arithmetic subgroup
 $\Gamma\subset \mathbf G(\mathbb Q)$ and denote by $\tau:D \to X=\Gamma\backslash D$ the
projection to the locally symmetric space $X$.

\subsection{  Parabolic subgroups } % sect. 2.2
\label{subsec-parabolics}

Fix a rationally defined parabolic subgroup $\mathbf P\subset \mathbf G.$  We have
the following groups:
\begin{enumerate}
\item  $\mathcal U_P=$ the unipotent radical of $\mathbf P$; $\mathfrak N_P=\text{
Lie}(\mathcal U_P)$ its Lie algebra
\item  $\mathbf R_d\mathbf P=$ the $\mathbb Q$ split radical of $\mathbf P$
\item  $\mathbf L_{\mathbf P}$ $=$ the Levi quotient; $\nu_P:\mathbf P\to
\mathbf {L_P}$ the projection
\item  $\mathbf M_{\mathbf P}=\ ^0\mathbf L_{\mathbf P}=\cap_{\chi}\ker(\chi^2)$ 
\item  $\mathbf S_{\mathbf P}=\mathbf R_d\mathbf P/\pmb{\mathcal U}_{\mathbf P}$
\item  $A_P=\mathbf S_{\mathbf P}(\mathbb R)^0$ the identity component of the set of real
points
\item  $\Gamma_P=\Gamma\cap P,$ $\Gamma_L=\Gamma_{L(P)}=\nu (\Gamma_P)\subset M_P$
\item  $K_P=K\cap P,$ $K_L=K_{L(P)}=\nu_P(K_P)\subset M_P$, $K'_P = K' \cap P=K_PA_G$
\end{enumerate}
The torus $\mathbf{S_P}$ may also be identified as the greatest $\mathbb Q$-split torus in the
center of $\mathbf{L_P}.$  It contains $\mathbf{S_G}$ and we denote the quotient by $\mathbf
{S'_P} = \mathbf {S_P}/\mathbf{S_G}$, with corresponding identity component $A'_P =
\mathbf{S'_P}(\mathbb R)^0 = A_P/A_G.$  We identify $A'_P$ with the subgroup $A_P \cap {}^0G$
to obtain a canonical decomposition $A_P = A'_PA_G.$

The group of real points of the Levi quotient is the
direct product, $L_P = M_P \times A_P$.  For any
$x\in P$ write $\nu_P(x)=\nu_M(x)\nu_A(x)$ for its $M_P$ and its $A_P$
components and write $\nu_{A'}(x)$ for the further projection of $\nu_A(x)$ to the quotient
$A'_P = A_P/A_G.$  The group $P$ acts transitively on $D$ with isotropy
$K'_P=A_GK_P=\text{Stab}_P(x_0)$.

     The choice of standard basepoint $x_0\in D$ with associated Cartan 
involution $\theta :\mathbf G\to \mathbf G$ determines a unique $\theta 
-$stable lifting [BS] \S 1.9 $i=i_{x_0}:L_P\to P.$
Denote the image by $L_P(x_0)=i(L_P).$ We obtain liftings of subgroups,
$A_P(x_0)=i(A_P)$ and $M_P(x_0)=i(M_P).$ Thus the choice $x_0 \in D$ of basepoint determines
a canonical {\it Langlands' decomposition} \begin{equation}\label{eqn-Langlands}
P = \mathcal U_P A_P(x_0) M_P(x_0) \end{equation}
and we write 
\begin{equation}\label{eqn-Langlandselement}
g = u_ga_gm_g \end{equation}
where $u_g=gi\nu_P(g^{-1})$, $a_g = i\nu_A(g)$, and $m_g = i\nu_M(g)$ for any $g\in P.$  The
groups $K_P \subset P$ and $K_{L(P)} = \nu_P(K_P)$ are canonically isomorphic, in fact,
\begin{equation}\label{eqn-KP}
K_P = i(K_{L(P)}) \subset M_P(x_0) \subset L_P(x_0).\end{equation}  By abuse of notation we
will usually write $K_P \subset L_P.$  If $x_1\in D$ is another basepoint with associated
Cartan involution $\theta_{x_1}:\mathbf G\to \mathbf G$  then, by (\ref{eqn-basepoint}),
the associated $\theta_{x_1}-$stable lifting $i_{x_1}:L_P\to P$ is given by
\begin{equation}\label{eqn-i}i_{x_1}(y)=gi(y)g^{-1} \end{equation} % 2.2.1
where $g\in P$ is any element such that $g\cdot x_0 =x_1 \in D.$
The {\it geodesic action} of Borel and Serre \cite{BS} \S 3  is the right action of
$A_P$ on $D$ which is given by
\begin{equation}
\label{eqn-geodesic}(zK'_P)\cdot a = zi(a)K'_P \in D=P/K'_P \end{equation} % 2.2.2
for any $a \in A_P$ and $z \in P$. It is well defined since $i(A_P)$ commutes
with $K'_P$, and it passes to an action of $A_P' = A_P/A_G.$  The geodesic action commutes
with the (left) action of $P$, and it is independent of the choice of basepoint, by
 (\ref{eqn-basepoint2}).  It is {\it not} an action by isometries.

For $g=u_ga_gm_g\in P$ as in (\ref{eqn-Langlandselement}), the element $a_g \in A_P$ is called
the {\it torus factor} of $g.$    We will often use without mention the following fact:  if
$\gamma=u_{\gamma}a_{\gamma}m_{\gamma} \in \Gamma \cap P$  then $a_{\gamma}=1.$

\subsection{Roots}\label{subsec-roots}
Fix once and for all a  minimal rational parabolic subgroup
$\mathbf {P_0}\subset \mathbf G.$ The parabolic subgroups $\mathbf P \supseteq \mathbf {P_0}$
are called {\it standard}.   Let $\mathbf S_0=i(\mathbf S_{\mathbf P_0})$ be the lift of
$\mathbf S_{\mathbf P_0}.$  Let $\Phi =\ _{\mathbb Q}\Phi (\mathbf
S_0,\mathfrak g)$ denote the rational relative roots of $\mathfrak g$  with respect to
$\mathbf S_ 0$.  The unipotent  radical $\mathcal U_{P_0}$ determines a linear order on the
root system $_{\mathbb Q}\Phi$  such that the positive roots are those occurring in $\mathfrak
N_0=\text{Lie}(\mathcal U_{P_0}).$ Let $\Delta$ denote the resulting collection of simple
roots.  Each subset $I\subset \Delta$ corresponds to a unique standard parabolic subgroup
$\mathbf P_0(I) \supset \mathbf P_0$ (\cite{BS} \S 4, \cite{Borel} \S 14.17, \S 21.11) such
that $\mathbf S_{\mathbf P_0(I)}\subset \ker(\alpha)$ for all $\alpha \in I.$

Suppose $\mathbf P\subset \mathbf G$ is any rational parabolic subgroup.  Then $\mathbf P$ 
is $\mathbf G(\mathbb Q)$-conjugate to a unique standard parabolic subgroup $
\mathbf P_0(I).$  Any choice of conjugating element $\mathbf P = g \mathbf P_0(I) g^{-1}$
gives rise to the same (canonical) isomorphism $\mathbf{S_P} \cong \mathbf{S}_{
\mathbf{P}_0(I)}.$  The elements of $\Delta - I$ give rise, (by
conjugation and restriction to $\mathbf S_{\mathbf P}$) to the set
$\Delta_{\mathbf P}$ of simple roots of $\mathbf S_{\mathbf P}$  occurring in $\mathfrak
N_P.$  The roots $\alpha \in \Delta_P$ are trivial on $\mathbf{S_G}$ and form a basis for the
character module $\chi_{\mathbb Q}(\mathbf{S'_P})=\chi^*(\mathbf
{S_P}/\mathbf{S_G})\otimes_{\mathbb Z} \mathbb Q.$  Rather than follow the common
practice of identifying $\Delta_P$ with $\Delta-I$ we will, for any $\alpha \in \Delta_P$
denote by $\alpha_0\in\Delta$ the unique simple root which agrees with $\alpha$ after
conjugation and restriction to $\mathbf{S_P}.$

\subsection{Two parabolic subgroups}\label{subsec-twoparab}
If $\mathbf P\subset \mathbf Q$  are rational parabolic subgroups then $\overline
{\mathbf P} = \nu_Q(\mathbf P)$ is a rational parabolic subgroup of $\mathbf
L_{\mathbf Q},$ with unipotent radical $\mathcal U_{\overline P}= \mathcal U_P/\mathcal U_Q.$
The $\theta$-stable lifts of the Levi quotients satisfy $L_P(x_0)\subset L_Q(x_0)$
and we have a diagram
\begin{equation*}\label{eqn-diag1}
\begin{CD} \mathcal U_Q\ &\ \subset\ &\ \mathcal U_P\ &\ \subset\ &\ P\ &\ \subset&Q\\
@VVV@VVV@VVV@VV{\nu_Q}V\\
1\ &\ \subset\ &\ \mathcal U_P/\mathcal U_Q\ &\ \subset\ &\ \overline {P}\ &\ 
\subset\ &\ L_Q\\
@VVV@VVV@VV{\nu_{\overline {P}}}V\\
1\ &\ \subset\ &\ 1\ &\ \subset\ &\ L_P \\ \end{CD} \end{equation*} %2.2.3
with $\nu_Q\nu_{\overline{P}} = \nu_P$.
The inclusion $R_d\mathbf Q\subset R_d\mathbf P$ induces an
injection $\mathbf S_Q\hookrightarrow \mathbf S_P$  which agrees with the inclusion $
S_P(x_0)\supset S_Q(x_0)$.  It follows from (\ref{eqn-Langlands2}) below that the geodesic
action of $A'_Q$ on $D$ agrees with the restriction 
(to $A'_Q\subset A'_P)$ of the geodesic action of $A'_P$ on $D$ (cf.  \cite{BS} prop.~3.11).
Each $\alpha \in \Delta_Q$ is the restriction to $\mathbf {S_Q}$ of a unique simple root
$i(\alpha) \in \Delta_P.$  The association $i:\Delta_Q \to \Delta_P$ is injective, so
$\Delta_P$ is the disjoint union
\begin{equation}\label{eqn-decomposition}
\Delta_P = i(\Delta_Q) \amalg J \text{ with }\mathbf{S_Q} = \bigl(\bigcap_{\alpha \in
J}\ker(\alpha)\bigr)^0.
\end{equation}
Among rational parabolic subgroups containing $\mathbf P$, the group $\mathbf Q$ is determined
by the set $J$, and we will write $\mathbf Q = \mathbf P(J).$  The subset $J\subset \Delta_p$
of simple roots may be identified with the set
\begin{equation}\label{eqn-JPbar}
J = \Delta_{\overline{P}}  \end{equation}
of simple roots $\Delta_{\overline{P}}$ of $\mathbf{S'_{\overline{P}}} = \mathbf{S_{\overline
{P}}}/\mathbf{S_Q}$ occurring in $\mathfrak N_{\overline{P}} = \text{Lie}(\mathcal
U_{\overline{P}}).$

A certain amount of confusion arises from the fact that $A'_Q$ has two natural complements
in $A'_P$.  One is the identity component $A'_{Q'}$ of the group of real points of the torus
$\mathbf{S'_{Q'}} = \mathbf{S_{Q'}}/\mathbf{S_G}$ where 
\begin{equation*}
\mathbf{S_{Q'}} = \bigl(\bigcap_{\alpha \in i(\Delta_Q)}\ker(\alpha)\bigr)^0 \subset
\mathbf{S_P}.
\end{equation*}
Then $\mathbf{S_{Q'}}$ is the (identity component of) the center of the Levi quotient of the
largest parabolic subgroup $\mathbf{Q'} \supset \mathbf P$ such that $\mathbf{Q} \cap
\mathbf{Q'}=\mathbf{P},$ which refer to as the  parabolic subgroup containing $\mathbf
{P}$ which is {\it complementary to $\mathbf{Q}.$}  We therefore refer to the {\it
complementary decomposition} $A'_P = A'_QA'_{Q'}.$ The other complement is 
\begin{equation*}
A_P^Q(x_0) = A'_P \cap M_Q(x_0)
\end{equation*}
whose Lie algebra $\mathfrak a_P^Q$ is the orthogonal complement to $\mathfrak a'_Q$ in
$\mathfrak a'_P$ with respect to any Weyl-invariant inner product on $\mathfrak a'_P.$
We will usually identify the quotient $A_P/A_Q = A'_P/A'_Q$ with this second complement,
$A_P^Q,$ and we will refer to the {\it orthogonal decompositions} $A'_P =
A'_Q A_P^Q$ and $A_P = A_QA_P^Q.$

 The canonical Langlands decompositions of $P$ and $Q$ are related as follows: 
Set $\mathcal U_{\overline{P}}(x_0) = i_{x_0}(\mathcal U_P / \mathcal U_Q).$  Note that
$M_P(x_0) \subset M_Q(x_0)$ and $A_Q = A'_QA_G.$  If 
\begin{equation}\label{eqn-Pbar}
\overline{P} = \mathcal U_{\overline{P}}[A_P^Q(x_0)A_G]M_P(x_0)\end{equation}
is the canonical Langlands decomposition of $\overline{P}$, then
\begin{equation}\label{eqn-Langlands2}
\begin{align}
P &= [\mathcal U_Q \mathcal U_{\overline{P}}(x_0)] [A_Q(x_0) A_P^Q(x_0)] M_P(x_0) \\
  &= \mathcal U_Q A_Q(x_0)[\mathcal U_{\overline{P}} (x_0)A_P^Q(x_0)M_P(x_0)].
\label{eqn-Langlands3} \end{align}\end{equation}
The first is the canonical Langlands decomposition of $P$ while the second is the
decomposition of $P$ which is induced from the canonical Langlands decomposition of $Q$.

\subsection{   Boundary strata  }\label{subsec-boundary components} %2.3

Fix a rational parabolic subgroup $\mathbf P \subset \mathbf G.$  Define
\begin{enumerate}
\item the Borel-Serre boundary {\it component} $e_P = D/A'_P$ (quotient under geodesic
action)
\item the Borel-Serre boundary {\it stratum} $Y_P=\Gamma_P \backslash e_P$
\item the reductive Borel-Serre boundary component 
\begin{equation*}D_P=\mathcal U_P\backslash e_P=P/K_PA_P\mathcal U_P
=L_P/K_PA_P=M_P/K_PA_G\end{equation*}
\item the reductive Borel-Serre boundary stratum
$X_P = \Gamma_P \backslash D_P = \Gamma_{L(P)} \backslash D_P.$
\end{enumerate}
The projection $\nu_P:P\to L_P$ induces a projection $\mu:e_P \to D_P$ which passes to a
projection $\mu:Y_P \to X_P.$  Writing $Y_P = \Gamma_P \backslash P/K_PA_P$ and $X_P =
\Gamma_{L(P)} \backslash L_P / K_PA_G$, the mapping $\mu$ is just $\mu(\Gamma_PxK_PA_P) =
\Gamma_{L(P)}\nu(x)K_PA_G.$

 As in \cite{stable} \S 4.2, the Langlands' decomposition (\ref{eqn-Langlands}) determines a
(basepoint-dependent) diffeomorphism,
\begin{equation}F=F_P:\mathcal U_P\times A'_P\times D_P\to D=P/K'_P\end{equation}
by
\begin{equation}\label{eqn-F}
F(u,a,mK_L)=ui(a)i(m)K_P \end{equation} %2.3.2
where $u\in \mathcal U_P$, $a\in A'_P$, and $m \in M_P.$
With respect to the coordinates defined by
the diffeomorphism $F$, the mapping $\mu$ is the projection to the third factor.  The  (left)
action of $g \in P$ and the (right) geodesic action of
$b\in A'_P$ on $D$ are given by
\begin{equation}\label{eqn-actionofP}
g.(u,a,mK_L)\cdot b = (gui\nu(g^{-1}),\nu_A(g).ab,\nu_M(g).mK_L) \end{equation} %2.3.4
(where $u \in \mathcal U_P$, $a,b \in A'_P$, and $m \in M_P$), as may be seen by
applying the function $F$ to both sides of this equation.  For any fixed $b\in A'_P$ the set
$F(\mathcal U_P \times \{b\} \times D_P) \subset D$ is called a {\it canonical cross section};
it is a single orbit of the group
\begin{equation*}
{}^0P = \bigcap_{\chi}\ker(\chi^2) = \mathcal U_P M_P,
\end{equation*}
the intersection being taken over all rationally defined characters $\chi:\mathbf{P} \to
\mathbf{GL_1}.$  The pullback by $F$ of the canonical Riemannian
metric on $D$ is given (\cite{stable} \S 4.3) by  the orthogonal sum,
\begin{equation}\label{eqn-metric}
F^*(ds^2) = \sum_{\beta \in \Phi} a^{-2\beta}h_{\beta}(z) \oplus
da^2 \oplus ds^2_M \end{equation} % 2.3.3
where $ds^2_M$ is the canonical Riemannian metric on $D_P$ as determined by
the Killing form for $M_P$,
where $\Phi$ denotes the set of roots of $\mathcal U_P$ with respect to $A_P$,
and $h_{\beta}(z)$ is a smoothly varying metric on the root space $\mathfrak
u_{\beta}$.

\subsection{  The flat connection }\label{subsec-connection} (\cite{GHM} \S 7.10)\ %2.4
For any point $x=gK'_P\in D$ with $g=uam$ decomposed according to (\ref{eqn-Langlands}),
define the submanifold
\begin{equation}\label{eqn-leaf}
\mathcal M_x=F_P(\{u\} \times \{a\} \times D_P) = u.i(a).i(M_P)K_P\subset D
\end{equation} %2.4.1
\begin{lem}\label{lem-leaf}  The manifold $\mathcal M_x$ is perpendicular to the
fibers of the mapping $\nu_P : D \to D_P.$  The restriction $\nu_P|\mathcal M_x$
is an isometry.  The submanifolds $\mathcal M_x$ form the horizontal submanifolds
of a flat connection on the fiber bundle $\mu:e_P\to D_P,$
which is independent of the choice of basepoint and is invariant under the action of
$\Gamma_P$ and which therefore
passes to a flat connection on $\mu:Y_P \to X_P.$ \end{lem}

\subsection{  Proof } %2.5
Perpendicularity follows from (\ref{eqn-metric}). Also, by (\ref{eqn-metric}), the mapping
$\nu_P$ is an isometry.
Finally the flat connection is $\Gamma_P$-invariant because by (\ref{eqn-actionofP}) the
action of $\gamma\in \Gamma_P$ on $D$ is given by
\begin{equation}\gamma\cdot(u,a,mK_L) = (\gamma u
i\nu(\gamma^{-1}),a,\nu_M(\gamma).mK_L)\end{equation}which does not mix the factors.\qed

\subsection{  Borel-Serre  compactification }\label{subsec-Borel-Serre} %2.6
In this section we recall basic facts from \cite{BS}.
Let $\mathbf P \subset \mathbf G$ be a rational parabolic subgroup.
The elements of $\Delta_P$
determine a canonical isomorphism  (\cite{BS}\S 4.2) $A'_P\cong (\ao,\ainfty)^{\Delta_P}$
which extends to a unique partial compactification,
\begin{equation*}\overline{A}'_P \cong (\ao,\ainfty]^{\Delta_P}\end{equation*}
So each $\alpha \in \Delta_P$ extends to a homomorphism of semigroups $\alpha:\overline{A}'_P
\to (\ao,\ainfty].$  Denote by 
\begin{equation}\label{eqn-Ap1}\begin{align}
{A}'_P(\agt 1) &= \{ a \in {A}'_P |\ \alpha(a) \agt 1 \text{ for all } \alpha \in
\Delta_P \} \\  \label{eqn-Ap2}
{A}'_P(\age 1) &= \{ a \in {A}'_P |\ \alpha(a) \age 1 \text{ for all } \alpha \in
\Delta_P \}\end{align}
\end{equation}
and similarly for 
$\overline{A}'_P(\agt 1)$ and $\overline{A}'_P(\age 1).$ 
The Borel-Serre partial compactification $\widetilde{D}$ of $D$ is obtained by
 adjoining, for each rational parabolic subgroup $\mathbf P\subset\mathbf G$
the rational boundary component $e_P=D/A'_P$ as the set of limits  of the $A'_P$ geodesic
orbits in $D,$ together with the Satake topology \cite{Satake} \S
2 (p.~562), \cite{BS} \S 7.1, \cite{Zucker3} \S 3.7.
It is covered by ``corners''; the corner associated to $P$ is 
\begin{equation}D(P)=D\times_{A'_P}\overline{A}'_P = \coprod_{Q \supseteq P}e_Q.\end{equation}
Then $D(P)$ is an open $\mathbf{P}(\mathbb Q)$-invariant neighborhood (in $\widetilde{D}$) of
the boundary component $e_P,$ on which $\mathbf{P}(\mathbb Q)$ acts in a continuous and
component-preserving way.  The diffeomorphism $F$ of equation (\ref{eqn-F}) extends to a
diffeomorphism of manifolds with corners,
\begin{equation}\label{eqn-Fbar}
\overline{ F}:\mathcal U \times \overline{A}_P \times D_P \cong D(P)\end{equation} % 2.6.1

The action $T_g:D \to D$ of any $g\in \mathbf G(\mathbb Q)$ extends continuously to a mapping
\begin{equation}\label{eqn-Qaction} \widetilde{T}_g:\widetilde{D} \to \widetilde{D}
\end{equation}
which takes the neighborhood $D(P)$ of $e_P$ isomorphically to the neighborhood $D({}^gP)$ of
$e_{{}^gP}$ (where ${}^gP = gPg^{-1}$).  (The proof of this fact is recalled in \S
\ref{lem-extension}, \S \ref{pf-extension}.)  It follows that the {\it Borel-Serre
compactification}  $\widetilde{X}= \Gamma \backslash \widetilde{D}$ is a (compact) manifold
with
corners, stratified with one stratum $Y_P = \Gamma_P \backslash e_P$ for each
$\Gamma$-conjugacy class of rational parabolic subgroups $\mathbf P.$  
The real analytic structure on $D$ extends to a semi-analytic structure on $\widetilde{D}$ and
passes to a subanalytic structure on $\widetilde{X}.$  Denote by
$\tilde\tau:\widetilde{D} \to \widetilde{X}$ the natural projection.

\subsection{Reductive Borel-Serre compactification}\label{subsec-RBS}
The reductive Borel-Serre partial compactification $\overline{D}$ of $D$ was first described
in \cite{Zucker1} \S 4.2 p.~190; see also \cite{GHM} \S 8.  It is the topological space 
obtained by collapsing each boundary component $e_P$ in $\widetilde{ D}$ to its reductive
quotient $D_P,$ (\S \ref{subsec-boundary components}) together with the quotient topology. 
(See also \cite{Zucker3} \S 3.7.)  The {\it geodesic projection}
\begin{equation}\label{eqn-piP} \pi_P:D \to D_P \end{equation}
is the composition $D \to e_P \to D_P.$  The closure $\overline{D}_P$ of $D_P$ in
$\overline{D}$ is the reductive Borel-Serre partial
compactification of $D_P.$  Let $\mu:\widetilde{ D }\to \overline{D}$ denote the quotient
mapping: it is continuous, its restriction to $D$ is the identity, and its restriction to each
boundary stratum agrees with the projection $\mu:Y_P \to X_P$ of \S\ref{subsec-boundary
components}.   Define 
\begin{equation}\label{eqn-D[P]}
D[P] = \mu(D(P)) = \bigcup_{Q \supseteq P} D_Q
\end{equation}
to be the image of the corner associated to $P$:  it is an open $\mathbf{P}(\mathbb
Q)$-invariant neighborhood of $D_P$ in $\overline{D}$ on which $\mathbf{P}(\mathbb Q)$ acts in
a component-preserving way.  The action $T_g:\widetilde{ D} \to \widetilde{ D}$
(\ref{eqn-Qaction}) of any $g\in \mathbf G(\mathbb Q)$
passes to a mapping $\overline{T}_g:\overline{D} \to \overline{D}$ which takes the
neighborhood $D[P]$ of $D_P$ isomorphically to the neighborhood $D[{}^gP]$ of $D_{{}^gP}.$
It follows that the {\it reductive Borel-Serre compactification} 
\begin{equation*} \overline{X} = \Gamma \backslash \overline{D}.\end{equation*}
is a compact singular space, canonically stratified with one boundary stratum $X_P=\Gamma_P
\backslash D_P$ for each
$\Gamma$-conjugacy class of rational parabolic subgroups $\mathbf P \subset \mathbf G.$
The closure $\overline{X}_P$ of $X_P$ in $\overline{X}$ is the reductive Borel-Serre
compactification of $X_P.$  There are $| \Delta_P |$ maximal boundary strata $X_Q$ such that
$\overline{X}_Q \supset \overline{X}_P$, each corresponding to a maximal (rational) parabolic
subgroup $\mathbf Q = \mathbf P(\Delta_P - \{\alpha\})$ for $\alpha \in \Delta_P$ (cf. \S
\ref{subsec-twoparab}).  Then $\overline{X}_P$ is the intersection
\begin{equation}\label{eqn-maximalintersection}
\overline{X}_P = \cap_Q \overline{X}_Q \end{equation}
of these $| \Delta_P |$ maximal boundary strata.

It is not difficult to see (\cite{Borel} \S 11.7 (iii)) that if
$\mathbf P$ and $\mathbf P'$ are $\mathbf G(\mathbb Q)$-conjugate but are not
$\Gamma$-conjugate, then 
\begin{equation}\label{eqn-nointersect}
\overline{X}_P \cap \overline{X}_{P'} = \phi.
\end{equation} 

The identity mapping $X\to X$ extends uniquely to a continuous surjection $\mu:\widetilde{X}
\to \overline{X}$ and  the subanalytic structure on $\widetilde{ X}$ passes to a subanalytic
structure on $\overline{X}.$
Denote by $\bar\tau:\overline{D} \to \overline{X}$ the projection.  Define $X(P) =
\tilde{\tau}(D(P))$ and $X[P] = \bar\tau(D[P]).$ 
The following diagram may be useful in helping to sort out these spaces,
\begin{equation}\begin{CD}
e_P\  &\subset\  & D(P)\  &\subset\  & \widetilde{D}\  @>{\tilde\tau}>> \widetilde{X}\
&\supset\ & X(P) &\supset\  & Y_P \\
@VVV @VVV @VV{\mu}V @V{\mu}VV @VVV @VV{\mu}V \\
D_P\  &\subset\  & D[P]\  &\subset\  & \overline{D}\  @>>{\bar\tau}> \overline{X}\
&\supset\ &X[P]  &\supset\
& X_P\ 
\end{CD}\end{equation} 

\section{Parabolic neighborhoods and root functions}\label{sec-parabolicneighborhood}
As in \S \ref{sec-notation}, $\mathbf G$ denotes a connected linear reductive algebraic group
defined over $\mathbb Q$, $D$ denotes the associated symmetric space,  $K' = A_GK(x_0)$ is the
stabilizer in $G$ of a fixed basepoint $x_0\in D$, $\Gamma \subset \mathbf G(\mathbb Q)$
is an arithmetic group and $X = \Gamma \backslash D.$  Although the constructions in this
section refer to the reductive Borel-Serre compactification $\overline{X}$ of $X$ (and the
reductive Borel-Serre partial compactification $\overline{D}$ of $D$), they may just as well
be applied to the Borel-Serre compactification $\widetilde{ X}$ (and the Borel-Serre partial
compactification $\widetilde{ D}$ of $D$).  Rather than repeat each statement for both
compactifications, we will present the RBS case only.

\subsection{Parabolic neighborhoods}\label{subsec-parabnbhd}
  Let $\mathbf P \subset \mathbf G$ be a rational parabolic subgroup.  Let
$\alpha:\overline{D} \to \Gamma_P \backslash \overline{D}$ and $\beta: \Gamma_P \backslash
\overline{D} \to \Gamma \backslash \overline{D} = \overline{X}$ be the projections.  We say
that an open set $V \subset \overline{D}$ is {\it $\Gamma$-parabolic} (with respect to
$\mathbf P$) if 
\begin{enumerate}
\item it is invariant under the geodesic action of the semigroup $A_P(\age 1)$ (\ref{eqn-Ap2})
and
\item if $\gamma \in \Gamma$ and $\gamma V \cap V \ne \phi$ then $\gamma \in \Gamma \cap P.$
\end{enumerate}    
Item (2) means that the covering $\beta:\Gamma_P \backslash \overline{D} \to
\overline{X}$ is one to one on the set $\alpha(V)$ so it takes $\alpha(V)$ homeomorphically to
its image $\bar\tau(V) \subset \overline{X}.$  In this case we will also refer to $\alpha(V)
\subset \Gamma_P \backslash \overline{D}$ (resp. $\bar\tau(V) \subset \overline{X}$) as {\it
$\Gamma$-parabolic} open sets.  
\quash{
\begin{equation*}\begin{CD}
\overline{D} @>{\alpha}>> \Gamma_P \backslash \overline{D} @>{\beta}>> \overline{X}\\
\bigcup && \bigcup && \bigcup \\
V @>>> \alpha(V) @>{\cong}>> \bar\tau(V) \\
\bigcup && \bigcup && \bigcup \\
D_P @>>> \alpha(D_P) @>{\cong}>> X_P
\end{CD}\end{equation*}
} %% end quash
\begin{equation*} \begin{CD}
\overline{D} &\ \supset\ & V & \ \supset\ & D_P \\
@V{\alpha}VV @VVV @VVV \\
\Gamma _P \backslash \overline{D} &\ \supset \ &  \alpha(V) &\ \supset \ & \Gamma_P \backslash
D_P \\
@V{\beta}VV @VV{\cong}V @VV{\cong}V \\
\overline{X} &\ \supset & \  \bar\tau(V) & \ \supset \ & X_P
\end{CD} \end{equation*}
Every stratum $X_P$ admits a fundamental system of $\Gamma$-parabolic neighborhoods.  In
section \ref{subsec-tilings} we will review a theorem of Saper \cite{Saper} (thm. 8.1) which
states the stronger fact that {\it the closure $\overline{X}_P$ of each stratum $X_P \subset
\overline{X}$ admits a fundamental system of $\Gamma$-parabolic
neighborhoods.}

\subsection{Root functions}\label{subsec-rootfunction}
  Let $\mathbf P \subset \mathbf G$ be a rational parabolic
subgroup.  Each character $\alpha \in \chi^*_{\mathbb Q}(\mathbf{S'_P})$ determines a mapping
\begin{equation}\label{eqn-rootfunction}
f_{\alpha}^P:D \to \mathbb R_{>0} \end{equation}
by $f_{\alpha}^P(F(u,a,mK_L))=\alpha(a)$ using (\ref{eqn-F}).  The mapping $f_{\alpha}^P$ is
independent of the choice of basepoint.  For any $g' = u'a'm'\in P,$  any $b'\in A'_P,$ and
any $\gamma \in \Gamma \cap P$ we have
\begin{equation}\label{eqn-rootfunctionP}
f_{\alpha}^P(\gamma g'x\cdot b') = \alpha(a'b')f_{\alpha}^P(x).
\end{equation}
If $\alpha \in \Delta_P$ is a simple root, we say $f_{\alpha}^P$ is a {\it root function}. 
If $\gamma \in \Gamma$, $\mathbf P' = \gamma \mathbf P \gamma^{-1}$ and if $\alpha' \in
\Delta_{P'}$ is the root corresponding to $\alpha \in \Delta_P$ then, for all $x\in D$,
\begin{equation}\label{eqn-roottranslate}
f_{\alpha'}^{P'}(\gamma x) = f_{\alpha}^P(x). \end{equation}
The root function $f_{\alpha}^P:D \to (0,\infty)$
extends to a continuous function $D[P] \to (0,\infty]$ (cf \S \ref{prop-extends} below)
which passes to a function  $\Gamma_P
\backslash D[P] \to (0,\infty]$ whose restriction to any $\Gamma$-parabolic
neighborhood $U \subset \overline{X}$ of $X_P$ we also denote by
\begin{equation*} f_{\alpha}^P:U \to (0,\infty] \end{equation*}
Similarly the geodesic projection $\pi_P:D \to D_P$  (cf (\ref{eqn-piP})) extends
continuously to
a projection $D(P) \to D_P$ and passes to projections $D[P] \to D_P$ and $\Gamma_P \backslash
D[P] \to \Gamma_P \backslash D$ whose restriction to any parabolic neighborhood $U \subset
\overline{X}$ we denote by
\begin{equation}\label{eqn-piP2} \pi_P: U \to X_P \end{equation}
The following lemma is a straightforward consequence of the definitions.

\begin{lem}\label{lem-convergence}  Let $U\subset \overline{X}$ be a parabolic neighborhood of
the stratum $X_P$.  Let $\{x_n\} \subset U$ be a sequence of points and let $y\in X_P$.  The
sequence $\{x_n\}$ converges to $y$ in $\overline{X}$ if and only if the following hold,
\begin{enumerate}
\item $\pi_P(x_n) \to y$ in $X_P$ and
\item $f_{\alpha}^P(x_n) \to \ainfty$ for all $\alpha \in \Delta_P.$  \qed
\end{enumerate} \end{lem}

\subsection{}\label{subsec-JPbar}
Suppose $\mathbf Q \supset \mathbf P$ is another rational parabolic subgroup of $\mathbf G$,
corresponding, say, to a subset $J\subset \Delta_P$ with $\mathbf{S_Q} \subset \ker(\alpha)$
for all $\alpha \in J,$ so that $\Delta_P = i(\Delta_Q) \amalg J$ as in
(\ref{eqn-decomposition}).  Let $\overline{P}=\nu_Q(P)\subset L_Q$ be the resulting
parabolic subgroup of $L_Q$.  It acts transitively on the boundary component $D_Q$. 

Let $x\in D$, say $x = u_Qa'_Qa_Gu_{\overline{P}}a_P^Qm_PK'_P$ is decomposed
according to (\ref{eqn-Langlands3}) with $a_Q=a'_Qa_G.$  Then $\pi_Q(x)=u_{\overline
{P}}a_P^Qa_Gm_PK'_P\in \overline{P}/K'_P = D_Q$ so the following equations hold:
\begin{equation}\begin{alignat}{2}
f_{\alpha}^P(x) &= \alpha(a'_Qa_P^Q) &\quad &\text{for all }\alpha \in \Delta_P
\label{eqn-2.4.1} \\
f_{\beta}^Q(x) &= i(\beta)(a'_Q)       &\quad &\text{for all }\beta \in \Delta_Q
\label{eqn-2.4.2}  \\
f_{\alpha}^{\overline{P}}(\pi_Q(x)) =f_{\alpha}^P(x) &= \alpha(a_P^Q) &\quad &\text{for all }
\alpha
\in J=\Delta_{\overline{P}}  \label{eqn-2.4.4}
\end{alignat}\end{equation}
since $\alpha(a'_Q) = 1$ for all $\alpha \in J.$
From this we may conclude:
\begin{prop}\label{prop-extends}  For all $\alpha \in
\Delta_P$, the root function $f_{\alpha}^P$ extends continuously to a function
$f_{\alpha}^P:D[P] \to (0,\infty]$ such that, for all $x\in D$ we have
\begin{equation}\label{eqn-rootpi}
f_{\alpha}^P(\pi_Q(x)) = \begin{cases} 
\begin{array}{cll} f_{\alpha}^P(x) &\text{for} & \alpha \in J \\
                   \ainfty         &\text{for} & \alpha \in \Delta_P - J.
\end{array} \end{cases} \end{equation}
The boundary component $D_Q \subset D[P]$ is the set of $x\in D[P]$ such that:
\begin{equation}\label{eqn-rootpi2}
\begin{cases}
f_{\alpha}^P(x) = \ainfty &\text{\rm for all }\alpha \in \Delta_P - J \\
f_{\alpha}^P(x) \alt \ainfty &\text{\rm for all }\alpha \in J.  \end{cases} 
\end{equation}
\end{prop}

\subsection{Remarks}\label{subsec-intuition}
Of course similar statements apply to the root function $f_{\alpha}^P: U \to (0,\infty]$ for
any parabolic neighborhood $U \subset \overline{X}$ of $X_P.$  We think of the ``negative
gradient'' of the root functions $f_{\alpha}^P$ as pointing in the ``normal
directions'' to $X_P.$  For $\alpha \in J,$  $- \text{grad } f_{\alpha}^P$ points from $X_P$
``into'' $X_Q.$

\begin{figure}[h!]               %% root functions
\begin{center}\begin{picture}(300,90)(80,-15)  %% (450,220)(0,0)
\setlength{\unitlength}{.09em}   %% {.1em} is the normal size
\linethickness{1pt}
\q(50,50)(75,25)(100,0)         %% left hand diagonal stratum
\q(100,0)(125,0)(175,0)         %% bottom stratum
\q(250,50)(275,25)(300,0)       %% same on left hand side
\q(300,0)(325,0)(375,0)

\linethickness{.1pt}
\q(83.764,16.236)(103.362,35.834)(122.96,55.433)
\q(87.823,12.177)(105.391,26.876)(122.96,41.575)
\q(91.882,8.118)(107.421,17.917)(122.96,27.716)
\q(95.941,4.059)(109.450,8.958)(122.96,13.858)

\q(322.96,0)(322.96,25)(322.96,55.433)
\q(317.22,0)(315.19,22.817)(313.16,45.634)
\q(311.48,0)(307.42,17.917)(303.36,35.834)
\q(305.74,0)(299.651,13.071)(293.562,26.035)

\put(100,0){\circle*{3}} \put(300,0){\circle*{3}}
\put(130,-12){$X_Q$}   \put(62,15){$X_{Q'}$}
\put(95,-12){$X_P$}

\put(330,-12){$X_Q$}   \put(262,15){$X_{Q'}$}
\put(295,-12){$X_P$}

\end{picture}\end{center}
\caption{Level curves of $f_{\alpha}^P$ for $\alpha \in i(\Delta_Q)$ and $\alpha \in J$
respectively.}
\label{fig-collar}\end{figure}
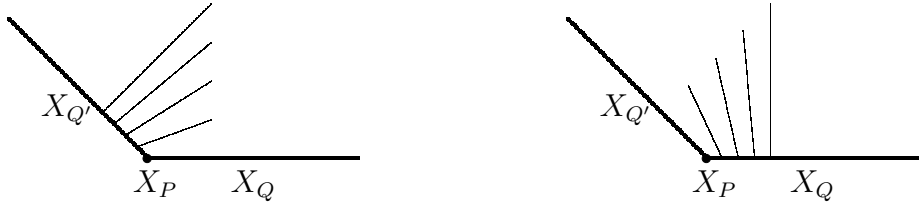

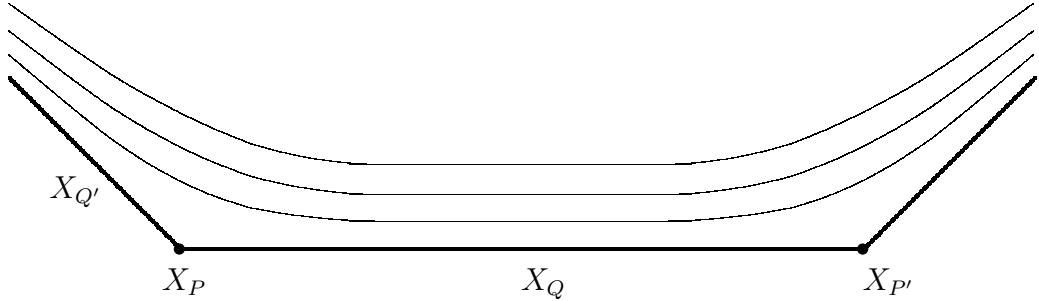
\begin{figure}[h!]               %% Zucker's point
\begin{center}\begin{picture}(450,125)(30,-15)  %% (450,220)(0,0)
\setlength{\unitlength}{.11em}   %% {.1em} is the normal size
\linethickness{1pt}
\q(50,50)(75,25)(100,0)
\q(100,0)(200,0)(300,0)
\q(300,0)(325,25)(350,50)

\linethickness{.1pt}

%% lower level                       middle level
%%\q(200,3.8)(220,4.7)(240,6)
%%\q(240,6)(260,8.2)(280,12.3)
\q(150,8.2)(200,8.2)(240,8.2)   \q(150,16)(200,16)(240,16)		
\q(240,8.2)(260,8.2)(280,12.3)  \q(240,16)(260,16)(280,21.6)	

\q(280,12.3)(285,13.8)(290,15.6)\q(280,21.6)(285,23.5)(290,25.4)	
\q(290,15.6)(295,17.6)(300,20)  \q(290,25.4)(295,27.6)(300,30) 	
\q(300,20)(305,22.6)(310,25.6)  \q(300,30)(305,32.6)(310,35.4)	
\q(310,25.6)(315,28.8)(320,32.3)\q(310,35.4)(315,38.4)(320,41.6)	
\q(320,32.3)(330,40)(350,57)	  \q(320,41.6)(330,48.5)(350,64)	
%\q(320,32.3)(340,48.3)(360,66)
%\q(360,66)(380,84)(400,103)

\q(80,32.3)(70,40)(50,57)	  \q(80,41.6)(70,48.5)(50,64)	
\q(90,25.6)(85,28.8)(80,32.3)   \q(90,35.4)(85,38.4)(80,41.6)	
\q(100,20)(95,22.6)(90,25.6)	  \q(100,30)(95,32.6)(90,35.4)		
\q(110,15.6)(105,17.6)(100,20)  \q(110,25.4)(105,27.6)(100,30)	
\q(120,12.3)(115,13.8)(110,15.6)\q(120,21.6)(115,23.5)(110,25.4)	

\q(160,8.2)(140,8.2)(120,12.3)  \q(160,16)(140,16)(120,21.6)	
\q(250,8.2)(200,8.2)(160,8.2)   \q(250,16)(200,16)(160,16)		
  
%% third level                
\q(150,24.8)(200,24.8)(240,24.8)
\q(240,24.8)(260,24.8)(280,31.2)

\q(280,31.2)(285,33.2)(290,35.3)
\q(290,35.3)(295,37.6)(300,40)
\q(300,40)(305,42.6)(310,45.3)
\q(310,45.3)(315,48.1)(320,51.2)
\q(320,51.2)(330,57.7)(350,72)

\q(80,51.2)(70,57.7)(50,72)
\q(90,45.3)(85,48.1)(80,51.2)
\q(100,40)(95,42.6)(90,45.3)
\q(110,35.3)(105,37.6)(100,40)
\q(120,31.2)(115,33.2)(110,35.3)

\q(160,24.8)(140,24.8)(120,31.2)
\q(250,24.8)(200,24.8)(160,24.8)
%% labels
\put(100,0){\circle*{3}} \put(300,0){\circle*{3}}
\put(200,-12){$X_Q$}   \put(62,15){$X_{Q'}$}
\put(95,-12){$X_P$}\put(300,-12){$X_{P'}$}

\end{picture}\end{center}
\caption{Level sets of $f_{\beta}^Q$}
\label{fig-Zuckerpoint}\end{figure}

 Zucker's vexatious point (\cite{Zucker1} \S 3.19) is that for $\mathbf P
\subset \mathbf Q$ and for $\beta\in\Delta_Q,$ the root functions $f_{\beta}^Q$ and
$f_{i(\beta)}^P$ do not necessarily agree:  see (\ref{eqn-2.4.1}) and (\ref{eqn-2.4.2}) above.
(In fact, they agree precisely when $i(\beta)(a_P^Q)=1$, which is to say, when $A'_{Q'}$ and
$A'_Q$ are orthogonal, where $\mathbf{Q'} \supset \mathbf P$ is the parabolic subgroup
complementary to $\mathbf Q$.) The nature of the level sets of $f_{\beta}^Q$ are
depicted in Figure \ref{fig-Zuckerpoint}.

This shortcoming will be circumvented by replacing the root function $f_{\alpha}^P$ with
Saper's partial distance function $r_{\alpha}^P$ (associated to a tiling), which is patched
together from the various relevant root functions; cf. (\ref{eqn-fixrootfunction}).

\section{Tilings}\label{sec-tilings}
In this section we recall a construction of Saper \cite{Saper}.  An equivalent construction of
Leuzinger \cite{Leuzinger1} could be used instead.  
\subsection{Tilings of $\overline{D}$}\label{subsec-tilings}
As in \S \ref{sec-notation}, $\mathbf G$ denotes a connected linear reductive algebraic group
defined over $\mathbb Q$, $D$ denotes the associated symmetric space,  $K' = A_GK(x_0)$ is the
stabilizer in $G$ of a fixed basepoint $x_0\in D$, $\Gamma \subset \mathbf G(\mathbb Q)$
is an arithmetic group and $X = \Gamma \backslash D.$ 
Let $\mathcal P_1$ denote the set of proper {\it maximal} rational parabolic subgroups of
$\mathbf G.$  For each $\mathbf Q \in \mathcal P_1$ choose $b_Q \in A'_Q.$  The collection
$\mathbf b=\{b_Q\}$ of such choices is called a {\it parameter}, the set of which we denote
by $\mathcal B.$  The parameters are partially
ordered with $\mathbf{b} \le \mathbf{c}$ iff $\alpha_Q(b_Q) \le \alpha_Q(c_Q)$ for all
$\mathbf Q \in \mathcal P_1$, where $\Delta_Q=\{\alpha_Q\}$ is the  simple root
associated with the maximal parabolic subgroup $\mathbf Q.$  A choice $\mathbf b\in\mathcal
B$ of parameter determines, for any rational parabolic subgroup $\mathbf P \subset \mathbf G$
a unique element $b_P \in A'_P$ such that, for each rational maximal parabolic subgroup
$\mathbf Q \supset \mathbf P,$ the element $b_Pb_Q^{-1}$ lies in $A_P^Q$  (cf \S
\ref{subsec-twoparab}).  In other words, $\log(b_Q)$ is the orthogonal projection of
$\log(b_P) \in \mathfrak a'_P$ with respect to any Weyl-invariant inner product on $\mathfrak
a'_P.$

Recall from \cite{Saper} that a {\it tiling with parameter $\mathbf b\in\mathcal B$} is a
cover of the reductive Borel-Serre partial compactification
\begin{equation}\label{eqn-tiling}
\overline{D} = \coprod_{\mathbf P\in \mathcal P}D^P \end{equation}
by disjoint sets (called tiles) such that
\begin{enumerate}
\item The central tile $D^0 = D^G$ is a closed, codimension 0 submanifold with corners
contained in $D$.  Its closed boundary faces $\{ \partial^PD^0 \}$ are indexed by $\mathbf P
\in \mathcal P$ with $\mathbf P \subset \mathbf Q \iff \partial^PD^0 \subset \partial^QD^0.$
\item  Each boundary face $\partial^PD^0$ is contained in the ``cross section''
$F(\mathcal U_P \times \{b_P\} \times D_P)$ where $F$ is defined in (\ref{eqn-F}).
\item Each tile $D^P = \partial^PD^0 \cdot \overline{A}_P(\agt 1)$ is obtained from
$\partial^PD^0$ by flowing out under the geodesic action of the cone $\overline{A}_P(\agt 1)$
(cf (\ref{eqn-Ap1}))\end{enumerate}

For any rational parabolic subgroup $\mathbf Q$, the intersections $\{D^P \cap
\overline{D}_Q\}$ (over all rational parabolic subgroups $\mathbf P \subseteq \mathbf Q$)
form a tiling of the reductive Borel-Serre partial compactification $\overline{D}_Q$, whose
central tile we denote by
\begin{equation}\label{eqn-centraltile}
D_Q^0=D^Q \cap D_Q.\end{equation} Then the tile $D^Q$ is given by
\begin{equation}\label{eqn-singletile}
D^Q = \{ x\in D[Q]\subset\overline{D} |\ \pi_Q(x) \in D_Q^0 \text{ and } f_{\alpha}^Q(x)\agt
\alpha(b_Q) \ \forall \alpha \in \Delta_Q\}
\end{equation}and
\begin{equation}\label{eqn-face}
\partial^QD^0 = \{x \in D|\ \pi_Q(x) \in D_Q^0 \text{ and }f_{\alpha}^Q(x) = \alpha(b_Q)
\forall \alpha \in \Delta_Q \}.
\end{equation}

A tiling, if it exists, is uniquely determined by its parameter $\mathbf b \in \mathcal B$, in
which case we say that the parameter is {\it regular}.  The parameter $\mathbf b$ is {\it
$\Gamma$-invariant} if, for all $\gamma \in \Gamma$, we have $b_{\gamma Q \gamma^{-1}} =
a_{\gamma Q \gamma^{-1}}(\gamma x_0)b_Q.$
The tiling $\{D^Q\}$ is {\it $\Gamma$-invariant} if $\gamma D^P = D^{\gamma P \gamma^{-1}}$
for all $\gamma \in \Gamma$.  A tiling is $\Gamma$-invariant if and only if its parameter
$\mathbf b$ is $\Gamma$-invariant (\cite{Saper} Corollary 2.7).  In (\cite{Saper} Thm. 10.1)
Saper proves the following.

\begin{thm}\label{thm-Saper} If the tiling parameter $\mathbf b\in \mathcal B$ is chosen
sufficiently \alarge \parens{with respect to the above partial ordering} and
$\Gamma$-invariant, then there exists a unique tiling with parameter $\mathbf b\in\mathcal B,$
and it is $\Gamma$-invariant.  Moreover, for any $\mathbf Q\in \mathcal P$ the union
\begin{equation}\label{eqn-tilenbhd}
T(\overline{D}_Q) = \coprod_{P \subseteq Q} D^P \end{equation}
is an open $\Gamma_Q$-invariant parabolic neighborhood of $\overline{D}_Q$ in $\overline{D}$
which may be made arbitrarily small by choosing the parameter $\mathbf b$ sufficiently
large. \end{thm}
Henceforth we shall refer to such a parameter as {\it regular and sufficiently large}.
Denote the closure of $T(\overline{D}_Q)$ by $\overline{T}(\overline{D}_Q),$ and the boundary
by $\partial T(\overline{D}_Q) = \overline{T}(\overline{D}_Q) - T(\overline{D}_Q).$  Following
\cite{Saper} Thm. 8.1 (ii), for each $\alpha \in\Delta_Q,$  define the {\it partial distance
function} $r_{\alpha}^Q:\overline{T}(\overline{D}_Q) \to [0,1)$ by
\begin{equation}\label{eqn-partiald}
r_{\alpha}^Q (x) = \begin{cases}f_{\alpha}^Q(x)^{-1}\alpha(b_Q) & \text{for }x\in
\overline{T}_Q \\
 f_{i(\alpha)}^P(x)^{-1}i(\alpha)(b_P) & \text{for }x\in \overline{T}_P
                   \end{cases}
\end{equation}
(where $\mathbf P \subset \mathbf Q$ and where $i:\Delta_Q \hookrightarrow \Delta_P$ is the
inclusion (\ref{eqn-decomposition})). 

\begin{lem}\label{lem-continuity} The following statements hold. 
\begin{enumerate}
\item  The mapping $r_{\alpha}^Q:T(\overline{D}_Q) \to
[0,1]$ is well-defined, continuous, and piecewise analytic.  
\item For all $\alpha \in \Delta_Q$, the geodesic action by $t\in A_Q(\age 1)$ satisfies
\begin{equation}\label{eqn-geodesic-r}
r_{\alpha}^Q(x\cdot t) = r_{\alpha}^Q(x)\alpha(t)^{-1}
\end{equation}
whenever $x \in \overline{T}(\overline{D}_Q).$  
\item If $x\in \overline{T}(\overline{D}_Q)$ then
\begin{align*}
x\in \overline{D}_Q \iff r_{\alpha}^Q(x)=0 &\ \text{\rm for all } \alpha \in \Delta_Q \\
x\in \partial T(\overline{D}_Q) \iff r_{\alpha}^Q(x)=1 &\ \text{\rm for some }\alpha \in
\Delta_Q.
\end{align*}
\item  If $\gamma \in \Gamma \cap Q$ then $r_{\alpha}^Q(\gamma x) = r_{\alpha}^Q(x).$
\item If $\gamma \in \Gamma$ and $\mathbf Q' = \gamma \mathbf Q \gamma^{-1}$ and if $\alpha'
\in \Delta_{\mathbf{Q'}}$ is the simple root corresponding to $\alpha \in \Delta_{\mathbf{Q}}$
then
\begin{equation*}
r^{Q'}_{\alpha'}(\gamma x) = r^Q_{\alpha}(x).
\end{equation*} 
\item If $\mathbf P \subset \mathbf Q$ and if $\Delta_P = i(\Delta_Q) \amalg J$ as in
\parens{\ref{eqn-decomposition}} then, for all $\alpha \in \Delta_Q$ and for all $x \in
T(\overline{D}_P) \subset T(\overline{D}_Q)$ we have
\begin{equation}\label{eqn-fixrootfunction}
r_{i(\alpha)}^P(x) = r_{\alpha}^Q(x).
\end{equation}
\end{enumerate}\end{lem}

\subsection{Proof}  
As in \S \ref{subsec-twoparab}, write $A'_P=A'_QA_P^Q$ with $\mathfrak a'_Q$ the orthogonal
complement to $\mathfrak a^Q_P$ in $\mathfrak a'_P.$  So the elements $b_P$ and $b_Q$
determined by the parameter $\mathbf b$ satisfy $b_P = b_Qb_P^Q$ for some $b_P^Q \in A_P^Q.$
Now suppose that $x = u_Qa'_Qa_Gu_{\overline{P}}a_P^Qm_PK_P \in P/K_P = D$ is decomposed
according to (\ref{eqn-Langlands3}).
Set $a'_P=a'_Qa_P^Q \in A'_P.$  If $x\in \partial^PD^0$ then by property (2), $a'_P = b_P$,
that is, $a'_Q=b_Q$ and $a_P^Q = b_P^Q.$  Flowing out under the geodesic action of $A'_Q$ we
see that $x \in \overline{T}_P \cap \overline{T}_Q \cap D$ implies that $a_P^Q = b_P^Q.$  For
such a point $x$ and for each $\alpha \in \Delta_Q$, we have
\begin{align*}
f_{i(\alpha)}^P(x)^{-1}i(\alpha)(b_P) &= i(\alpha)(a_Qa_P^Q)^{-1}i(\alpha)(b_Qb_P^Q) \\
&= \alpha(a_Q^{-1}b_Q) = f_{\alpha}^Q(x)^{-1}\alpha(b_Q)
\end{align*}
which shows that both equations agree on their common domain of definition, proving (1).  

By continuity, it suffices to prove (2) for points $x\in T(\overline{D}_Q).$  The
geodesic action by $t\in A_Q(\age 1)$ preserves the tiles in $T(\overline{D}_Q)$ so (2)
may be checked tile by tile.  If $x\in D^Q$ then $f_{\alpha}^Q(x\cdot t) =
f_{\alpha}^Q(x)\alpha(t)$ by (\ref{eqn-rootfunctionP}).  If $x\in D^P$ for some $\mathbf P
\subset \mathbf Q$, write $\Delta_P = i(\Delta_Q)\amalg J$ as in (\ref{eqn-decomposition}),
note that $A_Q(\age 1) \subset A_P(\age 1)$ and compute, for $\alpha \in \Delta_Q$,
\begin{equation*}f_{i(\alpha)}^P(x\cdot t) = f_{i(\alpha)}^P(x)i(\alpha)(t) =
f_{i(\alpha)}^P(x)\alpha(t)\end{equation*}
which proves the second statement.  Part (3) follows from (\ref{eqn-rootfunctionP})
for points $x\in \overline{T}_Q$ and from (\ref{eqn-roottranslate}) for points $x\in
\overline{T}_P.$ Part (4) follows from Lemma
\ref{lem-convergence} for points $x \in \overline{T}_Q$ and from (\ref{eqn-rootpi}) for points
$x\in \overline{T}_P.$ Part (5) follows from (\ref{eqn-roottranslate}) and part (6) is an
immediate consequence of the definition.
\qed

\subsection{Tiling of $\overline{X}$}\label{subsec-tilingofXbar}
Suppose $\mathbf b\in\mathcal B$ is a sufficiently \alarge regular parameter and
(\ref{eqn-tiling}) is
the associated $\Gamma$-invariant tiling.  Let $\bar\tau:\overline{D}\to \overline{X}$ denote
the projection to the reductive Borel-Serre compactification of $X=\Gamma \backslash D.$  If
$\mathbf P, \mathbf{P'}$ are rational parabolic subgroups of $\mathbf G$ then either $\bar
\tau(D^P) \cap \bar\tau(D^{P'}) = \phi$ or else they coincide.  Hence
the collection of images 
\begin{equation*}X^P=\bar\tau(D^P)
\end{equation*} form a decomposition of
$\overline{X}$ whose ``tiles'' are indexed by the set of $\Gamma$-conjugacy classes of
rational parabolic subgroups of $\mathbf G.$  Let $X_P^0 =X^G= \bar\tau (D_P^0)$ be the
``central tile.''  Denote by
\begin{equation}\label{eqn-nbhd}
T(\overline{X}_P) = \bar\tau\big( T(\overline{D}_P) \big) = \coprod_{\{R\} \subseteq P}
X^R
\end{equation}
the resulting neighborhood of $\overline{X}_P$ in $\overline{X}$, and by $\partial
T(\overline{X}_P)= \overline{T}(\overline{X}_P)-T(\overline{X}_P)$ its boundary.  (Here, $R$
runs through a set of representatives, one from each $\Gamma$-conjugacy class $\{R\}$ of
parabolic subgroups contained in $P$).    
For all $\alpha \in \Delta_P$ the  functions $r_{\alpha}^P$ pass to piecewise
analytic functions on $T(\overline{X}_P)$, which we also denote by $r_{\alpha}^P.$

%%%%%%%%%%%%%%%%%%%%%%%%%%%%%%%%  pictures  %%%%%%%%%%%%%%%%%%%%%%%%%

\begin{figure}[h!]               %% tiles
\begin{center}\begin{picture}(450,130)(30,-15)  %% (450,220)(0,0)
\setlength{\unitlength}{.11em}   %% {.1em} is the normal size
\linethickness{1pt}
\q(50,50)(75,25)(100,0)         %% left hand diagonal stratum
\q(100,0)(200,0)(300,0)         %% bottom stratum
\q(300,0)(325,25)(350,50)       %% right hand diagonal stratum

\linethickness{.1pt}

%% third level                
\q(150,24.72)(200,24.72)(240,24.72)
\q(240,24.72)(260,24.72)(268,27)
\q(268,27)(276,29.76)(286.06,33.63) %%RH curve

\q(286.06,33.62)(306.06,53.62)(336.06,83.62)  %%RH diagonal
\q(114,33.62)(94,53.62)(64,83.62)             %%LH diagonal

\q(132,27)(124,29.76)(114,33.63)
\q(160,24.72)(140,24.72)(132,27)
\q(250,24.72)(200,24.72)(160,24.72)          %%LH curve

\q(286,33.63)(286,10)(286,0)                %%right corner
\q(286,33.63)(296,23.63)(309.81,9.81)

\q(114,33.63)(114,10)(114,0)
\q(114,33.63)(104,23.63)(90.19,9.81)

%% labels
\put(100,0){\circle*{3}} \put(300,0){\circle*{3}}
\put(200,-12){$X_Q$}   \put(62,15){$X_{Q'}$}
\put(200,9){$X^Q$}    \put(74,35){$X^{Q'}$}
\put(200,45){$X^0$}
\put(95,-12){$X_P$}\put(300,-12){$X_{P'}$}
\put(98,9){$X^P$}
\put(286,9){$X^{P'}$}

\end{picture}\end{center}\caption{Tiles}
\label{fig-tiles}\end{figure}

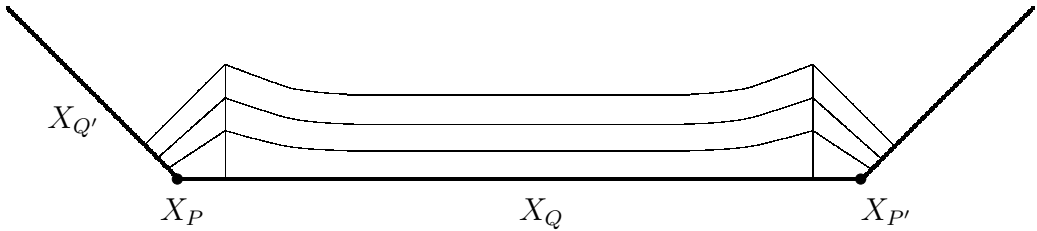
\begin{figure}[h!]               %% level curves of $r_{\alpha}$
\begin{center}\begin{picture}(450,85)(30,-15)  %% (450,220)(0,0)
\setlength{\unitlength}{.11em}   %% {.1em} is the normal size
\linethickness{1pt}
\q(50,50)(75,25)(100,0)         %% left hand diagonal stratum
\q(100,0)(200,0)(300,0)         %% bottom stratum
\q(300,0)(325,25)(350,50)       %% right hand diagonal stratum

\linethickness{.1pt}

%% third level                
\q(150,24.72)(200,24.72)(240,24.72)
\q(240,24.72)(260,24.72)(268,27)
\q(268,27)(276,29.76)(286.06,33.63) %%RH curve

\q(286,33.63)(296,23.63)(309.81,9.81)
\q(286,33.63)(286,10)(286,0) 

\q(132,27)(124,29.72)(114,33.63)
\q(160,24.72)(140,24.72)(132,27)
\q(250,24.72)(200,24.72)(160,24.72) %% LH curve

\q(114,33.63)(114,10)(114,0)
\q(114,33.63)(104,23.63)(90.19,9.81)

%% second level
\q(150,16.05)(200,16.05)(240,16.05)
\q(240,16.05)(260,16.05)(268,18)
\q(268,18)(276,20.31)(286,23.80) %% RH curve

\q(286,23.80)(296,14.9)(306,6)

\q(132,18)(124,20.31)(114,23.80)
\q(160,16.05)(140,16.05)(132,18)
\q(250,16.05)(200,16.05)(160,16.05) %% LH curve

\q(114,23.80)(104,14.9)(94,6)

%% first level
\q(150,8.28)(200,8.28)(240,8.28)
\q(240,8.28)(260,8.28)(268,9.61)
\q(268,9.61)(276,11.32)(286,14.19)

\q(286,14.19)(294.5,8.59)(303,3)

\q(132,9.61)(124,11.32)(114,14.19)
\q(160,8.28)(140,8.28)(132,9.61)
\q(250,8.28)(200,8.28)(160,8.28) %% LH curve

\q(114,14.19)(105.5,8.59)(97,3)
\put(100,0){\circle*{3}} \put(300,0){\circle*{3}}
\put(200,-12){$X_Q$}   \put(62,15){$X_{Q'}$}
\put(95,-12){$X_P$}\put(300,-12){$X_{P'}$}

\end{picture}\end{center}
\caption{Level sets of $r^Q_{\alpha}.$}
\label{fig-bentlevel}\end{figure} 

\subsection{Retraction and exhaustion}\label{subsec-retraction}
  Saper proves (\cite{Saper} \S 6.1) that there exists a unqiue
$\Gamma$-equivariant continuous and piecewise analytic ``geodesic'' retraction $R:D \to D^0$
which is the identity on $D^0$ such that, for all $y \in D^P$ and for all $t\in A_P(\age 1)$
the following holds:
\begin{equation}\label{eqn-retraction}
R(y\cdot t) = R(y).
\end{equation}
Then $R$ preserves tiles and it passes to a retraction which we also denote by $R:X \to X^0;$
it has the same property (\ref{eqn-retraction}). 

Define  $W:\overline{X} \to [0,1]$ by
\begin{equation*}
W(x) = \begin{cases}  1 & \text{ if }x \in X^0 \\
                      \underset{\alpha \in \Delta_Q}{\sup}\{ 1-r_{\alpha}^Q(x) \} &
\text{ if } x
\in X^Q \end{cases} \end{equation*}
We refer to $W$ as an {\it exhaustion function} because
\begin{equation*}
W^{-1}(0) = X^0 \text{  and  } W^{-1}(1) = \overline{X} - X.
\end{equation*}
The function $W$ is continuous (and piecewise analytic):  If $\mathbf P \subset \mathbf Q$,
write $\Delta_P = i(\Delta_Q) \amalg J$ as in (\ref{eqn-decomposition}).  Let $\overline{P} =
\nu_Q(P) \subset L_Q$ be the resulting parabolic subgroup of $L_Q$; then $\Delta_{\overline P}
= J.$  If $x\in \overline{X}^P \cap \overline{X}^Q$ then $\pi_Q(x) \in \overline{X}_Q^P \cap
\overline{X}_Q^0$ so $r_{\alpha}^P(x) = r_{\alpha}^{\overline{P}}(\pi_Q(x)) = 1$ for all
$\alpha \in J$ as in \S \ref{subsec-JPbar}.  Hence
\begin{equation*}
\underset{\alpha \in \Delta_P}{\sup}\{ 1-r_{\alpha}^P(x)\} = \underset{\alpha\in
\Delta_Q}{\sup}\{1 - r_{\alpha}^Q(x) \}.
\end{equation*} 
For each boundary stratum $X_Q$ the same constructions define a tile-preserving retraction 
\begin{equation}\label{eqn-RQ}
R_Q:X_Q \to X_Q^0\end{equation}
and an exhaustion function $W_Q: \overline{X}_Q \to [0,1]$ with $W_Q^{-1}(0) = X_Q^0$ and
$W_Q^{-1}(1) = \overline{X}_Q - X_Q.$  In fact, the stratum closure $\overline{X}_Q$ is tiled
by the collection of intersections $X_Q^P = \overline{X}_Q \cap X^P$ for $\mathbf P \subseteq
\mathbf Q$ and
\begin{equation}\label{eqn-exhaustion}
W_Q(x) = \begin{cases}  1 & \text{ if }x \in X_Q^0 \\
                      \underset{\alpha \in J}{\sup}\{ 1-r_{\alpha}^P(x) \} &
\text{ if } x
\in X_Q^P \end{cases} \end{equation}
where $\Delta_P = i(\Delta_Q) \amalg J.$

%%  \centerline{Diagram of level curves of $W$}

\subsection{Remarks}  We risk a certain amount of confusion by having defined
$r_{\alpha}^P(x)$ so as to decrease as $x \to X_P$ whereas the root function $f_{\alpha}^P(x)$
increases as $x \to X_P.$  
Although Saper \cite{Saper} actually constructs a tiling of the Borel-Serre
compactification $\widetilde{X}$  the same approach gives a tiling of the reductive
Borel-Serre compactification $\overline{X}$.  The collection $\{T(\overline{X}_P), \pi_P,
r_P=\text{max}_{\alpha \in \Delta_P}\{ r^P_{\alpha}\} \}$ of tubular
neighborhoods, tubular projections, and tubular distance functions are very much like a
``system of control data'' \cite{Mather, Gibson, SMT} for the stratified space $\overline{X}$,
but there are several important differences.  The functions $r_P$ are continuous and piecewise
analytic but are not smooth.  Whenever $\mathbf Q \subseteq \mathbf P$ we have $\pi_Q\pi_P =
\pi_Q$ however we do not have $r_Q\pi_P = r_Q$.  For this price we gain an especially strong
form of ``local triviality'' for the stratification of $\overline{X}$:  the neighborhood
$T(\overline{X}_P)$ is (homeomorphic to) a mapping cylinder neighborhood of the {\it closure}
of the stratum $X_P.$  In fact, it is possible to use the various geodesic actions to
construct a (piecewise analytic) homeomorphism between $T(\overline{X}_P)$ and the (open)
mapping cylinder of the projection $\pi_P: \partial T(\overline{X}_P) \to \overline{X}_P.$
(The open mapping cylinder of a mapping $\pi:A \to B$ is the quotient $(A \times [0,1) \coprod
B)/ \sim$ under the relation $(a,0) \sim \pi(a).$)  Analogous statements for other Satake
compactifications (such as the Baily-Borel compactification) are false.

\section{A Little Shrink}\label{sec-shrink}
\subsection{}  As in \S \ref{sec-notation}, $\mathbf G$ denotes a connected linear reductive
algebraic group defined over $\mathbb Q$, $D$ denotes the associated symmetric space,  $K' =
A_GK(x_0)$ is the stabilizer in $G$ of a fixed basepoint $x_0\in D$, $\Gamma \subset \mathbf
G(\mathbb Q)$ is an arithmetic group and $X = \Gamma \backslash D.$   In this section we
construct a homeomorphism $\overline{X} \to \overline{X}$ which moves a neighborhood of the
boundary towards the boundary.  When composed with a Hecke correspondence, this will have the
effect of chopping the fixed point set into pieces, each of which is contained in a single
stratum of $\overline{X}.$  The resulting behavior is much easier to analyze.  This ``shrink''
homeomorphism may be considered to be a topological analog to Arthur's truncation procedure.
\subsection{}\label{subsec-5.2}
Let $\mathbf Q \in \mathcal P_1$ be a standard proper maximal rational parabolic
subgroup of $\mathbf G.$  Fix $t \in A'_Q(\agt 1)$ so $\alpha(t) > 1$, where $\alpha \in
\Delta_Q$ is the unique simple root.  The geodesic action
of $t$ on $D$ extends continuously to the neighborhood $D[Q]$ of $D_Q$ in
the reductive Borel-Serre partial compactification $\overline{D}$ of $D$.
This geodesic action even extends continuously to the neighborhood
\begin{equation*}
D\{Q\} =\bigcup_{\mathbf P \subset \mathbf Q}D[P]
\end{equation*}
of the closure $\overline{D}_Q$ (where the union is taken over all
rational parabolic subgroups $\mathbf P \subset \mathbf G$ which
are contained in $\mathbf Q$). For if $\mathbf P \subset \mathbf Q$ 
then there is a canonical inclusion $i:A_Q \to A_P$.  The geodesic
action of the image $i(t) \in A_P$ agrees with the geodesic action
of $t\in A_Q$, so it also extends continuously to $D[P].$  We continue to
denote this action by $x \mapsto x\cdot t$ for $x\in D\{Q\}.$

Now fix a sufficiently \alarge $\Gamma$-invariant regular parameter $\mathbf b \in \mathcal B$
with its resulting tiling (\ref{eqn-tiling}) and partial distance functions
(\ref{eqn-partiald}) satisfying  Lemma \ref{lem-continuity}.  Let $T(\overline{D}_Q) \subset
D\{Q\}$ denote the neighborhood of $\overline{D}_Q$ (\ref{eqn-tilenbhd}) consisting of the
union of all tiles which intersect $\overline{D}_Q$ nontrivially, and let
$\overline{T}(\overline{D}_Q)$ denote its closure.  The geodesic action by $t\in A'_Q(>1)$
preserves $\overline{T}(\overline{D}_Q)$ since $\alpha(t) \agt 1.$  

Fix once and for all a smooth non-increasing function $\rho:[0,1] \to [0,1]$
with $\rho(r) = 1\ \iff\ r \le \frac{1}{2}$ and with $\rho(r)=0\ \iff r=1.$

\begin{figure}[h!]\begin{center}
\begin{picture}(180,95)(-5,-10)
\setlength{\unitlength}{.09em} 
\put(0,0){\vector(1,0){180}}
\put(0,0){\vector(0,1){70}}
\qbezier(0,60)(20,60)(80,60)
\qbezier(80,60)(105,60)(120,30)
\qbezier(120,30)(135,0)(160,0)

\put(186,-2){$r$} \put(-5,77){$\rho(r)$}
\put(80,0){\circle*{2}}  \put(160,0){\circle*{2}}
\put(78,-13){$\frac{1}{2}$}
\put(158,-11){$1$}
\put(0,60){\circle*{2}}\put(-10,55){$1$}

\end{picture}\end{center}\caption{The function $\rho$}
\end{figure}
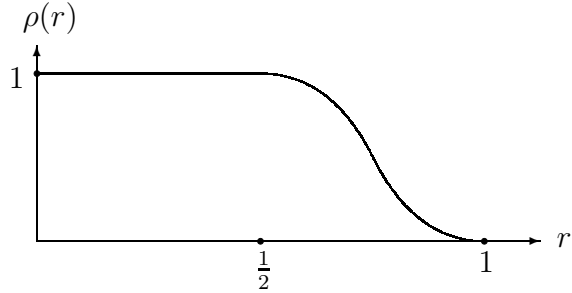

Let $r_{\alpha}^Q:T(\overline{D}_Q) \to [0,1)$ be the  partial distance function
(\ref{eqn-partiald}) which corresponds to the unique simple root $\alpha \in \Delta_Q.$
For $t \in A'_Q(\agt 1)$ define $Sh(Q,t):\overline{T}(\overline{D}_Q)\to
\overline{T}(\overline{D}_Q)$ by
\begin{equation}\label{eqn-shrink}
Sh(Q,t)(x) = x\cdot t^{\rho(r_{\alpha}^Q(x))}.
\end{equation}
Then by (\ref{eqn-rootfunctionP}) and (\ref{eqn-geodesic-r}),
\begin{align}\label{eqn-shrink3}
r_{\alpha}^Q(Sh(Q,t)(x)) &= r_{\alpha}^Q(x) \alpha(t)^{-\rho(r_{\alpha}^Q(x))}\\
f_{\alpha}^Q(Sh(Q,t)(x)) &= f_{\alpha}^Q(x) \alpha(t)^{\rho(r_{\alpha}^Q(x))}
\end{align}
for all $x \in \overline{T}(\overline{D}_Q).$  The quantity $\alpha(t)^{
\rho(r_{\alpha}^Q(x))}$ is bounded between $1$ and $\alpha(t).$  It equals $1$ if and
only if $x\in \partial T(\overline{D}_Q),$ that is, if and only if $r_{\alpha}^Q(x)= 1.$

If $\mathbf Q'$ is another maximal rational parabolic subgroup of $\mathbf G$, let us write
$\mathbf Q' \sim \mathbf Q$ if $\mathbf Q'$ is $\mathbf G(\mathbb Q)$ conjugate to $\mathbf
Q$.  In this case, any choice $g\in \mathbf G(\mathbb Q)$ of conjugating element induces the
same isomorphism $\mathbf{S_Q} \cong \mathbf{S_{Q'}}$ so we obtain a corresponding $t' \in
A'_{Q'}$ and a corresponding mapping $Sh(Q',t'):\overline{T}(\overline{D}_{Q'}) \to
\overline{T}(\overline{D}_Q).$  Patching these homeomorphisms together, define the {\it
shrink $Sh_{\mathbf Q}(t):\overline{D} \to \overline{D}$ corresponding to conjugates of the
standard parabolic subgroup $\mathbf Q$} by
\begin{equation*}
Sh_{\mathbf Q}(t)(x) = \begin{cases}
Sh(Q',t')x & \text{ if } x\in \overline{T}(\overline{D}_{Q'}) \text{ for some } \mathbf{Q'}
\sim
\mathbf{Q} \\
x & \text{otherwise} \end{cases}\end{equation*}
Then $Sh_{\mathbf Q}(t)$ is well defined and continuous because $T(\overline{D}_{Q'}) \cap
T(\overline{D}_Q) = \phi$ whenever $\mathbf Q' \sim \mathbf Q$ (and $\mathbf Q' \ne \mathbf
Q$), cf.\cite{Borel} \S 11.17 ({\it iii}).  Moreover, if $\gamma \in \Gamma$ then
\begin{equation*}
Sh_{\mathbf Q}(t)(\gamma x) = \gamma Sh_{\mathbf Q}(t)(x)
\end{equation*}
by Lemma \ref{lem-continuity}.  So (dividing by $\Gamma$), the homeomorphism $Sh_{\mathbf
Q}(t)$ passes to a homeomorphism which we denote in the same way, $Sh_{\mathbf
Q}(t):\overline{X} \to \overline{X}.$

\subsection{}
Suppose $\mathbf P \subseteq \mathbf Q$ is a rational parabolic subgroup of 
$\mathbf G$; set $\Delta_P = i(\Delta_Q) \amalg J$ as in (\ref{eqn-decomposition}).
It follows from (\ref{eqn-geodesic-r}) that for all $\beta\in J$ and for all $x\in
\overline{T}(\overline{D}_P)$ we have
\begin{equation}\label{eqn-shrink2}
r_{\beta}^P Sh(Q,t)(x) = r_{\beta}^P(x) \end{equation}
since $\beta(t)=1$ for any $\beta\in J.$  Now suppose $\mathbf{Q_1}, \mathbf{Q_2}$ are two
standard maximal rational parabolic subgroups of $\mathbf G$ whose intersection $\mathbf P =
\mathbf{Q_1} \cap \mathbf{Q_2}$ is parabolic.  Let $\alpha_i \in \Delta(Q_i)$ be
the unique nonzero roots.  Choose $t_i \in A'_i = A'_{Q_i}$ with
$\alpha_i(t_i) \agt 1$ and let 
\begin{equation*}Sh_i= Sh(Q_i,t_i): \overline{T}(\overline{D}_{Q_i}) \to
\overline{T}(\overline{D}_{Q_i})\end{equation*}
denote the resulting two shrinks.  It follows by taking $\mathbf P = \mathbf Q_1 
\cap \mathbf Q_2$ in (\ref{eqn-shrink2}) that {\it the mappings 
$Sh_1$ and $Sh_2$ commute on their common domain of definition,}
\begin{equation*}
\overline{T}(\overline{D}_P) = \overline{T}(\overline{D}_{Q_1}) \cap
\overline{T}(\overline{D}_{Q_2}).
\end{equation*}

\subsection{}\label{subsec-5.4}
Let $\mathbf {P_0} \subset \mathbf G$ be the standard minimal rational parabolic subgroup with
$\mathbf{S'_0} = \mathbf{S_{P_0}}/\mathbf{S_G}$ and with simple roots $\Delta =
\{\alpha_1,\alpha_2,\ldots,\alpha_r\}$ numbered in any order.  Each $j$ (with $1 \le j \le r$)
corresponds to a standard maximal proper rational parabolic subgroup $\mathbf{Q_j}$ with split
torus $\mathbf{S'_j}=\mathbf{S_{P_j}} / \mathbf{S_G}$ and identity component $A'_j =
\mathbf{S'_j}(\mathbb R)^0.$  Choose $\mathbf t \in A'_0$ to be
dominant and regular.  In other words, with respect to the canonical
{\it complementary} decomposition (cf
\ref{subsec-twoparab}) $A'_0 = \prod_j A'_j$ (with $\ker(\alpha_j) = \prod_{i \ne
j}\mathbf{S'_i}$) we may write $\mathbf t = t_1t_2\ldots t_r$ where $t_j \in A'_j$ and
$\alpha_j(t_j) = \alpha_j(\mathbf t) \agt 1.$  Define
$Sh(\mathbf t):\overline{D} \to \overline{D}$ to be the composition
\begin{equation*}
Sh(\mathbf t) = Sh_{\mathbf{Q_1}}(t_1)\circ Sh_{\mathbf{Q_2}}(t_2)\circ \ldots \circ
Sh_{\mathbf{Q_r}}(t_r)
\end{equation*}  
where each $Sh_{\mathbf {Q_j}}(t_j):\overline{D} \to \overline{D}$ is the shrink defined
above, corresponding to conjugates of the standard parabolic subgroup $\mathbf Q_j.$

\begin{prop}\label{prop-shrink}
  The mapping $Sh(\mathbf{t}):\overline{D} \to \overline{D}$ is independent of the
ordering $\Delta = \{\alpha_1,\alpha_2,\ldots,\alpha_r\}$ of the simple roots.  It is a
$\Gamma$-equivariant homeomorphism and passes to a homeomorphism $Sh(\mathbf t):\overline{X}
\to \overline{X}$ with the following properties:
\begin{enumerate}
\item It preserves the tiles and the strata, that is, for each rational
parabolic subgroup $\mathbf P \subseteq \mathbf G$, we have $Sh(\mathbf t)(X^P) = X^P$ and
$Sh(\mathbf t)(X_P)=X_P.$  
\item Within each tile, it is given by a geodesic action:  for each $x\in X^P$ there exists $b
= b_x \in A_P(\age 1)$ so that $Sh(\mathbf t)(x) = x\cdot b.$  
\item It is the identity on each central tile $X_P^0$ and $\pi_P(Sh(\mathbf t)(x)) = \pi_P(x)
\in X_P^0$ for all $x \in X^P.$
\item It commutes with the geodesic projection, that is, for any rational parabolic subgroup
$\mathbf P \subset \mathbf G$ and for each $x\in \overline{T}(\overline{X}_P)$ we have
$\pi_P(Sh(\mathbf t)(x)) = Sh(\mathbf t)(\pi_P(x)).$
\item  It is (globally) homotopic to the identity.
\item For any rational parabolic subgroup $\mathbf P \subset \mathbf G$ and for each $\alpha
\in \Delta_P$ and for each $x \in \overline{T}(\overline{X}_P),$ by equation
\parens{\ref{eqn-shrink3}} we have:
\begin{equation}\label{eqn-rshrink}
r_{\alpha}^P(Sh(\mathbf t)(x)) = r_{\alpha}^P(x) \alpha_0(\mathbf t)^{-\rho(r_{\alpha}^P(x))}
\end{equation}
where $\alpha_0 \in \Delta$ is the unique root which agrees with $\alpha$ after conjugation
and restriction to $\mathbf{S_P}.$  If $x \in X^P$ is constrained to lie in the single tile
$X^P$ then also
\begin{equation}\label{eqn-fshrink}
f_{\alpha}^P(Sh(\mathbf t)(x)) = f_{\alpha}^P(x) \alpha_0(\mathbf t)^{\rho(r_{\alpha}^P(x))}
\qed\end{equation}
\end{enumerate}\end{prop}
\noindent
We remark that the mapping $Sh(\mathbf t)$ {\it does} depend on the choice of regular
parameter (which determines the size of the tiles).

\section{Morphisms and Hecke correspondences}\label{sec-Hecke}
\subsection{}  As in \S \ref{sec-notation}, $\mathbf G$ denotes a connected linear reductive
algebraic group defined over $\mathbb Q$, $D$ denotes the associated symmetric space, and $K'
= A_GK(x_0)$ is the stabilizer in $G$ of a fixed basepoint $x_0\in D.$   Let $\Gamma, \Gamma'
\subset\mathbf G(\mathbb Q)$ be arithmetic subgroups and set $X = \Gamma \backslash D$ and $X'
=\Gamma' \backslash D.$
\begin{defn} \label{defn-morphism}
A mapping $f:X' \to X$ is a {\it morphism} if there exists $g\in \mathbf
G(\mathbb Q)$ such that \begin{enumerate}
\item $g\Gamma'g^{-1} \subset \Gamma$
\item $[\Gamma: g\Gamma' g^{-1}] < \infty $
\item $f(\Gamma'xK') = \Gamma gxK'$ for any $x\in G.$
\end{enumerate}\end{defn}
The morphism $f$ is determined by the pair $(\Gamma',g)$ by (3); it is well defined by (1).
For any $\gamma\in \Gamma$, $\gamma' \in \Gamma'$ the pair $(\Gamma', \gamma g \gamma')$
determine the same morphism.  If $\Gamma$ is torsion-free then $f$ is an unramified covering
of degree $[\Gamma: g\Gamma' g^{-1}]$ and it is locally an isometry with respect to the
invariant Riemannian metrics on $X'$ and $X$ induced from the Killing form.  Denote by
$\text{Mor}(X',X)$ the set of morphisms $X' \to X.$
\begin{lem}\label{lem-extension}  Each morphism $f\in \text{Mor}(X',X)$  admits unique
continuous extensions $\widetilde{f}:\widetilde {X}' \to \widetilde {X}$ to the Borel-Serre
compactification and $\overline{f}:\overline{X}' \to \overline{X}$ to the reductive
Borel-Serre compactification.  The mappings $\widetilde{f}$ and $\overline{f}$ are finite, and
they restrict to morphisms on each boundary stratum.  If $\overline{f}(X'_P) = X_Q$, if $U$
and $U'$ are $\Gamma'$ and $\Gamma$-parabolic neighborhoods of $X'_P$ and $X_Q$ in
$\overline{X}'$ and $\overline{X}$ respectively, then
\begin{equation}\label{eqn-barfpi}
\overline{f}(\pi_P(x)) = \pi_Q(\overline{f}(x)) \end{equation}
for all $x \in U' \cap \overline{f}^{-1}(U).$
\end{lem}
We say that a mapping $\overline{X}' \to \overline{X}$ is a morphism if it is the unique
continuous extension of some $f \in \text{Mor}(X',X).$

\subsection{Proof of Lemma \ref{lem-extension}}\label{pf-extension}  
Suppose the morphism $f:X' \to X$ is given by the pair $(\Gamma',g).$ Let $T_g:D \to D$ denote
the action of $g$ on $D$. It moves the basepoint $x_0$ to a new basepoint $x_1 = gx_0$ with
stabilizer $K'(x_1) = {}^gK' =gK'g^{-1}.$  If $\mathbf P$ is a rational parabolic subgroup and
if $\mathbf Q =\mathbf{{}^gP}=g\mathbf P g^{-1}$ set $K'_P(x_0) = K' \cap P$ and $K'_Q(x_1) =
K'(x_1)\cap Q.$ Then $T_g$ may also be described as the mapping
\begin{equation}\label{eqn-morphism}
D = P/K'_P(x_0) \to Q/K'_Q(x_1) = D\end{equation}
which is given by $xK'_P(x_0) \mapsto gxg^{-1}K'_Q(x_1)$ by (\ref{eqn-basepoint2}).  This
intertwines the geodesic actions of $A'_P$ and $A'_Q$, that is, 
\begin{equation}\label{eqn-intertwine}
T_g(x\cdot a) =gxg^{-1}gi_{x_0}(a)g^{-1}K'_Q(x_1)=gxg^{-1}i_{x_1}(\hat a)K'_Q(x_1)= 
T_g(x) \cdot \hat{a}\end{equation}
where $a \mapsto \hat{a}$ is the canonical identification $A'_P \cong A'_Q$ of \S
\ref{subsec-roots}.  It follows that $T_g$ extends to a
mapping $\widetilde T_g:\widetilde{ D} \to \widetilde{ D}$ on the Borel-Serre partial
compactification, which takes the boundary component $e_P=P/K'_PA'_P$ to $e_Q=Q/K'_QA'_Q$ by
conjugation and satisfies
\begin{equation}\label{eqn-intertwine2}
\pi_Q \widetilde{T}_g (x) = \widetilde{T}_g (\pi_P(x)).
\end{equation} 

 The mapping $\widetilde T_g$ passes to a mapping $\widetilde{f}:\Gamma'\backslash \widetilde{
D}\to\Gamma \backslash \widetilde{ D}$ which is the desired extension.  It maps $Y'_P =
\Gamma'_P\backslash e_P$ to $Y_Q = \Gamma_Q \backslash e_Q$ by $\widetilde{f}(\Gamma'_P
xK'_P(x_0)A'_P)= \Gamma_Q gxg^{-1}K'_Q(x_1)A'_Q,$ which is a mapping of degree
$[\Gamma_Q:g\Gamma'_Pg^{-1}] <\infty.$  (Here, $\mathbf Q = {}^g\mathbf P$.)  The extension
$\widetilde{f}$ may map several strata of $\widetilde {X}'$ to a single stratum of $\widetilde
{X}$:  let $\mathbf {Q_1},\mathbf{Q_2}, \ldots, \mathbf {Q_m}$ be a set of representatives for
the ${}^g\Gamma'$-conjugacy classes of rational parabolic subgroups which are
$\Gamma$-conjugate to $\mathbf Q$, and set $\mathbf P_j = g^{-1}\mathbf Q_jg.$  Then
$\widetilde {f}$ maps each stratum $Y'_{P_j} \subset \widetilde{ X}'$ to the stratum
$Y_Q\subset \widetilde{X}$ by a morphism which may be described in a manner similar to
(\ref{eqn-morphism}).  This shows that $\widetilde{f}$ is finite and that its restriction to
each boundary stratum is a morphism.

Similarly the mapping $\widetilde{T}_g$ passes to a mapping $\overline{T}_g:\overline{D} \to
\overline{D}$ on the reductive Borel-Serre compactification of $D$, which further passes to a
mapping $\overline{f}: \Gamma' \backslash \overline{D} \to \Gamma \backslash \overline{D}.$
Then $\overline{f}$ maps $X'_P = \Gamma'_P \backslash P /K_PA_P\mathcal U_P $ to $X_Q =
\Gamma_Q \backslash Q /K_Q A_Q \mathcal U_Q$ by 
\[ \overline{f}(\Gamma'_P xK_P(x_0)A_P\mathcal U_P) = \Gamma_Q gxg^{-1} K_Q(x_1)
A_Q \mathcal U_Q. \]
The degree of this mapping is not obviously finite because the intersection $\Gamma'_P \cap
\mathcal U_P$ is nontrivial.  By (\ref{eqn-basepoint}) conjugation by $g$ takes $L_P(x_0)$ to
$L_Q(x_1).$  Let $K_{L(P)}(x_0)A_P$ be the stabilizer in $L_P$ of the basepoint $\pi_P(x_0)
\in D_P$ and let $K_{L(Q)}(x_1)A_Q$ be the stabilizer in $L_Q$ of the basepoint $\pi_Q(x_1)
\in D_Q.$  Set $\Gamma'_{L(P)} = \nu_P(\Gamma'_P) \subset L_P$ and
$\Gamma_{L(Q)} = \nu_Q(\Gamma_Q) \subset L_Q.$  Then $X'_P = \Gamma'_{L(P)} \backslash L_P /
K_{L(P)}(x_0)A_P$ and $X_Q = \Gamma_{L(Q)}\backslash L_Q / K_{L(Q)}(x_1)A_Q$ with respect to
which we may express $\overline{f}$ as follows:
\[\overline{f}(\Gamma'_{L(P)} x K_{L(P)}(x_0)A_P) = 
\Gamma_{L(Q)}gxg^{-1}K_{L(Q)}(x_1)A_Q \]
which has degree $[\Gamma_{L(Q)}: g\Gamma'_{L(P)}g^{-1}] < \infty.$  As in the preceding
paragraph, the mapping $\overline{f}$ will take each of the finitely many strata $X'_{P_j}$ to
the stratum $X_Q$ (for $1 \le j \le m$) by a similarly defined finite morphism.  \qed

\begin{defn}\label{def-Hecke}
  A {\it correspondence} on $X = \Gamma \backslash G/K'$ is an arithmetic subgroup
$\Gamma' \subset \mathbf G(\mathbb Q)$ together with two morphisms $c_1,c_2\in
\text{Mor}(C,X)$, where $C = \Gamma' \backslash D.$    
 A point $x\in C$ is {\it fixed} if $c_1(x) = c_2(x).$  Two correspondences
$(c_1,c_2):C\rightrightarrows X$ and $(c'_1,c'_2):C'\rightrightarrows 
X$ are said to be {\it isomorphic} if there is an 
invertible morphism $\alpha :C\to C'$ such that $c'_j\circ\alpha 
=c_j$ (for $j=1,2$). \end{defn}
Each $g\in \mathbf G(\mathbb Q)$ gives rise to a {\it Hecke correspondence} $
C=C[g]\rightrightarrows X$ as follows:  set $\Gamma[g]=\Gamma\cap g^{-1}\Gamma g,$
$C=\Gamma[g]\backslash G/K',$ and define
\begin{equation}\label{eqn-Hecke} 
(c_1,c_2)(\Gamma[g]xK')=(\Gamma xK',\Gamma gxK').  \end{equation}

Modifying $g$ by an element of $\mathbf{S_G}(\mathbb Q)$ does not change the Hecke
correspondence.  By Lemma \ref{lem-extension} each correspondence $C \rightrightarrows X$ has
a unique continuous extension $\overline{C} \rightrightarrows \overline{X}$ to the reductive
Borel-Serre compactification, and each isomorphism $\alpha:C \to C'$ of correspondences $C
\rightrightarrows X$, $C' \rightrightarrows X$ extends uniquely to an isomorphism
$\overline{C} \to \overline{C}'$ of the extended correspondences.
 
\begin{lem}\label{lem-isomorphism} %4.2
Let $X = \Gamma \backslash D$ and let $g \in \mathbf{G}(\mathbb Q).$  The isomorphism class of
the resulting Hecke correspondence $\overline{C}=\overline{C}[g]\rightrightarrows
\overline{X}$ depends only on the double coset $\Gamma g\Gamma\in\Gamma\backslash 
\mathbf G(\mathbb Q)/\Gamma$.   \end{lem}
\subsection{  Proof  }  %4.3
If suffices to verify the statement for the correspondence $C \rightrightarrows X$ since the
extension to the reductive Borel-Serre compactification exists uniquely.  Let
$\gamma_1,\gamma_2\in\Gamma$ and let $g'=\gamma_1g\gamma_2$ be another element in the
same double coset $\Gamma g \Gamma.$  Set $\Gamma[g']=\Gamma\cap g^{
\prime -1}\Gamma g',$ $C'=\Gamma[g']\backslash G/K',$ and define
$(c'_1,c'_2):C'\rightrightarrows X$ by $c'_1(\Gamma[g']xK')=\Gamma
xK',$ $c'_2(\Gamma[g']xK')=\Gamma g'xK'$.
One verifies by direct calculation that the morphism $f:C'\to C$ which is given by
\begin{equation}\label{eqn-isomorphism}
f(\Gamma[g']xK')=\Gamma[g]\gamma_2xK' \end{equation} %4.3.1
is a well-defined isomorphism of correspondences, with inverse given by
 $f^{-1}(\Gamma[g]xK')=\Gamma[g']\gamma_2^{-1} xK'$.  \qed

\subsection{Remark}  It can be shown that the mapping $(c_1,c_2):C \to X \times X$ is
generically one-to-one.  In the event that every element of $g^{-1}\Gamma g \Gamma$ is neat,
then this mapping is globally an embedding.

\quash{
If every element of $g^{-1}\Gamma g \Gamma $ is
neat, then the mapping $(c_1,c_2):C \to X \times X$ is an embedding.

Next we will show that the mapping $(c_1,c_2):C\to X\times X$ is an embedding.
Suppose two points $\Gamma[g]xK',\Gamma[g]x'K'\in C$ map to the same point in $
X\times X$.  
Thus, $c_1(\Gamma[g]xK')=c_1(\Gamma[g]x'K')$, or $\Gamma xK'=\Gamma x'K',$ so
$x'=\gamma_1xk_1 $
for some $\gamma_1\in\Gamma$ and some $k_1\in K'.$ Similarly, $c_2
(\Gamma[g]xK')=c_2(\Gamma[g]x'K')$ gives
$gx'=\gamma_2gxk_2 $
for some $\gamma_2\in\Gamma$ and some $k_2\in K'$.  Comparing these gives
\begin{equation}g^{-1}\gamma_2^{-1}g\gamma_1=xk_2k_1^{-1}x^{-1} \end{equation} %4.3.6
The right hand side is in a compact subgroup $xK'x^{-1}\subset G$ while the 
left hand side is neat by assumption, so 
both elements equal $1$.  This gives
\begin{equation}\gamma_1=g^{-1}\gamma_2g\in\Gamma\cap g^{-1}\Gamma g \end{equation} %4.3.7
so $\gamma_1\in\Gamma[g].$ Therefore $\Gamma[g]x'K'=\Gamma[g]\gamma_1xK'=
\Gamma[g]xK'$ and the two points are the
same.  \qed    }  %% end quash

%\subsection{   Standard form } %4.5
In the next section we will show that every correspondence is a
covering of a Hecke correspondence.

\begin{prop}\label{prop-covering correspondence} % 4.6
Let $\Gamma'\subset \mathbf G(\mathbb Q)$ be an arithmetic subgroup, $C' =
\Gamma' \backslash D$ and $(c'_1,c'_2):C' \rightrightarrows X$ be a correspondence.
Then there
is a Hecke correspondence $(c_1,c_2):C[g]\rightrightarrows
X$ and a subgroup $\Gamma^{\prime\prime} \subset \Gamma[g]$ such that the
correspondence $C'\rightrightarrows X$ is isomorphic to the correspondence
\begin{equation}
\begin{CD}
C^{\prime\prime}@>{h}>> C[g] \rightrightarrows X\end{CD}\end{equation}
where $C^{\prime\prime}=\Gamma^{\prime\prime}\backslash G/K$
and $h(\Gamma^{\prime\prime}xK) = \Gamma[g]xK.$  \end{prop}

\subsection{  Proof }  %4.7
Suppose $c'_1(\Gamma'xK')=\Gamma g_1xK'$ and $c'_2(\Gamma'x
K')=\Gamma g_2xK'$ where $g_j\Gamma'g_j^{-1}{\subset}\Gamma$ are subgroups of finite index. 
Then $g=g_2g_1^{-1}$ determines a Hecke correspondence 
$C[g]=\Gamma [g]\backslash D\rightrightarrows X$.  Define $\Gamma^{\prime\prime}=g_1\Gamma'
g_1^{-1}$ and $C^{\prime\prime} = \Gamma^{\prime\prime}\backslash D.$  
Since $\Gamma^{\prime\prime}\subset \Gamma[g]$ we obtain a correspondence  
in ``standard form'',
\begin{equation}\begin{CD}
C^{\prime\prime} @>h>> C[g] \rightrightarrows X\end{CD}\end{equation}
with
$h(\Gamma^{\prime\prime}xK') = \Gamma[g]xK'$.  Define $f:C' \to
C^{\prime\prime}$ by $f(\Gamma'xK') = \Gamma^{\prime\prime}g_1xK'$.  Then $f$
is well defined, and it is easily seen to be an isomorphism of correspondences.  \qed

\subsection{Narrow tilings}\label{subsec-narrow}
 Let $(c_1,c_2):\overline{C} \rightrightarrows \overline{X}$ be a
Hecke correspondence defined by some element $g\in \mathbf G(\mathbb Q),$ so $C =
\Gamma'\backslash D$ with $\Gamma' = \Gamma \cap g^{-1}\Gamma g.$  Let $b\in \mathcal B$ be a
sufficiently \alarge $\Gamma$-invariant regular parameter.  Then it is also
$\Gamma'$-invariant, it gives rise to tilings $\{ C^Q\}$ of $\overline{C}$ and $\{X^Q\}$ of
$\overline{X},$ and the mapping $c_1:\overline{C} \to \overline{X}$ takes tiles to tiles
(although the same cannot necessarily be said of $c_2$).  Let us say this tiling is {\it
narrow} with respect to the Hecke correspondence if, for every stratum $C_Q$ of
$\overline{C}$, the following holds:
\begin{equation*}
c_1(\overline{T}(\overline{C}_Q)) \cap c_2(\overline{T}(\overline{C}_Q)) \ne \phi
\iff c_1(C_Q) = c_2(C_Q)
\end{equation*}
and if, in this case, $c_1(\overline{T}(\overline{C}_Q)) \cup
c_2(\overline{T}(\overline{C}_Q))$ is a $\Gamma$-parabolic neighborhood of $\overline{X}_Q$ in
$\overline{X}.$  
\begin{prop}\label{prop-narrow}
Fix a Hecke correspondence $\overline{C} \rightrightarrows \overline{X}.$  If the
$\Gamma$-invariant regular parameter $b \in \mathcal B$ is chosen sufficiently \alarge then
the resulting tiling $\{C^Q\}$ of $\overline{C}$ is narrow for that Hecke correspondence.
\end{prop}
\subsection{Proof}
Let $C_Q$ be a stratum of $\overline{C}$ and suppose $c_1(C_Q) = X_P$ and $c_2(C_Q)= X_{P'}.$
Then $\mathbf Q$ is $\Gamma$-conjugate to $\mathbf P$ while $g\mathbf Q g^{-1}$ is
$\Gamma$-conjugate to $\mathbf{P'}$  (cf. Lemma \ref{lem-Xi}). In particular, $\mathbf P$ and
$\mathbf{P'}$ are $\mathbf G(\mathbb Q)$-conjugate, which implies that either $X_P=X_{P'}$
(and $\mathbf P = \mathbf{P'}$) or $\overline{X}_P \cap \overline{X}_{P'} = \phi$
(\ref{eqn-nointersect}).  In the latter case there exist neighborhoods $U$ of $\overline{X}_P$
and $U'$ of $\overline{X}_{P'}$ which do not intersect.  Choose the tiling parameter so
\alarge that $\overline{T}(\overline{C_Q}) \subset c_1^{-1}(U) \cap c_2^{-1}(U').$  Since
there are finitely many strata $C_Q$ in $\overline{C}$, this amounts to finitely many
conditions on the tiling parameter.  Similarly, if $c_1(C_Q) = c_2(C_Q)=X_Q$, choose any
parabolic neighborhood $U\subset \overline{X}$ of $\overline{X}_Q$ and then choose the tiling
so small that $\overline{T}(\overline{C}_Q) \subset c_1^{-1}(U) \cap c_2^{-1}(U).$  This will
guarantee that $c_1(\overline{T}(\overline{C}_Q)) \cup c_2(\overline{T}(\overline{C}_Q))$ will
be a $\Gamma$-parabolic set in $\overline{X}.$
\qed

\section{ Restriction to the boundary }\label{sec-restriction}

\subsection{Parabolic Hecke correspondence}\label{subsec-restriction}As in \S
\ref{sec-notation}, $\mathbf G$ denotes a connected linear reductive algebraic group defined
over $\mathbb Q$, $D$ denotes the associated symmetric space,  $K' = A_GK(x_0)$ is the
stabilizer in $G$ of a fixed basepoint $x_0\in D$, $\Gamma \subset \mathbf G(\mathbb Q)$
is an arithmetic group and $X = \Gamma \backslash D.$ 
Fix a rational parabolic subgroup $\mathbf P \subset \mathbf G$ and let $X_P=\Gamma_P
\backslash D_P\subset \overline{X}$ be the corresponding stratum in the reductive Borel-Serre
compactification of $X$.  Each $y \in \mathbf P(\mathbb Q)$ determines a correspondence on a
$\mathbf{P}(\mathbb Q)$-invariant neighborhood of $X_P$ which we now describe.  Set $\Gamma'_P
= \Gamma_P[y]=\Gamma_P \cap y^{-1}\Gamma_P y.$  Define the {\it parabolic Hecke
correspondence}
\begin{equation}\label{eqn-Pcorrespondence}
(c_1,c_2): \Gamma'_P\backslash D[P] \rightrightarrows \Gamma_P \backslash D[P]
\end{equation}
{\it determined by $y \in \mathbf P(\mathbb Q)$}
to be the unique continuous extension of the correspondence $\Gamma'_P\backslash D
\rightrightarrows \Gamma_P \backslash D$ which is given by
\begin{equation}\label{eqn-correspondence}
\Gamma'_PxK'_P \mapsto (\Gamma_PxK'_P, \Gamma_PyxK'_P)
\end{equation}
where we identify $D = P/K'_P.$  It follows from (\ref{eqn-intertwine}) (by taking $\mathbf{P}
= \mathbf{Q}$) that this correspondence commutes with the geodesic action of $A'_P$, that is,
\begin{equation}\label{eqn-commutes geodesic}
c_i(x \cdot a) = c_i(x)\cdot a\end{equation}
(for $i=1,2$) for any $x \in \Gamma'_P\backslash D[P]$ and
for any $a \in A'_P.$  Therefore the parabolic Hecke correspondence preserves the corner
structure near $C_P$, that is, if $\mathbf Q \supset \mathbf P$ is a rational parabolic
subgroup then each mapping $c_i$ takes the stratum $\Gamma'_P \backslash D_Q \subset \Gamma'_P
\backslash D[P]$ to the stratum $\Gamma_P \backslash D_Q \subset \Gamma_P \backslash D[P].$ 

There is also an associated (global) correspondence, $\overline{C} = \Gamma' \backslash
\overline{D} \rightrightarrows \overline{X} = \Gamma \backslash \overline{D}$ (where $\Gamma'
= \Gamma \cap y^{-1} \Gamma y$).  If $V \subset D[P] \subset \overline{D}$ is a
$\Gamma$-parabolic neighborhood of $D_P$ then it is also a $\Gamma'$-parabolic neighborhood of
$D_P$, as is $y^{-1} \cdot V.$  Hence $V \cap y^{-1} V$ is also a $\Gamma'$-parabolic
neighborhood of $D_P.$  It follows that, if $U \subset \Gamma_P \backslash D[P] \subset
\Gamma_P \backslash \overline{D}$ is a $\Gamma$-parabolic neighborhood of $X_P$ then $U' =
c_1^{-1}(U) \cap c_2^{-1}(U) \subset \Gamma'_P \backslash \overline{D}$ is a
$\Gamma'$-parabolic neighborhood of $C_P = \Gamma'_P \backslash D_P.$  We will say that any
correspondence isomorphic to such a $U' \rightrightarrows U$ is {\it modeled} on the
parabolic Hecke correspondence (\ref{eqn-Pcorrespondence}).
\begin{equation*}\begin{CD}
U' &\ \subset\ & \Gamma'_P \backslash D[P] &\ \rightrightarrows\ &
\Gamma_P \backslash D[P] &\ \supset\ & U \\
|| && \bigcap && \bigcap && || \\
U' &\ \subset\ & \Gamma'_P \backslash \overline{D} &\ \rightrightarrows\ &
\Gamma_P \backslash \overline{D} &\ \supset\ & U \\
@V{\cong}VV @VV{\beta'}V @VV{\beta}V @VV{\cong}V \\
U' &\ \subset\ & \overline{C} &\ \rightrightarrows\ & \overline{X} &\ \supset\ & U
\end{CD}\end{equation*}

\subsection{}
Now suppose $g\in \mathbf G(\mathbb Q)$ gives rise to the Hecke correspondence
$C=C[g]\rightrightarrows X$ with its canonical extension $(c_1,c_2):\overline{C}
\rightrightarrows
\overline{X}$, where $C = \Gamma[g]\backslash D$ as in \S \ref{eqn-Hecke}.  If $\mathbf P$ is
a rational parabolic subgroup of $\mathbf G$ and if $X_P$ denotes the corresponding RBS
stratum, then we may consider the part $c_1^{-1}(X_P)$ of $\overline{C}$ which lies over this
stratum.  It will consist of {\it several} RBS boundary strata $C_Q$ of $C$.  Some of these
boundary strata may be mapped back to $X_P$ by the mapping $c_2$.  In this case, we shall say
that the Hecke correspondence $\overline{C}$ has a {\it restriction} to $X_P$ consisting of
the union of those boundary strata $C_Q$ such that $(c_1,c_2)|C_Q:C_Q \rightrightarrows X_Q.$  

\begin{prop}\label{prop-restriction}
Let $\Gamma \subset \mathbf G(\mathbb Q)$ be a neat arithmetic group.  Let
$(c_1,c_2):C=C[g]\rightrightarrows X$ be the Hecke correspondence which is determined by an
element $g\in \mathbf G(\mathbb Q).$  Let $\mathbf P$ be a rational parabolic subgroup of
$\mathbf G$, with corresponding boundary stratum $X_P \subset \overline{X}$.  Decompose the
intersection $\Gamma g \Gamma \cap P$ into a union of $\Gamma_P$-double cosets,
\begin{equation}\label{eqn-cosets}
\Gamma g \Gamma \cap P \ = \coprod_{j=1}^m \Gamma_P g_j \Gamma_P \end{equation}
with $g_j \in \mathbf P(\mathbb Q)$.  Then $m < \infty$ and, over a sufficiently small
parabolic neighborhood of $X_P$, the Hecke correspondence $\overline{C} \rightrightarrows
\overline{X}$  breaks into a disjoint union of $m$ correspondences which are given by $g_j$
and which are modeled on the parabolic Hecke correspondences
\begin{equation*}
\Gamma_P[g_j] \backslash D[P] \rightrightarrows \Gamma_P \backslash D[P] \end{equation*}
for $j = 1, 2, \ldots, m,$ where $\Gamma_P[g_j] = \Gamma \cap g_j^{-1}\Gamma g_j \cap P.$ 
\end{prop}

(A similar procedure is described in the adelic setting in \cite{Harder}.)  The proof will
take the rest of \S \ref{sec-restriction}.  First we establish a one to one
correspondence between the components of the restriction of the Hecke correspondence
 to $X_P$ and the double cosets which appear in (\ref{eqn-cosets}).
% \begin{defn}{}\label{def-Xi} %6.3
Let $\Gamma_Py\Gamma_P\subset\Gamma g\Gamma\cap P$ be a double coset from
(\ref{eqn-cosets}).   Write $y=\gamma_2g\gamma_1$  for some $\gamma_1,\gamma_2\in\Gamma$.
Define 
\begin{equation}\Xi (\Gamma_P y \Gamma_P)=\gamma_1\mathbf P\gamma_1^{-1}
= g^{-1}\gamma_2^{-1}\mathbf P \gamma_2 g.\end{equation}

\begin{lem}\label{lem-Xi} %6.4
The mapping $\Xi$ gives a well defined one to one correspondence between 
\renewcommand{\theenumi}{\alph{enumi}}
\renewcommand{\labelenumi}{(\theenumi)}
\renewcommand{\theenumii}{\roman{enumii}}
\renewcommand{\labelenumii}{(\theenumii)}
\begin{enumerate}
\item Double cosets $\Gamma_Py\Gamma_P\subset\Gamma g\Gamma\cap
P$
\item $\Gamma[g]$-conjugacy classes of rational parabolic subgroups
$\mathbf Q=\Xi (\Gamma_P y \Gamma_P)\subset \mathbf G$ such that 
\begin{enumerate}
\item $\mathbf Q$ is $ \Gamma$-conjugate to $\mathbf P$ and
\item $g\mathbf Qg^{-1}$
is $\Gamma$-conjugate to $\mathbf P$ \end{enumerate}  
\item Boundary strata $C_Q\subset\overline{C}$ such that $(c_1,c_2)
:C_Q\rightrightarrows X_P.$ \end{enumerate} 
In particular, this set is finite.\end{lem}
\renewcommand{\theenumi}{\arabic{enumi}}
\renewcommand{\labelenumi}{\theenumi.}

\subsection{    Proof of lemma \ref{lem-Xi} } % 6.5
First compare the sets (b) and (c).  The boundary strata in $
\overline {C}'$  are in one to one correspondence with $\Gamma[g]$-conjugacy 
classes of rational parabolic subgroups of $\mathbf G$, while the boundary 
strata in $\overline {X}$ are in one to one correspondence with $\Gamma$ conjugacy
classes of rational parabolic subgroups.  Condition (i) is equivalent to the
statement that $c_1$ maps $C_Q$ to $X_P$ while condition (ii) is equivalent
to the statement that $c_2$ maps $C_Q$ to $X_P$.

Now verify that the mapping $\Xi$ is well defined, i.e.~that the $\Gamma[g]$-conjugacy class
of $\Xi(\Gamma_P g_j \Gamma_P) = \gamma_1\mathbf P \gamma_1^{-1}$ is independent of the
choices.  Let $y' = \gamma'_2 g \gamma'_1 \in \Gamma_P g_j \Gamma_P$ be another element in
the
same double coset, and set $\mathbf Q'=\gamma'_1 \mathbf P (\gamma'_1)^{-1}$.  
Since $y'\in\Gamma_Py\Gamma_P,$ there exists $\gamma_a,\gamma_b\in
\Gamma_P$ such that $y'=\gamma_ay\gamma_b.$ So 
$y'=\gamma_2'g\gamma_1'=\gamma_a\gamma_2g\gamma_1\gamma_b$, which gives
\begin{equation*}h:=g^{-1}\gamma_2^{\prime -1}\gamma_a\gamma_2g=\gamma_1\gamma_b\gamma_
1^{\prime -1}\in g^{-1}\Gamma g\cap\Gamma =\Gamma[g]. \end{equation*} %6.5.3
Then,
\begin{equation}h^{-1}\mathbf Qh=(\gamma_1\gamma_b\gamma^{\prime -1}_1)^{-1}(\gamma_
1\mathbf P\gamma_1^{-1})(\gamma_1\gamma_b\gamma^{\prime -1}_1)=\mathbf Q'\end{equation}
which verifies that $\mathbf Q$ and $\mathbf Q'$ are $\Gamma[g]$ conjugate.

Next we show that $\Xi$ is surjective.  Suppose that $\mathbf Q$ and $g\mathbf Q
g^{-1}$ are both $\Gamma$-conjugate to $\mathbf P$.  Say, $Q=\gamma_1 \mathbf P
\gamma_1^{-1}$ and $g\mathbf Q g^{-1} = \gamma_2^{-1}\mathbf P \gamma_2$ for some
$\gamma_1,\gamma_2 \in \Gamma.$  Then
\begin{equation}\gamma_2g\gamma_1\mathbf P\gamma_1^{-1}g^{-1}\gamma_2^{-1}=\mathbf
P\end{equation}
so the element $y=\gamma_2g\gamma_1 \in \mathbf P(\mathbb Q)\cap \Gamma g \Gamma$ (since a
parabolic subgroup is its own normalizer) and  $\Xi(\Gamma_P y \Gamma_P) = \mathbf Q.$

Finally we show that $\Xi$ is injective.  Suppose $y, y' \in \Gamma g \Gamma \cap P$, say
$y=\gamma_2g\gamma_1$ and $y' = \gamma'_2 g \gamma'_1.$  Set
\begin{alignat*}{2}
\mathbf Q &= \gamma_1 \mathbf P \gamma_1^{-1} &\ 
&= g^{-1}\gamma_2^{-1}\mathbf P \gamma_2 g\\
\mathbf Q' &= \gamma'_1 \mathbf P (\gamma'_1)^{-1} &\
   &= g^{-1}(\gamma'_2)^{-1}\mathbf P \gamma_2 g \end{alignat*}
Suppose $\mathbf Q$ and $\mathbf Q'$ are $\Gamma[g]$-conjugate, say $\mathbf Q = \gamma
\mathbf Q' \gamma^{-1}$ for some $\gamma \in \Gamma \cap g^{-1}\Gamma g.$  Comparing these
gives two relations,
\begin{alignat*}{2}
\mathbf Q &= \gamma_1\mathbf P \gamma_1^{-1}  &\ 
&= \gamma \gamma'_1 \mathbf P (\gamma'_1)^{-1} \gamma^{-1} \\
\mathbf Q &= g^{-1}\gamma_2^{-1}\mathbf P \gamma_2 g &\  
&= \gamma g^{-1}(\gamma'_2)^{-1} \mathbf P \gamma'_2 g \gamma^{-1} \end{alignat*}
from which it follows that $h_1 = (\gamma'_1)^{-1}\gamma^{-1} \gamma_1^{-1} \in P$ and
$h_2 = \gamma_2 g \gamma g^{-1}(\gamma'_2)^{-1} \in P.$  Moreover, $h_1,h_2 \in \Gamma$.  But
$h_2y'h_1 = y$ hence $\Gamma_P y' \Gamma_P = \Gamma_P y \Gamma_P$ as claimed.  \qed

\subsection{Proof of Proposition \ref{prop-restriction}}
By lemma \ref{lem-Xi}, the restriction of the Hecke correspondence $\overline{C} =
\overline{C}[g]\rightrightarrows \overline{X}$ to the stratum $X_P$ breaks into a union of $m$
correspondences, indexed by the elements $g_1,g_2,\ldots, g_m$.  Fix $j$ (with $1 \le j \le
m$) and set $\Gamma[g_j] = \Gamma \cap g_j^{-1}\Gamma g_j$ and $\Gamma_P[g_j] = \Gamma[g_j]
\cap P.$  The following commutative diagram of correspondences provides an explicit
isomorphism between the parabolic Hecke correspondence given by $g_j$ with the corresponding
piece of the Hecke correspondence given by $g$:
\begin{equation*}\begin{CD}
\Gamma_P[g_j]\backslash D[P] &\ \rightrightarrows\ & \Gamma_P \backslash D[P] \\
@V{\beta_j}VV @VV{\beta}V \\
\Gamma[g_j]\backslash \overline{D} &\ \rightrightarrows\ &\Gamma \backslash \overline{D} \\
@V{f}VV @| \\
\Gamma[g] \backslash \overline{D} &\ \rightrightarrows\ & \Gamma \backslash \overline{D}
\end{CD}\end{equation*}
The first line is the parabolic Hecke correspondence (\ref{eqn-Pcorrespondence}) defined by
$g_j$, i.e., it is the continuous extension of the mapping $\Gamma_P[g_j] \mapsto
(\Gamma_PxK'_P, \Gamma_P g_jxK'_P).$  The second line is the Hecke correspondence
(\ref{eqn-Hecke}) defined by $g_j$, i.e. it is the continuous extension of $\Gamma[g_j]xK'
\mapsto (\Gamma xK', \Gamma gxK').$  The vertical mapping $\beta$  (resp. $\beta_j$) is
described in \S \ref{subsec-parabnbhd}; it is a homeomorphism over any
$\Gamma$-parabolic neighborhood of $X_P$ (resp. over any $\Gamma[g_j]-$parabolic neighborhood
of $\Gamma_P[g_j] \backslash D_P$).  The top square of this diagram commutes
by direct computation.  The third line is the given Hecke correspondence (\ref{eqn-Hecke}).
The vertical mapping $f$ is the isomorphism of Hecke correspondences given in lemma
\ref{lem-isomorphism} and equation (\ref{eqn-isomorphism}).  In other words, write $g_j =
\gamma_1g\gamma_2$ and $f(\Gamma[g_j]xK') = \Gamma[g]\gamma_2xK'.$  The bottom square
also commutes.  This completes the construction of an explicit isomorphism with the parabolic
Hecke correspondence, and hence of the proof of proposition \ref{prop-restriction}.  \qed

\section{Counting the fixed points}\label{sec-counting}
\subsection{}\label{subsec-fpboundary}
As in \S \ref{sec-notation}, $\mathbf G$ denotes a connected linear reductive algebraic group
defined over $\mathbb Q$, $D$ denotes the associated symmetric space,  $K' = A_GK(x_0)$ is the
stabilizer in $G$ of a fixed basepoint $x_0\in D$, $\Gamma \subset \mathbf G(\mathbb Q)$
is an arithmetic group and $X = \Gamma \backslash D.$ 
Fix a rational parabolic subgroup $\mathbf P \subseteq \mathbf G$, set $\Gamma_P
= \Gamma \cap P$, $\Gamma_L = \nu_P(\Gamma_P) \subset \mathbf{L}(\mathbb Q)$ and
denote by $X_P = \Gamma_P \backslash D_P = \Gamma_L \backslash D_P$ the corresponding
stratum in the reductive Borel-Serre compactification $\overline{X}.$  An
element $y \in \mathbf{P}(\mathbb Q)$ gives rise to a parabolic Hecke correspondence
$(c_1,c_2):\Gamma'_P \backslash D[P] \rightrightarrows \Gamma_P \backslash D[P]$
where $\Gamma'_P = \Gamma_P[y] = \Gamma_P \cap y^{-1}\Gamma_P y. $ Let $\Gamma'_L =
\nu_P(\Gamma'_P).$  The restriction $C_P \rightrightarrows X_P$ of this parabolic
correspondence to the boundary stratum $C_P = \Gamma'_L \backslash D_P$ is given by
\begin{equation}\label{eqn-boundarycorrespondence}
 (c_1,c_2)(\Gamma'_L xK_PA_P) = (\Gamma_L x K_PA_P, \Gamma_L \bar y x K_PA_P)\end{equation}
where $\bar y = \nu_P(y).$

\subsection{Characteristic element}\label{subsec-char-elt}
Let us suppose that $w\in C_P$ is a {\it fixed point} of the parabolic
Hecke correspondence, that is, $c_1(w) = c_2(w).$  Choose any lift
$\tilde w\in D_P$ of $w$ and write $\tilde w =  z K_PA_P
\in L_P/K_PA_P.$  Since the  point $w=\Gamma'_L z K_PA_P$ is fixed, we have
\begin{equation}\label{eqn-e1} 
\Gamma_L  z K_PA_P = \Gamma_L \bar y  z K_PA_P.
\end{equation}
Hence there exists $\gamma \in \Gamma_L$ such that  the element $e = \gamma \bar y$ 
fixes the point $z\in D_P$, that is, 
\begin{equation}\label{eqn-e2}
e  z K_PA_P=\gamma \bar y  z K_PA_P =  zK_PA_P\in D_P. \end{equation}
\begin{defn}
The element $e = \gamma \bar y \in \mathbf{L_P}(\mathbb Q)$ is called
{\it a characteristic element} for the fixed point $w$, or {\it the
characteristic element} corresponding to the lift $\tilde w$ of $w$.
\end{defn}

Denote by $F_P(e) \subset C_P$ the set of fixed points in $C_P$ for which $e$ is a
characteristic element.  We refer to $F_P(e)$ as a {\it fixed point constituent}; it may
consist of several connected components.  Let $D_P^e$ denote the fixed points of the mapping
$T_e : D_P \to D_P.$  Then $F_P(e)$ is the image of $D_P^e$ under the projection $D_P \to
C_P.$  Write $e=a_em_e\in A_PM_P.$  We say an element  of $L_P$ is {\it elliptic} (or is
elliptic modulo $A_P$) if it is $\mathbf{L_P}(\mathbb R)$-conjugate to an element of $K_PA_P.$
Let $L_e \subset L_P$ denote the centralizer of $e$ in $L_P.$

\begin{prop}\label{prop-char} Suppose the arithmetic group $\Gamma$ is neat.  Then the
following statements hold.
\begin{enumerate}
\item  The characteristic element $e \in L_P$ is semisimple and is elliptic
\parens{modulo $A_P$}.  The group $L_e$ is reductive, algebraic, and defined over $\mathbb Q.$
The torus factors $a_y = a_e\in A_P$ are equal \parens{cf \S \ref{subsec-parabolics}}.
The fixed point constituent $F_P(e)$ is a smooth submanifold of $C_P.$
\item If $y\in P$ is changed by multiplication by an element $u\in \mathcal U_P$, or if
$\gamma \in \Gamma_L$ is replaced by another element of $\Gamma_L$ which also satisfies
\parens{\ref{eqn-e2}}, or if a different representative $ z' \in L_P$ of $\tilde w\in D_P$ is
chosen, then the characteristic element $e\in L_P$ does not change.
\item If a different lift $\tilde w'\in D_P$ of $w$ is chosen, or if $y$ is changed within its
double coset $\Gamma_Py\Gamma_P$ then $e$ changes at most by $\Gamma_L$-conjugacy.
\item The characteristic element $e$ is a rigid invariant of the fixed point set:  if $w_t
\in C_P$ is a one parameter family of fixed points \parens{with $t\in [0,1]$} and if $z_t \in
D_P$ is a lift to a one parameter family of points in $D_P$ then the resulting characteristic
elements $e_t$ do not vary with $t$.
\item 
The group $L_e$ acts transitively on $D_P^e.$  Set $\Gamma_e = \Gamma_L \cap L_e$, $\Gamma'_e
= \Gamma'_L \cap L_e$, and $K'_e = L_e \cap ( z (K_PA_P) z^{-1}).$ Then the action
of $L_e$ on $D_P^e$ induces diffeomorphisms
\begin{equation}\label{eqn-centralizer}
F_P(e) \cong \Gamma'_e \backslash L_e / K'_e \text{ and } c_i(F_P(e)) \cong \Gamma_e
\backslash L_e /K'_e.
\end{equation}
The projection $F_P(e) \to c_i(F_P(e))$ is a covering of degree
\begin{equation}\label{eqn-degree-d}
d = [ \Gamma_e : \Gamma'_e] = [ \Gamma_L \cap \bar y^{-1} \Gamma_L \bar y :
\nu_P(\Gamma_P \cap y^{-1} \Gamma_P y)].  \end{equation}
\item Conversely, let  $\mathcal C \subset \Gamma_L \bar y \Gamma_L$ be any
$\Gamma_L$-conjugacy class which is elliptic (modulo $A_P$). Then $\mathcal C \cap \Gamma_L
\bar y$ consists of a single element $e'$, and there exists a fixed point $w' \in C_P$ for
which $e'$ is a characteristic element.  In particular, $F_P(e') \ne \phi.$
\end{enumerate}
\end{prop}
So the constituents of the fixed point set in $C_P$ are in one to one correspondence
with $\Gamma_L$-conjugacy classes of elliptic (modulo $A_P$) elements $e\in \Gamma_L \bar y
\Gamma_L$.

\subsection{Proof}\label{subsec-pf-prop-char}
By (\ref{eqn-e2}), $e = \gamma \bar y \in  z K_PA_P z^{-1}$ which is compact
modulo $A_P$, so $e$ is also semisimple.  The elements $e$ and $\bar y$ are in the same
$\Gamma_L$-double coset so they have the same torus component $a_e = a_y \in A_P.$  Since
$\Gamma$ is neat, the group $\Gamma_e$ acts freely on $D_P^e.$  This
proves (1).  Next, consider (2) and suppose different choices $y'=yu$, $\gamma' \in \Gamma_L$,
and $ z'\in L_P$  were made, with $u \in \mathcal U_P$, with $\tilde w =  z K_PA_P =
 z'K_PA_P \in D_P$ and, as in (\ref{eqn-e2}), $ z'K_PA_P = \gamma' \bar y' 
z'K_PA_P$ (where $\bar y' = \nu_P(y')$).  Then $\bar y' = \bar y$ and  
\[  \gamma \bar y  z K_PA_P =  z K_PA_P =\gamma ' \bar y' z' K_PA_P
= \gamma' \bar y z K_P A_P \]
so $\gamma ^{-1} \gamma '\in (\bar y  z) K_PA_P(\bar y  z)^{-1}$.
Since $\Gamma_L$ is torsion-free, this implies $\gamma = \gamma'$, hence the characteristic
element $e = \gamma \bar y$ is unchanged.
This proves (2).  Since $\Gamma_L$ is discrete, the characteristic element is constant in a
continuous family of fixed points, which proves (4).

Now consider changing $y$ within its double coset $\Gamma_P y \Gamma_P$ and consider changing
the lift $\tilde w \in D_P$ of the fixed point.  Let $\hat y = \gamma_1y\gamma_2$ with
$\gamma_1,\gamma_2 \in \Gamma_P.$  Set $\bar\gamma_1 = \nu_P(\gamma_1)$, $\bar\gamma_2 =
\nu_P(\gamma_2)$ and $\bar{\hat y} = \nu_P(\hat y).$   As
in \S \ref{subsec-fpboundary}, the element $\bar{\hat y}$ determines a Hecke correspondence
$(\hat c_1,\hat c_2):\widehat{C}_P \rightrightarrows X_P$ as follows:  Set $\widehat{\Gamma}_P
=\Gamma_P \cap \hat y^{-1} \Gamma_P \hat y$, $\widehat{\Gamma}_L = \nu_P(\widehat{\Gamma}_P)$
and $\widehat{C}_P = \widehat{\Gamma}_L\backslash D_P = \widehat{\Gamma}_L \backslash
L_P/K_PA_P$.  Then $(\hat c_1,\hat c_2)(\widehat{\Gamma}_LxK_PA_P) = (\Gamma_LxK_PA_P,
\Gamma_L \bar{\hat y} x K_PA_P).$  As in equation (\ref{eqn-isomorphism}), an isomorphism of
correspondences $f:\widehat{C}_P \to C_P$ is given by
\begin{equation*}
f(\widehat{\Gamma}_LxK_PA_P) = \Gamma'_L\bar\gamma_2xK_PA_P.
\end{equation*}
Choose any lift ${\hat z}K_PA_P \in D_P$ of the fixed point $\hat w = f^{-1}(w)$ (with
${\hat z} \in L_P$).  We obtain a new characteristic element $\hat e = \hat\gamma
\bar{\hat y}$ (for some $\hat\gamma \in \Gamma_L$) such that
\begin{equation}\label{eqn-eprime2}
\hat\gamma\bar{\hat y} {\hat z} K_PA_P = {\hat z} K_P A_P.
\end{equation}
We need to show that $\hat e=\hat\gamma \bar{\hat y}$ is $\Gamma_L$-conjugate to $e = \gamma
\bar y.$ Since $f(\hat w)=w$ we have
\begin{equation*}%\label{eqn-eprime3}
\Gamma'_L\bar \gamma_2 {\hat z} K_PA_P = \Gamma'_L z K_PA_P
\end{equation*}
so there exists $h \in \Gamma'_L$ such that
$h\bar \gamma_2 {\hat z} K_PA_P =  z K_PA_P $ or
\begin{equation}\label{eqn-eprime4}
{\hat z}K_PA_P = \bar\gamma_2^{-1}h^{-1} z K_PA_P.
\end{equation}
Substituting (\ref{eqn-eprime4}) into both sides of (\ref{eqn-eprime2}) and using
(\ref{eqn-e2}) gives
\begin{equation*}%\label{eqn-eprime5}
h \bar \gamma_2 \hat\gamma \bar{\hat y} \bar \gamma_2^{-1}h^{-1}  z K_PA_P =  z
K_PA_P= \gamma \bar y  z K_PA_P
\end{equation*}
or
\begin{equation*}%\label{eqn-eprime6}
h \bar\gamma_2 \hat\gamma \bar\gamma_1 \bar y h^{-1}  z K_P A_P = \gamma \bar y  z
K_PA_P.
\end{equation*}
Since $h,h^{-1} \in \nu_P(y^{-1}\Gamma_Py) =\bar y^{-1}\Gamma_L \bar y$ there exists $h' \in
\Gamma_L$ such that $\bar y
h^{-1} = h' \bar y,$ which gives
\begin{equation*}%\label{eqn-eprime7}
\gamma^{-1}h \bar\gamma_2 \hat\gamma \bar\gamma_1 h' (\bar y  z K_PA_P) = \bar y 
zK_PA_P.
\end{equation*}
This implies that $\gamma^{-1}h\bar \gamma_2\hat\gamma \bar \gamma_1 h' = 1$ since it is both
in the group $(\bar y  z)K_PA_P (\bar y  z)^{-1}$ and in $\Gamma_L.$  Therefore
$1.\bar y = \gamma^{-1}h \bar\gamma_2 \hat\gamma \bar\gamma_1 \bar y h^{-1}$ or
\begin{equation*}%\label{eqn-eprime8}
\gamma \bar y = h \bar\gamma_2 \hat\gamma \bar\gamma_1 \bar y \bar \gamma_2 \bar\gamma_2^{-1}
h^{-1} = (h\bar\gamma_2) \hat\gamma \bar{\hat y} (h\bar\gamma_2)^{-1}.
\end{equation*}
Thus, the characteristic elements $\gamma \bar y$ and $\hat\gamma \bar{\hat y}$ are
$\Gamma_L$-conjugate, which proves (3).  

Now let us prove (5).  It is easy to see that $L_e$ acts on $D_P^e.$  To see that this action
is transitive, let $v_1,v_2\in D_P^e$, say $v_1 =  z_1 K_PA_P$ and $v_2 =  z_2
K_PA_P$ with $ z_i \in L_P$ (for $i = 1,2$).  Then there exists $k_1,k_2 \in K_PA_P$ so
that $e z_1 =  z_1 k_1$ and $e z_2 =  z_2 k_2$, hence $k_1$ and $k_2$
are $L$-conjugate (by $ z_2 z_1^{-1}$).  It follows from \cite{Borel} \S 24.7
that $k_1$ and $k_2$ are also $K_PA_P$ conjugate.  Say, $k_2 = mk_1m^{-1}$ for some $m \in
K_PA_P.$  Define $x=  z_2m  z_1^{-1}$.  Then $v_2 = xv_1$ and moreover, $x\in
L_e$ since
\begin{equation*}
xex^{-1} =  z_2m z_1^{-1}e z_1 m^{-1}  z_2^{-1} =
 z_2 m k_1m^{-1} z_2^{-1}= e.
\end{equation*}
This completes the verification that $L_e$ acts transitively on $D_P^e$.  

Using the chosen lift $\tilde w = zK_PA_P \in D_P^e$ as a basepoint, we obtain a
diffeomorphism $L_e /K'_e\cong D_P^e$ where  $K'_e = L_e \cap(z K_PA_P z^{-1})$ is
the stabilizer (in $L_e$) of $\tilde w.$  This induces a surjection $(L_e \cap \Gamma'_L
)\backslash L_e / K'_e \to F_P(e)$ which we will now show to be injective.  

Suppose $x_1,x_2\in L_e$ and that $x_1\tilde w, x_2\tilde w \in D_P^e$ map to the same point
in $C_P$, that is, $\Gamma'_L x_1 z K_PA_P =\Gamma'_L x_2 z K_PA_P.$  Then there exists
$\gamma\in \Gamma'_L$ so that 
\begin{equation}\label{eqn-tildez}
\gamma x_1  z K_PA_P = x_2  z K_PA_P.
\end{equation}  
We need to show that $\gamma \in L_e.$  Acting by $e$ on the left hand side of
(\ref{eqn-tildez}) and using (\ref{eqn-e2}) gives the quantity
\begin{equation*}
e\gamma e^{-1}e x_1  z K_PA_P = e\gamma e^{-1} x_1e z K_PA_P = e\gamma e^{-1}
x_1 z K_PA_P \end{equation*}
while acting by $e$ on the right hand side of (\ref{eqn-tildez}) gives 
\begin{equation*}
e x_2 z K_PA_P = x_2  z K_PA_P = \gamma x_1 z K_PA_P. \end{equation*}
So $\gamma^{-1} e \gamma e^{-1} \in (x_1 z)K_PA_P (x_1 z)^{-1}$.  But
$\gamma^{-1} e \gamma e^{-1} \in \Gamma_L$ so this element is trivial, that is, $\gamma e = e
\gamma,$ hence $\gamma \in L_e.$  Therefore $F_P(e) = \Gamma'_e \backslash D_P^e.$  The
equality $c_1(F_P(e)) = c_2(F_P(e)) = \Gamma_e \backslash L_e /K'_e$ is similar.
Equation (\ref{eqn-degree-d}) will be proven in \S \ref{subsec-addendum}.  
\quash{
From part (5) we have:  $F_P(e) = \Gamma'_e \backslash D_P^e,$ and its image under
$c_1$ (which equals its image under $c_2$) is $\Gamma_e \backslash D_P^e.$  So part (6) will
follow once we know that the inclusion $\Gamma'_e \subset \Gamma_e$ is an equality.  If
$\gamma \in \Gamma_e$ then $\gamma = e^{-1} \gamma e.$  But the left side of this equation is
in $\Gamma_L$ and the right side is in $\bar y ^{-1} \Gamma_L \bar y$, hence $\gamma \in
\Gamma'_L$ as claimed.  } %% end quash

Now let us verify part (6).  Suppose $e'' = \gamma_2 \bar y \gamma_1 \in \mathbf{L_P}(\mathbb
Q)$ is elliptic (modulo $A_P$).  Then $e''$ is $\Gamma_L$-conjugate to the element
$e'=\gamma_1\gamma_2\bar y$ which is also elliptic modulo $A_P.$  There exists
$ z \in L_P$ so that $e' \in z K_PA_P  z^{-1}$ hence $e' z K_PA_P =
 z K_PA_P.$  In other words, $e'$ is a characteristic element for the point
$w'=\Gamma'_L z K_PA_P \in C_P$ (which is easily seen to be fixed under the Hecke
correspondence).  \qed

\subsection{Addendum}\label{subsec-addendum}
  (We refer to the notation of \S \ref{prop-char} and \S \ref{subsec-pf-prop-char}.)
Unfortunately the arithmetic group $\Gamma_L[\bar y] = \Gamma_L \cap \bar y^{-1} \Gamma_L \bar
y$ may be larger than $\Gamma'_L= \nu_P(\Gamma_P \cap y^{-1} \Gamma_P y)$, so the
correspondence (\ref{eqn-boundarycorrespondence}) is not necessarily a {\it Hecke}
correspondence, but rather it is a covering of the following Hecke correspondence:
\begin{align*} (\tilde c_1, \tilde c_2):\widetilde{C}_P = \Gamma_L[\bar y] \backslash D_P
&\rightrightarrows X_P \\
\Gamma_L[\bar y] x K_PA_P &\mapsto (\Gamma_L x K_PA_P, \Gamma_L \bar y x K_PA_P).\end{align*}
This covering $\phi:\Gamma'_L \backslash D_P \to \Gamma_L[\bar y] \backslash D_P$ has degree
$d = [\Gamma_L[\bar y] : \Gamma'_L].$  A point $w \in C_P$ is fixed iff the point $\phi(w) \in
\widetilde{C}_P$ is fixed.

Let $\widetilde{F}_P(e) \subset \widetilde C_P$ be the set of fixed points within this
(smaller) Hecke correspondence with characteristic element $e.$  We claim that the restriction
of $\tilde c_i$ to $\widetilde {F}_P(e) \subset \widetilde{C}_P$ is a diffeomorphism:
$\widetilde{F}_P(e) \cong c_i(F_P(e))$ (for $i = 1,2$).  As in (\ref{eqn-centralizer}) it is
clear that $\widetilde{F}_P(e) \cong (\Gamma_L[\bar y] \cap L_e) \backslash L_e /K'_e$ so it
suffices to verify that the inclusion $\Gamma_L[\bar y] \cap L_e \subset \Gamma_e$ is an
isomorphism.  If $\gamma_1 \in \Gamma_e$ then $\gamma_1 = e^{-1} \gamma_1 e.$  But the left
side of this equation is in $\Gamma_L$ and the right side is in $\bar y^{-1} \Gamma_L \bar y$,
hence $\gamma_1 \in \Gamma_L[\bar y]$, which proves the claim.

In summary, we have a diagram
\begin{equation*} \begin{CD}
C_P @>{\phi}>> \widetilde{C}_P & \rightrightarrows & X_P\\
\bigcup && \bigcup && \bigcup \\
F_P(e) @>>> \widetilde{F}_P(e) @>{\cong}>> c_i(F_P(e))
\end{CD}\end{equation*}
and in particular the degree $d$ of the covering $\phi$ coincides with the degree of the
mapping $F_P(e) \to c_i(F_P(e))$ which gives (\ref{eqn-degree-d}).

\subsection{Remark}\label{subsec-orientation}  The codimension of $F_p(e)$ in $D_P$ is odd if
and only if the action of
$e$ reverses orientations in the normal bundle of $F_P(e)$ in $D_P.$
This is because $e$ preserves an appropriately chosen normal slice through any point in
$F_P(e)$, the boundary of which is a  sphere on which $e$ then acts as a
diffeomorphism without fixed points, so its Lefschetz number is 0.  In the odd codimension
case this sphere is even dimensional, so by the Lefschetz fixed point theorem, its action on
the top degree cohomology is given by multiplication by $-1$, that is, it reverses the
orientation.  In the even codimension case, it preserves orientation.

%%%%%%%%%%%%%%%%%%%%%%%%%%%%%%%%%%%%%%%%%%%%%%%%%%%%%%%%%%%%%%%%%%%%%

\section{Hyperbolic properties of Hecke correspondences}\label{sec-hyperbolic}
\subsection{Expanding and contracting roots}\label{subsec-hyperbolic}  As in \S
\ref{sec-notation}, $\mathbf G$ denotes a connected linear reductive algebraic group defined
over $\mathbb Q$, $D$ denotes the associated symmetric space,  $K' = A_GK(x_0)$ is the
stabilizer in $G$ of a fixed basepoint $x_0\in D$, $\Gamma \subset \mathbf G(\mathbb Q)$
is an arithmetic group and $X = \Gamma \backslash D.$  Throughout this section we fix a Hecke
correspondence $(c_1,c_2):\overline{C}\rightrightarrows\overline{X}$ defined by some element
$g\in \mathbf G(\mathbb Q).$  So $C = \Gamma'\backslash D$ with $\Gamma' = \Gamma \cap
g^{-1}\Gamma g.$  

  Let $\mathbf P \subset \mathbf G$ be a rational parabolic subgroup and suppose that
$c_1(C_P) = c_2(C_P) = X_P.$   By Proposition \ref{prop-restriction},  near $C_P$ the
correspondence is modeled on a parabolic Hecke correspondence $\Gamma'_P \backslash D[P]
\rightrightarrows \Gamma_P \backslash D[P]$ (\ref{eqn-Pcorrespondence}) which is determined by
some $y\in \mathbf P(\mathbb Q)$ (where $\Gamma'_P = \Gamma_P \cap y^{-1} \Gamma_P y$). 
Suppose $y = u_y a_y m_y$ is the Langlands decomposition (\ref{eqn-Langlandselement}) of $y\in
\mathbf P(\mathbb Q).$  If $y$ is allowed to vary within the double coset $\Gamma_P y
\Gamma_P$ then the element $a_y \in A_P$ will remain fixed, so we may write $a_P = a_y.$  We
refer to $a_P$ as the {\it torus factor} associated to the Hecke correspondence near $C_P.$
The torus factor may be used
to define a partition of the simple roots $\Delta_P$ into three subsets
\begin{align*}
\Delta^{+}_P &= \{ \alpha \in \Delta_P|\ \alpha(a_P) \alt 1 \} \\
\Delta^{-}_P &= \{ \alpha \in \Delta_P|\ \alpha(a_P) \agt 1 \} \\
\Delta^{0}_P &= \{ \alpha \in \Delta_P|\ \alpha(a_P) = 1 \}
\end{align*}
consisting of those simple roots which are {\it expanding, contracting}, or {\it neutral},
respectively, near the stratum $C_P.$  (See also \S \ref{subsec-hypercorres})
The terminology is motivated by the following fact, whose proof follows immediately from
(\ref{eqn-rootfunctionP}) and the definition (\ref{eqn-correspondence}) of the correspondence:
for all $\alpha \in \Delta_P$ and for all $z \in \Gamma'_P \backslash D[P]$ we have:
\begin{equation}\label{eqn-monitor}
f_{\alpha}^P(c_2(z)) = \alpha(a_P)f_{\alpha}^P(c_1(z)).  
\end{equation}

Now suppose $\mathbf P \subset \mathbf Q$ are rational parabolic subgroups of $\mathbf G$,
with $\Delta_P = i(\Delta_Q) \amalg J$ as in (\ref{eqn-decomposition}).  Suppose that
$c_1(C_P) = c_2(C_P),$ giving rise to a torus factor $a_P \in A_P$ and a decomposition of
$\Delta_P$ into expanding, contracting and neutral roots as above.  Then $c_1(C_Q) = c_2(C_Q)$
(by \S \ref{subsec-restriction}) so we obtain a torus factor $a_Q \in A_Q$ and a decomposition
of $\Delta_Q$ into expanding, contracting and neutral roots also. 
\begin{prop}\label{prop-neutral} Suppose that $J \subset \Delta_P^0.$   Then \begin{enumerate}
\item The torus factors $a_P=a_Q$ are equal; in particular $a_P$ lies in the sub-torus $A_Q
\subset A_P.$
\item The expanding, contracting, and neutral simple roots for $\mathbf P$ are given by:
\begin{equation}\label{eqn-decomp2}
\Delta_P^{+} = i(\Delta_Q^{+}),\ \Delta_P^{-} = i(\Delta_Q^{-}),\ \Delta_P^0 = i(\Delta_Q^0)
\amalg J.
\end{equation}
\item  For all $z \in \Gamma'_P \backslash D[P]$ and for all $\beta\in\Delta_Q$ we have,
\begin{equation*}
\frac{f_{i(\beta)}^P(c_2(z))}{f_{i(\beta)}^P(c_1(z))} = i(\beta)(a_P)=\beta(a_P)=
\frac{f_{\beta}^Q(c_2(z))}{f_{\beta}^Q(c_1(z))} 
\end{equation*}
provided the denominators do not vanish.
\end{enumerate}\end{prop}
\noindent
In this case we say that $\mathbf Q$ is a {\it neutral} parabolic subgroup containing $\mathbf
P$ and we write $\mathbf P \prec \mathbf Q.$  Intuitively, the Hecke correspondence is neutral
in those directions normal to $C_P$ which point into $C_Q$; cf. \S \ref{subsec-intuition}.

\subsection{Proof} Locally near $C_Q$ the Hecke correspondence is isomorphic to a parabolic
Hecke correspondence given by some $y'=u_{y'}a_{y'}m_{y'} \in Q$ with torus factor $a_Q =
a_{y'}\in A_Q.$  In a neighborhood of $C_P$ the correspondence is isomorphic to the parabolic
Hecke correspondence given by some $y=u_ya_ym_y \in P$ (with $a_P = a_y$).  Moreover $y$ may
be chosen to lie in the double coset $\Gamma_Q y' \Gamma_Q$ since the correspondence
$\overline{C}_P \rightrightarrows \overline{X}_P$ is the restriction to $\overline{C}_P$ of
the correspondence $\overline{C}_Q \rightrightarrows \overline{X}_Q;$ cf. Proposition
\ref{prop-restriction}.  By assumption, $\alpha(a_P)=1$ for all $\alpha \in J$ which implies
that $a_P \in A_Q.$  It follows that $a_P= a_Q$ because the homomorphism $Q \to A_Q$ (which
associates to any $z \in Q$ its torus factor $a_z$) kills $\Gamma_Q.$  Therefore, for any
$\beta \in \Delta_Q$ we have:  $\beta(a_Q) = i(\beta)(a_P).$  This proves (1) and (2).  The
first equality in part (3) is just (\ref{eqn-monitor}).  The last equality in part (3) follows
from part (1) and from (\ref{eqn-monitor}) (with $f_{\alpha}^P$ replaced by $f_{\beta}^Q$).
\qed

\subsection{Maximal neutrality}\label{subsec-dagger}
  Suppose $\mathbf P\subset \mathbf Q \subset \mathbf R$ are
rational parabolic subgroups and that $c_1(C_P) = c_2(C_P).$  Write $\Delta_P = i(\Delta_Q)
\amalg I$ and $\Delta_Q = j(\Delta_R) \amalg J$ for the disjoint union of
(\ref{eqn-decomposition}).  Suppose moreover that $\mathbf P \prec \mathbf Q$ and $\mathbf Q
\prec \mathbf R,$ that is, that $I \subset \Delta_P^0$ and $J \subset \Delta_Q^0.$  Then it
follows from Proposition \ref{prop-neutral} that $\mathbf P \prec \mathbf R.$  Hence there is
a greatest neutral parabolic subgroup $\mathbf P^{\dagger}$ containing $\mathbf P$; in fact it
is $\mathbf P^{\dagger} = \mathbf P(\Delta_P^0)$ in the notation of \S \ref{subsec-twoparab}
and \S \ref{subsec-roots}. It is easy to see that
\begin{equation}\label{eqn-noneutrals}
\Delta_{P^{\dagger}}^0 = \phi \text{ and } \Delta_P^{\pm} = i'(\Delta_{P^{\dagger}}^{\pm})
\end{equation}
(where $i': \Delta_{P^{\dagger}} \hookrightarrow \Delta_P$ is the natural inclusion).
Moreover, \begin{equation}\label{eqn-transitive}
\mathbf P \prec \mathbf Q \implies \mathbf P^{\dagger} = \mathbf Q^{\dagger}.  \end{equation}

\section{Structure of the fixed point set}
\subsection{}  As in \S \ref{sec-notation}, $\mathbf G$ denotes a connected reductive linear
algebraic group defined over $\mathbb Q$, $D = G/K'$ is its associated symmetric space with
basepoint $x_0 \in D$ and stabilizer $K' = A_GK(x_0)$, $\Gamma \subset \mathbf G(\mathbb Q)$
denotes an arithmetic subgroup, and $X = \Gamma \backslash D$.  Throughout this section we fix
a Hecke correspondence $(c_1,c_2):\overline{C}
\rightrightarrows
\overline{X}$ defined by some element $g\in \mathbf G(\mathbb Q).$  So $C = \Gamma'\backslash
D$ with $\Gamma' = \Gamma \cap g^{-1}\Gamma g.$  We also fix a $\Gamma$-equivariant tiling of
$D$ which is narrow with respect to the Hecke correspondence, and denote by $\{C^P\}$
and $\{X^P\}$ the resulting tilings of $\overline{C}$ and $\overline{X}$ respectively.
Let $F\subset \overline{C}$ be a connected component of the set of fixed points.  The
following lemma says that if $F$ spans two strata $C_P \subset \overline{C}_Q$ then the Hecke
correspondence is neutral in those directions which point from $C_P$ into $C_Q$; cf. \S
\ref{subsec-intuition}.

\begin{lem}\label{lem-geodesicinvariant}
Let $\mathbf P \subset \mathbf Q\subseteq \mathbf G$ be rational parabolic subgroups and write
$\Delta_P =i(\Delta_Q) \amalg J$ as in \parens{\ref{eqn-decomposition}}.  Suppose $F \cap
C_Q\cap T(\overline{C}_P) \ne \phi$ \parens{that is, $C_Q$ contains fixed points which lie in
the $\Gamma$-parabolic neighborhood $T(\overline{C}_P)$ of $\overline{C}_P$}.  Then
\begin{enumerate}
\item $J \subset \Delta_P^0$ \parens{hence the conclusions of Proposition \ref{prop-neutral}
hold}.
\item $F \cap T(\overline{C}_P)$ is invariant under the geodesic action of $A'_P(\age 1).$
\item $\pi_P(F \cap T(\overline{C}_P)) \subset F$, that is, each fixed point in this
$\Gamma'$-parabolic neighborhood projects to a fixed point in $C_P.$
\end{enumerate}\end{lem}

\subsection{Proof} Part (1) follows from (\ref{eqn-rootpi2}) and (\ref{eqn-monitor}) by taking
$z\in F \cap C_Q \cap T(\overline{C}_P)$ to be a fixed point. Part (2) follows from
(\ref{eqn-commutes geodesic}).  Part (3) follows by continuity.  \qed

\begin{prop}\label{prop-structure}  Let $F \subset \overline{C}$ be a connected component of
the fixed point set.
Let $\mathbf Q$ be a rational parabolic subgroup and suppose $F \cap C_Q \ne \phi.$  Let $a_Q
\in A_Q$ be the torus factor for this stratum.  Let $\mathbf Q^{\dagger} = \mathbf
Q(\Delta_Q^0)$ be the maximal neutral parabolic subgroup containing $\mathbf Q.$  Then
\begin{enumerate}
 \item The whole connected component $F$ of the fixed point set is contained in the closure
\begin{equation*}
F \subset \overline{C}_{Q^{\dagger}}
\end{equation*}
of the single stratum $C_{Q^{\dagger}}.$  
\item The Hecke correspondence $\overline{C}\rightrightarrows \overline{X}$ restricts to a
correspondence $\overline{C}_{Q^{\dagger}} \rightrightarrows \overline{X}_{Q^{\dagger}}$
on this stratum-closure.  Within this restricted correspondence, near each point
$c\in F,$ every simple root is neutral: If $\mathbf P \subset \mathbf Q^{\dagger}$, if $F \cap
C_P \ne \phi$ if $i':\Delta_{Q^{\dagger}} \hookrightarrow \Delta_P$ is the inclusion, and if
$\Delta_P= i'(\Delta_{Q^{\dagger}}) \amalg J$ as in \parens{\ref{eqn-decomposition}} then $J =
\Delta_P^0.$ 
\item  There exists a neighborhood $U(F) \subset \overline{C}$ of the fixed point set such
that for all $\alpha\in \Delta_Q$ and for all $x \in U(F)$, if $x \in T(\overline{C}_Q)$ and
if $c_2(x) \in T(\overline{X}_Q)$ then 
\begin{equation}\label{eqn-exp-contr}
r_{\alpha}^Q(c_2(x)) = \alpha(a_Q)^{-1} r_{\alpha}^Q(c_1(x)).
\end{equation}
\end{enumerate}\end{prop}

\subsection{Remarks}  Part (1) does {\it not} imply that $F \cap C_{Q^{\dagger}} \ne \phi.$ In
fact, the fixed point component $F$ may be ``reducible'': it does not necessarily coincide
with the closure of its intersection $ F \cap C_P$ with any single stratum $C_P.$  (See \S
\ref{subsec-reducible}.)

\subsection{Proof}  Suppose that $F$ has a nontrivial intersection with some
other stratum, say $F\cap C_R \ne \phi.$  Suppose for the moment that $\mathbf R \supset
\mathbf Q$ and that  $F \cap C_Q$ contains limit points from $F \cap C_R$, that is,
\begin{equation}\label{eqn-limitpoint}
(F \cap C_Q) \cap \overline{F \cap C_R} \ne \phi.
\end{equation}
Then Lemma \ref{lem-geodesicinvariant} part (1) implies that $\mathbf R$ is a neutral
parabolic subgroup containing $\mathbf Q$ so (\ref{eqn-transitive}) implies that $\mathbf
R^{\dagger} = \mathbf Q^{\dagger},$ hence $F \cap C_R \subset F \cap
\overline{C}_{R^{\dagger}} = F \cap \overline{C}_{Q^{\dagger}}.$ Now we drop the assumption
(\ref{eqn-limitpoint}).  Since $F$ is connected, the stratum $C_R$ is related to the stratum
$C_Q$ through a chain of strata $C_{R_i}$ (say, $1 \le i \le m$), each having nontrivial
intersection with $F$, with each step in the chain related to the next by
\[(F \cap C_{R_i}) \cap \overline{F \cap C_{R_{i+1}}} \ne \phi \text{ or }
(F \cap C_{R_{i+1}}) \cap \overline{F \cap C_{R_i}} \ne \phi. \]  Repeated
application of (\ref{eqn-transitive}) implies that 
\begin{equation}\label{eqn-chain}
\mathbf R^{\dagger} = \mathbf R_1^{\dagger} = \cdots = \mathbf R_m^{\dagger}=\mathbf
Q^{\dagger}.\end{equation}
So once again, $F \cap C_R \subset F \cap \overline{C}_{Q^{\dagger}}.$ This verifies part (1).

Consider part (2).
Since the stratum $C_Q$ is preserved by the Hecke correspondence, the same holds for each
larger stratum, especially $C_{Q^{\dagger}}.$  Suppose $F \cap C_P \ne \phi.$  By
(\ref{eqn-chain}), $\mathbf P^{\dagger} = \mathbf Q^{\dagger}$ so by (\ref{eqn-noneutrals}),
$\Delta_P = i'(\Delta_{Q^{\dagger}}) \amalg \Delta_P^0.$

Next we verify part (3). Suppose $\mathbf P \subset \mathbf Q$ and suppose $F \cap C_P \ne
\phi.$  By Proposition \ref{prop-neutral} part (3), for all $\alpha \in \Delta_Q$ and for all
$w \in T(\overline{C}_Q) \cap c_2^{-1}(T(\overline{X}_Q)),$
the root function $f_{i(\alpha)}^{P}$ satisfies 
\begin{equation}\label{eqn-more}
f_{i(\alpha)}^{P}(c_2(w)) = \alpha(a_Q)f_{i(\alpha)}^{P}(c_1(w))
\end{equation}
where, as in (\ref{eqn-decomposition}), we have written $\Delta_{P} = i(\Delta_Q) \amalg J.$
However this does not yet prove (\ref{eqn-exp-contr}).  The problem is that
the partial distance function $r_{\alpha}^Q(w)$ is patched together
(\ref{eqn-partiald}) from these root functions $f_{i(\alpha)}^P$ in a way that depends on
which tile contains the point $c_i(w)$.  So we need to show that the Hecke correspondence
preserves the tile boundaries in some neighborhood $U(F)$ of the fixed point set.

\begin{lem}\label{lem-neutral}  Suppose $F \subset \overline{C}$ is a connected component of
the fixed point set of the Hecke correspondence $\overline{C} \rightrightarrows \overline{X}.$
Then there exists a neighborhood $U(F) \subset \overline{C}$ of $F$ such that for all $w \in
U(F)$ and for any rational parabolic subgroup $\mathbf P \subseteq \mathbf G,$  
\begin{equation}\label{eqn-preservestiles}
c_1(w) \in X^P \iff c_2(w) \in X^P.
\end{equation}
\end{lem}

\subsection{Proof}  Assume not.  Then there is a sequence of points $x_i \in \overline{C}$
converging to $F$ so
that for each $i$, $c_1(x_i)$ and $c_2(x_i)$ are in different tiles.  By taking subsequences
if necessary we may assume the sequence $x_i$ converges to some point $x_0\in F$, that $x_i$
are all contained in a single tile $C^P$ (so $c_1(x_i) \in X^P$) and that $c_2(x_i)$ all lie
in a single tile $X^Q$.  Since $c_1(x_0) \in \overline{C}^P \cap \overline{C}^Q$ is a fixed
point, the Hecke correspondence must preserve the strata $C_P$ and $C_Q$ (meaning that
$c_1(C_P)=c_2(C_P)=X_P$ and $c_1(C_Q)=c_2(C_Q) = X_Q$) and  we may
assume that either $\mathbf P \subseteq \mathbf Q$ or $\mathbf Q \subseteq \mathbf P.$  Since
the tiling is narrow this implies that $F \cap C_P \ne \phi,$ that $F \cap C_Q \ne \phi,$ and
that either $F\cap C_P$ contains limit points from $F\cap C_Q$ (if $\mathbf P \subseteq
\mathbf Q$) or else $F \cap C_Q$ contains limit points from $F \cap C_P$ (if $\mathbf Q
\subseteq \mathbf P$).

Let us first consider the case that $\mathbf P \subseteq \mathbf Q = \mathbf G$, meaning that
$X^Q = X^0.$ Let $a_P \in A_P$ denote the torus factor for the Hecke correspondence near $C_P$
as in \S \ref{subsec-hyperbolic}.  Then it follows from Lemma \ref{lem-geodesicinvariant} that 
$\mathbf Q = \mathbf G$ is a neutral parabolic subgroup containing $\mathbf P$, that is,
$\alpha(a_P) = 1$ for all $\alpha \in \Delta_P.$  Within any $\Gamma$-parabolic neighborhood
$W$ of $\overline{X}_P$, the tile $X^P$ is given by (\ref{eqn-singletile}):
\[ X^P = \{ x \in W |\   \pi_P(x)\in X_P^0\text{ and } f_{\alpha}^P(x) \agt \alpha(b_P)
\text{ for all } \alpha \in \Delta_P \}. \]
Since $c_2(x_i) \in X^0$ it follows that for at least one $\alpha \in \Delta_P$ we have:
\[ \alpha(b_P) \age f_{\alpha}^P(c_2(x_i)) = \alpha(a_P)f_{\alpha}^P(c_1(x_i)) =
f_{\alpha}^P(c_1(x_i)) \agt \alpha(b_P) \]
(using equation (\ref{eqn-monitor})) which is a contradiction.

Next consider the case $\mathbf P \subset \mathbf Q \ne \mathbf G.$  For sufficiently large
$i$ the points $x_i$ will lie in some $\Gamma'$-parabolic neighborhood of $\overline{C}_Q$,
and the same argument applied to the sequence $z_i = \pi_Q(x_i) \rightarrow z_0 = \pi_Q(x_0)
\in F \cap {C}_Q$ also leads to a contradiction.

The case $\mathbf Q \subseteq \mathbf P$ may be handled by reversing the roles of $\mathbf P$
and $\mathbf Q$ in these arguments.  This completes the proof of Lemma \ref{lem-neutral} and
also the proof of Proposition \ref{prop-structure}.\qed

\section{Modified Hecke correspondence}\label{sec-modifiedHecke}\subsection{}
As in \S \ref{sec-notation}, $\mathbf G$ denotes a connected linear reductive algebraic group
defined over $\mathbb Q$, $D$ denotes the associated symmetric space,  $K' = A_GK(x_0)$ is the
stabilizer in $G$ of a chosen basepoint $x_0\in D$, $\Gamma \subset \mathbf G(\mathbb Q)$
is an arithmetic group and $X = \Gamma \backslash D.$ 
Throughout this section we fix a Hecke correspondence $(c_1,c_2):\overline{C}
\rightrightarrows \overline{X}$ defined by some element $g \in \mathbf G(\mathbb Q)$, with  $C
= \Gamma' \backslash D$ and $\Gamma' = \Gamma \cap g^{-1} \Gamma
g.$  Let $F = \{ w \in \overline{C}|\ c_1(w) = c_2(w) \}$ denote the fixed point set.
Fix a sufficiently \alarge regular $\Gamma$-equivariant parameter $b \in \mathcal B$ which is
so \alarge that the resulting tilings $\{D^P\}$ of $\overline{D}$, $\{X^P \}$ of
$\overline{X}$ and $\{ C^P \}$ of $\overline{C}$ are narrow (\S \ref{subsec-narrow}) with
respect to the Hecke correspondence.  Fix $\mathbf t\in A_{P_0}(\agt 1)$ dominant and regular,
with resulting shrink homeomorphism $Sh(\mathbf t):\overline{X} \to \overline{X}$ as in \S
\ref{sec-shrink}. Define the (shrink-) modified correspondence
\begin{equation}\label{eqn-modified}
(c'_1,c'_2):\overline{C} \rightrightarrows \overline{X}\end{equation}
 by $c'_1 = c_1$ and $c'_2 = Sh(\mathbf t)\circ c_2.$  Let $F' = \{ w \in \overline{C} |\
c'_1(w) = c'_2(w) \}$ denote the fixed point set of the modified correspondence.

\begin{prop}\label{prop-nofixedpt}
If $\mathbf t \in A_{P_0}(\agt 1)$ is chosen regular and sufficiently close to $\mathbf 1$,
then 
\begin{equation}\label{eqn-shrinkfix}
F'\cap C_Q = F \cap C_Q^0
\end{equation}
for each stratum $C_Q \subset \overline{C}$, and  
\begin{equation}\label{eqn-shrinkfix2}
c_1(F') \cap X_Q = c_1(F) \cap X_Q^0
\end{equation}
for each stratum $X_Q \subset \overline{X}.$
\end{prop}

\subsection{Proof}
The correspondence $\overline{C}$ has finitely many boundary
strata $C_P$ with the property that $c_1(C_P) = c_2(C_P).$  For each such stratum $C_P$,
according to proposition \ref{prop-restriction}, the Hecke correspondence is locally
isomorphic near $C_P$ to a parabolic Hecke correspondence $\Gamma'_P\backslash D[P]
\rightrightarrows \Gamma_P \backslash D[P]$ which is given by some $y \in P(\mathbb Q)$ and to
which we may uniquely associate a  torus factor $a_y=a_P \in A_P$ as in \S
\ref{subsec-restriction}.  Conjugating all these torus factors back to $\mathbf {S_{P_0}}$
gives a collection $\{a_1,a_2,\ldots,a_N \} \subset A_{P_0}$ of finitely many {\it standard
torus factors} (some of which may coincide and some of which may equal $\mathbf 1$) associated
to the Hecke correspondence $g$.  If $\mathbf{t} \in A_{P_0}(\agt 1)$ is chosen to be regular
and sufficiently close to $\mathbf 1$ then we can guarantee that the following condition
holds:  For all $\alpha \in \Delta$ and for all $i=1,2,\ldots,N$, if $\alpha(a_i) \alt 1$ then
$\alpha(a_i\mathbf t) \alt 1$ while if $\alpha(a_i) \age 1$ then $\alpha(a_i\mathbf t) \agt
1.$  Therefore, for any $\rho$ with $0 < \rho \le 1$, for all $\alpha \in \Delta$ and for all
$i= 1,2,\ldots, N$, the following holds:
\begin{equation}\label{eqn-smallt}\begin{cases}
 \alpha(a_i)\alt 1 &\implies \alpha(a_i)\alpha(\mathbf t)^{\rho} \alt 1\\
\alpha(a_i) \age 1 & \implies \alpha(a_i)\alpha(\mathbf t)^{\rho} \agt 1.
\end{cases}\end{equation}
Having made these choices, let us now prove Proposition \ref{prop-nofixedpt}.
Certainly $F \cap C_Q^0 = F' \cap C_Q^0$ because the shrink acts as the identity on $C_Q^0.$
So we only need to show that $F' \cap C_Q \subset
C_Q^0$, that is, we must show that the fixed points of the modified Hecke correspondence which
appear in the stratum $C_Q$ are all contained in the central tile of that stratum.  Suppose
otherwise and let $w\in C_Q$ be a fixed point of the modified correspondence which lies in
some tile $C^P$ for $\mathbf P \subset \mathbf Q,$ and $\mathbf P \ne \mathbf Q.$  Since the
shrink preserves tiles, it follows that $c_1(C^P) \cap c_2(C^P) \ne \phi$.  The tiling is
narrow so this implies that $c_1(C_P) = c_2(C_P).$  Set $\Delta_P = i(\Delta_Q)
\amalg J$ as in (\ref{eqn-decomposition}); then $J\ne \phi.$  By Proposition
\ref{prop-restriction}, locally near $C_P$ we may replace the Hecke correspondence by a
parabolic correspondence:  in other words, we may assume that
$g\in \mathbf P.$  Let $a_P \in A_P$ be the torus factor for the correspondence near $C_P$,
that is, if  $g= u_ga_gm_g\in \mathcal U_P A_PM_P$ is the Langlands decomposition then $a_P =
a_g.$  The point $c_1(w)= c'_1(w)=c'_2(w)$ lies in $X_Q \cap X^P \subset \overline{X}$.  Since
the shrink preserves tiles, $c_2(w) \in X_Q \cap X^P$ also. 

For any $\alpha \in \Delta_P$,
\begin{xalignat*}{2}
 f_{\alpha}^P(c'_2(w)) &= f_{\alpha}^P(c_2(w))\alpha_0(\mathbf t)^{\rho(r_{\alpha}^P(c_2(w)))}
&\text{by (\ref{eqn-fshrink})} \\
&= \alpha(a_P)\alpha_0(\mathbf t)^{\rho(r_{\alpha}^P(c_2(w)))}f_{\alpha}^P(c_1(w)) &\text{by
} (\ref{eqn-monitor})
\end{xalignat*}
where $\alpha_0 \in \Delta$ is the unique simple root which, after conjugation and restriction
to $\mathbf{S_P}$, agrees with $\alpha.$

This gives a contradiction:  First note that $\rho(r_{\alpha}^P(c_2(w))) \ne 0,$ for otherwise
we would have $r_{\alpha}^P(c_2(w))=1$ or $c_2(w) \notin X^P.$ As the shrink preserves tiles,
this would imply that $ c'_2(w) \notin X^P$ which is absurd.  So by (\ref{eqn-smallt}) the
factor $\alpha(a_P)\alpha_0(\mathbf t)^{\rho(r_{\alpha}^P(c_2(w)))} \ne 1.$  If we choose
$\alpha \in J$ then the assumption $c_1(w) \in X_Q$ implies that $f_{\alpha}^P(c_1(w)) \ne 0.$  
Therefore the point $w$ cannot be fixed by the modified correspondence, which proves
(\ref{eqn-shrinkfix}).

There are finitely many strata $C_P$ such that $c_1(C_P) = X_Q.$  To prove
(\ref{eqn-shrinkfix2}) it suffices to show, for each of these strata, that
\[ c_1(F' \cap C_P) \cap X_Q = c_1(F \cap C_P) \cap X_Q^0.\]
Write $F \cap C_P = (F \cap C_P^0) \cup \widetilde{F}$ as a disjoint union.  Then $\widetilde
{F}$ is contained in a union of tiles $C^{P'}$ with $\mathbf{P'} \subset \mathbf P$ a proper
inclusion.  Since $c_1$ takes tiles to tiles, it follows that $c_1(\widetilde{F}) \cap X_Q^0 =
\phi$ hence
\begin{align*}
c_1(F\cap C_P) \cap X_Q^0 = c_1(F\cap C_P^0) \cap X_Q^0 &= c_1(F \cap C_P^0) \\
&= c_1(F' \cap C_P) = c_1(F' \cap C_P) \cap X_Q \end{align*}
by (\ref{eqn-shrinkfix}).  \qed

\subsection{Tangential distance}\label{subsec-tangential}
Choose a regular $\Gamma$-invariant parameter  so that the associated tiling is narrow (\S
\ref{subsec-narrow}) with respect to the Hecke correspondence.   Choose $\mathbf t \in
A_{P_0}(\agt 1)$ to be regular and sufficiently close to $\mathbf 1$ as in Proposition
\ref{prop-nofixedpt}.  Suppose $C_Q$ is a stratum of $\overline{C}$ for which $F \cap C_Q \ne
\phi.$  Fix $e\in L_Q$ and let $F_Q(e)\subset C_Q$ denote the set of fixed points in $C_Q$ for
which $e$ is a characteristic element as in \S \ref{subsec-char-elt}.  Set
\begin{equation*}
F_Q^0(e) = F_Q(e) \cap C_Q^0 \text{ and } E = c_1(F_Q^0(e)) = c_1(F_Q(e)) \cap X_Q^0.
\end{equation*}
 Let $R_{X_Q}:X_Q \to X_Q^0$ and $R_{C_Q}:C_Q \to C_Q^0$ be
the retraction(s) and let $W_Q:\overline{X_Q} \to [0,1]$ be the exhaustion function of \S
\ref{subsec-retraction}.  A choice of $G$-invariant Riemannian metric on
$D$ induces Riemannian metrics on $C$, $X$, $C_Q$, and $X_Q.$ Define the tangential distance
$d_E:X_Q \to
[0,\infty)$ by
\begin{equation}\label{eqn-tangential}
d_E(x) =  W_Q(x) + \text{dist}_{X_Q}(R_{X_Q}(x),E)
\end{equation} 
where $\text{dist}_{X_Q}$ denotes the distance in $X_Q$ with respect to the Riemannian metric.
Then $d_E^{-1}(0) = E.$   Although the restriction of the Hecke correspondence to
$\overline{C}_Q$ is locally an isometry, composing
with $Sh(\mathbf t)$ has the following effect:  points near the boundary of $X_Q$ are moved
even closer to the boundary of $X_Q$ and hence they are moved away from $E.$
 This is the intuition behind the following lemma.

\begin{lem}\label{lem-dist-neutral}
There exists a neighborhood $V\subset C_Q$ of $F_Q^0(e)=F_Q(e)\cap C_Q^0$ such that 
\begin{equation}\label{eqn-dist-neutral}
d_E(c'_2(w)) \ge d_E(c'_1(w)) 
\end{equation}
for all $w\in V.$
\end{lem}

\subsection{Proof} The stratum closure $\overline{C}_Q$ is tiled by the collection of
intersections $C_Q^P = \overline{C}_Q \cap C^P$ with $\mathbf P \subseteq \mathbf Q.$  Let
$\overline{C}_Q^P$ denote the closure of such a tile.
Let $U_1 \subset \overline{C}_Q$ be a neighborhood of the closure $\overline{F}_Q(e)$ so that
for any rational parabolic subgroup $\mathbf P \subset \mathbf Q$, 
\begin{equation*}
\overline{F}_Q(e) \cap \overline{C}_Q^P = \phi \iff U_1 \cap \overline{C}_Q^P = \phi.
\end{equation*}
By Lemma \ref{lem-neutral} we may also assume that the Hecke correspondence preserves tile
boundaries in $U_1.$  The mapping $c_1$ preserves tiles, and the points in $F_Q(e)$ are fixed,
hence
\begin{equation*}
E = c_1(F_Q^0(e)) = c_2(F_Q^0(e)).\end{equation*}
By Proposition \ref{prop-char} (6), and for $i=1,2$, the mapping $c_i$ is one-to-one on
$F_Q(e).$  Moreover, it is locally an isometry.  It follows that, 
\begin{equation}\label{eqn-dist-neutral2}
\text{dist}_{X_Q}(c_1(w),E) = \text{dist}_{C_Q}(w, F_Q^0(e)) =
\text{dist}_{X_Q}(c_2(w), E)
\end{equation}
for $w\in C_Q$ in some neighborhood $U_2 \subset C_Q$ of $F_Q^0(e).$  The desired neighborhood
is
\begin{equation*} V = U_1 \cap U_2 \subset C_Q.
\end{equation*}
  
By  Proposition \ref{prop-structure} the restricted correspondence $\overline{C_Q}
\rightrightarrows \overline{X}_Q$ is neutral near $F_Q(e).$  So if $\mathbf P \subset \mathbf
Q$ and if $F_Q(e) \cap C^P \ne \phi$ and if $\Delta_P = i(\Delta_Q) \amalg J$ as in
(\ref{eqn-decomposition}) then by (\ref{eqn-exp-contr}),
\begin{equation}\label{eqn-merde}
r_{\alpha}^P(c_2(w)) = r_{\alpha}^P(c_1(w))
\end{equation}
for all $w\in C_Q^P$ and for all $\alpha\in J.$  Moreover,  by lemma 
\ref{lem-neutral} the correspondence preserves tiles near $F_Q(e)$, that is, for all $w\in
U_1$  and for all $\mathbf P \subseteq \mathbf Q$ we have
\begin{equation*}
w \in \overline{C}_Q^P \iff c_1(w) \in \overline{X}_Q^P \iff c_2(w) \in \overline{X}_Q^P.
\end{equation*}
By (\ref{eqn-commutes geodesic}) the correspondence commutes with the geodesic action of
$A_P.$ Therefore
\begin{equation*}
R_{X_Q}(c_1(w)) = c_1(R_{C_Q}(w)) \text{ and } R_{X_Q}(c_2(w)) = c_2(R_{C_Q}(w))
\end{equation*}

Now suppose $w\in V$ and $w\in C_Q^P$ for some $\mathbf P \subseteq \mathbf Q.$  If $\mathbf
P = \mathbf Q$ (that is, if $w\in C_Q^0$ lies in the central tile) then $c_i(w)\in C_Q^0$ as
well, in which case $R_{X_Q}(c_i(w)) =c_i(w)$, $W_Q(c_i(w))=0$, and $c_i(w) = c'_i(w).$  Then
(\ref{eqn-dist-neutral}) follows from (\ref{eqn-dist-neutral2}) and in fact, equality holds.

Now suppose $w\in V\cap C_Q^P$ for some $\mathbf P \ne \mathbf Q.$  Then $\mathbf P \subset
\mathbf Q,$ $c_1(C_P) = c_2(C_P)=X_P$ by Proposition \ref{prop-restriction}, and locally
near $X_P$ this correspondence is isomorphic to a parabolic Hecke correspondence, that is, we
may assume that $g \in \mathbf P(\mathbb Q).$  In the tile $C^P$ the retraction $R$
commutes with the geodesic action of $A_P$, cf. (\ref{eqn-retraction}), and so does the Hecke
correspondence, (\ref{eqn-commutes geodesic}), hence
\begin{equation*}
R_{X_Q}(c'_2(w)) = R_{X_Q}(Sh(\mathbf t) c_2(w)) = R_{X_Q}(c_2(w))=c_2(R_{C_Q}(w)).
\end{equation*}
So the second terms in (\ref{eqn-tangential}) are
equal:
\begin{align*}
\text{dist}_{X_Q}(R_{X_Q}c'_2(w),E) &= \text{dist}_{X_Q}(c_2(R_{C_Q}(w)),c_2(F_Q^0(e))) =
\text{dist}_{C_Q}(R_{C_Q}(w), F_Q^0((e)) \\
&= \text{dist}_{X_Q}(c_1(R_{C_Q}(w)),E) = \text{dist}_{X_Q}(R_{X_Q}c_1(w),E)
\end{align*}  
because both morphisms $c_1$ and $c_2$ are local isometries.  Now consider the first terms in
(\ref{eqn-tangential}).  Fix $w\in C_Q^P.$  For $\alpha\in \Delta_P$ set $\rho(\alpha) =
\rho(r_{\alpha}^P(c_2(w))).$  Using (\ref{eqn-exhaustion}), (\ref{eqn-rshrink}), and
(\ref{eqn-merde}) we find,
\begin{xalignat*}{2}
W_Q(c'_2(w)) &= 1- \underset{\alpha\in J}{\inf}\{ 
r_{\alpha}^{P}Sh(\mathbf t)c_2(w)\} &\   \\
&= 1 - \underset{\alpha\in J}{\inf}\{r_{\alpha}^P(c_2(w))
                \alpha_0(\mathbf t)^{-\rho(\alpha)} \}
\\ &= 1 - \underset{\alpha\in J}{\inf}\{r_{\alpha}^P(c_1(w)) \alpha_0
(\mathbf t)^{-\rho(\alpha)} \} \\
& \ge 1 - \underset{\alpha\in J}{\inf} \{ r_{\alpha}^P(c_1(w))\} \\
&= W_Q(c'_1(w))
\end{xalignat*}
which completes the proof of (\ref{eqn-dist-neutral}).  \qed

\subsection{Hyperbolic correspondences}\label{subsec-hypercorres}
Recall (\cite{LFP}) that the correspondence $(c_1,c_2):\overline{C} \rightrightarrows
\overline{X}$ is {\it weakly hyperbolic} near a connected component $F \subset
\overline{C}$ of the fixed point set, if there is a neighborhood $N(F')\subset \overline{X}$
of the image $F'=c_1(F)=c_2(F)$ and an {\it indicator mapping} $t=(t_1,t_2):N(F)\to \mathbb
R_{\ge 0} \times \mathbb R_{\ge 0}$ such that
\begin{enumerate}
\item the mapping $t$ is proper and subanalytic
\item the preimage of the origin $t^{-1}(0)=F'$ consists precisely of $F'$
\item there is a neighborhood $N(F)\subset \overline{C}$ so that $c_i(N(F)) \subset N(F')$
(for $i=1,2$), and so that for any $x\in N(F),$
\begin{align*}
t_1(c_1(x)) &\le t_1(c_2(x)) \\
t_2(c_1(x)) &\ge t_2(c_2(x))
\end{align*}\end{enumerate}

\subsection{The modified correspondence is hyperbolic}

Choose a tiling parameter $\mathbf b \in \mathcal B$ so that the associated tiling is narrow
with respect to the Hecke correspondence $(c_1,c_2):\overline{C} \rightrightarrows
\overline{X}.$  Choose $\mathbf t \in A_{P_0}(\agt 1)$ to be dominant, regular, and
sufficiently close to $\mathbf 1$ as in Proposition \ref{prop-nofixedpt} and equation
(\ref{eqn-smallt}).  Let $(c'_1,c'_2): \overline{C} \rightrightarrows \overline{X}$ be the
modified correspondence.

Suppose $C_Q$ is a stratum of $\overline{C}$ for which $F \cap C_Q \ne \phi.$  Then $c_i(C_Q)
= X_Q$ and we may write $\Delta_Q = \Delta^{+}_Q \cup \Delta^{-}_Q \cup \Delta^0_Q$ 
according to whether the simple root is expanding, contracting, or neutral near $X_Q$ as in \S
\ref{subsec-restriction}.  Fix $e\in L_Q$ and let $F_Q(e) \subset C_Q$ denote the set of fixed
points in $C_Q$ for which $e$ is a characteristic element as in \S \ref{subsec-char-elt}.  Let
$E = c_1(F_Q(e)) \cap X_Q^0.$  Let $V\subset C_Q$ be the neighborhood of $F_Q^0(e)=F_Q(e)\cap
C_Q^0$ described in Lemma \ref{lem-dist-neutral}.  
  Define $t=(t_1,t_2): T(\overline{X}_Q) \to \mathbb R_{\ge 0} \times \mathbb R_{\ge 0}$ by
\begin{align*}
t_1(x) &= \sum_{\alpha\in\Delta^{+}_Q} r_{\alpha}^Q(x) + d_E(\pi_Q(x))
 \\
t_2(x) &= \sum_{\alpha \in \Delta^{-}_Q} r_{\alpha}^Q(x) +
\sum_{\alpha\in\Delta^0_Q}r_{\alpha}^Q(x)
\end{align*}

\begin{thm}\label{thm-hyperbolic}
The mapping $(t_1,t_2)$ is an indicator mapping, with respect to which the modified Hecke
correspondence is hyperbolic in $\pi_Q^{-1}(V)$.
\end{thm}

\subsection{Proof}  The idea is that the composition with the $Sh(\mathbf t)$ converts
neutral directions (normal to a given stratum) into contracting directions but it does not
change the nature of the expanding or contracting (normal) directions.  It converts
distances within the stratum (which are preserved by the Hecke correspondence and hence
neutral) into expanding directions.

Since $c_1(C_Q) = c_2(C_Q) = X_Q,$ by Proposition \ref{prop-restriction} we may, locally near
$C_Q$, replace the Hecke correspondence with a parabolic Hecke correspondence determined by
some $g = u_ga_gm_g \in \mathbf Q(\mathbb Q).$  Let $a_Q = a_g \in A_Q$ be the torus factor
for the correspondence near $C_Q.$  It is easy to check that $t_1^{-1}(0) \cap
t_2^{-1}(0) = E.$  For any $w \in {T}(\overline{C}_Q)$ we have
\begin{xalignat*}{2}
r_{\alpha}^Q(c'_2(w)) &= r_{\alpha}^Q(Sh(\mathbf t)c_2(w)) 
= \alpha_0(\mathbf t)^{-\rho(r_{\alpha}^Qc_2(w))} r_{\alpha}^Qc_2(w) &\text{ by
(\ref{eqn-rshrink})}
\\
&=  \alpha_0(\mathbf t)^{-\rho(r_{\alpha}^Qc_2(w))} \alpha(a_Q)^{-1}r_{\alpha}^Q(c_1(w))
&\text{ by (\ref{eqn-exp-contr})} 
\end{xalignat*}
If $\alpha\in \Delta_Q^{-} \cup \Delta_Q^0$ then $\alpha_0(\mathbf
t)^{-\rho(r_{\alpha}^Qc_2(w))} \alpha(a_Q)^{-1} \le 1$ since both factors are $\le 1.$  This
proves that \begin{equation*}t_2c'_2(w) \le t_2c'_1(w)
\end{equation*}
If $\alpha\in \Delta_Q^{+}$ then $\alpha(a_Q) \alt 1$ so $\alpha_0(\mathbf
t)^{-\rho(r_{\alpha}^Qc_2(w))} \alpha(a_Q)^{-1} > 1$ by (\ref{eqn-smallt}).  Moreover, for all
$w\in \pi_Q^{-1}(V)$ we have
\begin{xalignat*}{2}
d_E\pi_Qc'_2(w) &= d_E\pi_Q Sh(\mathbf t)c_2(w) & \\
                &= d_E Sh(\mathbf t)\pi_Qc_2(w) &\text{by \S \ref{prop-shrink}(4)} \\
                &= d_E Sh(\mathbf t)c_2\pi_Q(w) &\text{by (\ref{eqn-intertwine2})} \\
                &\ge d_E c'_1\pi_Q(w) &\text{by  (\ref{eqn-dist-neutral})} \\
                &= d_E \pi_Qc'_1(w)
\end{xalignat*}
which completes the proof that
\begin{equation*}
  t_1c'_2(w) \ge t_1c'_1(w)  
\end{equation*}
and so also completes the proof that the modified Hecke correspondence is hyperbolic.  \qed

\newcommand{\Sd}{\mathbf{S^{\bullet}}}

\section{Local weighted cohomology with supports}\label{sec-localweighted}

\subsection{Quadrants}\label{subsec-quadrants}
 (See \cite{GKM} \S 7.14 p. 534.)\label{subsec-local2}
As in previous sections we suppose $\mathbf G$ is a connected reductive linear algebraic group
defined over $\mathbb Q$ and we denote the greatest $\mathbb Q$-split torus in its center 
by $\mathbf{S_G}.$  Let $\mathbf P$ be a rational parabolic subgroup with $\mathbf{S_P}$
the greatest $\mathbb Q$-split torus in the center of its Levi quotient $\mathbf{L_P}.$
Let $\Delta_P$ denote the simple positive roots of $\mathbf{S_P}$ occurring in 
$\mathfrak N_P = \text{Lie}(\mathcal U_P)$.  The elements $\alpha\in\Delta_P$ are
trivial on $\mathbf {S_G}$ and form a basis of 
$\chi^*_{\mathbb Q}(\mathbf{S'_P}) \subset \chi^*_{\mathbb Q}(\mathbf{S_P})$ where
$\mathbf{S'_P} = \mathbf{S_P}/\mathbf{S_G}.$  For any subset $J\subset \Delta_P$ as in \S
\ref{subsec-twoparab} let $\mathbf Q = \mathbf{P(J)}$ be the parabolic subgroup containing
$\mathbf P$ for which the corresponding torus
\[ \mathbf{S_J} = \mathbf{S_{P(J)}} \]
is the identity component of the intersection $\bigcap_{\alpha \in J}\ker(\alpha).$     Let
$\{t_{\alpha}\}$ be
the basis of the cocharacter group $\chi_*(\mathbf{S'_P})\otimes\mathbb Q$ which is dual to
the basis $\Delta_P$ so that $\langle\alpha,t_{\beta}\rangle = \delta_{\alpha,\beta}$
(with respect to the canonical pairing $\langle\cdot,\cdot\rangle$).  If $\overline{J} =
\Delta_P - J$ denotes the complement, then the cocharacter group $\chi_{*}^{\mathbb
Q}(\mathbf{S_J})$ is spanned by $\chi_*^{\mathbb Q}(\mathbf{S_G})$ and the set $\{ t_{\alpha}
\left|\right.\  \alpha \in \overline{J}\}.$

 Fix $\nu_P \in \chi^*_{\mathbb Q} (\mathbf {S_P})$ and $J\subset \Delta_P.$  Let $\gamma
\in \chi_{\mathbb Q}^{*}(\mathbf{S_P})$ and suppose that $\gamma|\mathbf{S_G} =
\nu_P|\mathbf{S_G}.$  Then $\gamma - \nu_P$ may be regarded as an element of $\chi_{\mathbb
Q}^*(\mathbf{S'_P})$ so we may define
\begin{align}\label{eqn-Inu}
I_{\nu_P}(\gamma) &= \left\{ \alpha \in \Delta_P |\ \langle\gamma -\nu_P, t_{\alpha}\rangle
<0\right\} \\
\label{eqn-Inu2}
\chi^*_{\mathbb Q}(\mathbf{S_P})\qdr{J} &= \left\{ \gamma \in
\chi^*_{\mathbb Q}
(\mathbf{S_P}) |\ I_{\nu_P}(\gamma) = J \text{ and } \gamma |\mathbf{S_G} = \nu_P|
\mathbf{S_G}.\right\}
\end{align}
This last set is called the {\it quadrant of type $J$}.  The disjoint union of the
$2^{|\Delta_P|}$ quadrants,\begin{equation*}
\coprod_{J \subseteq \Delta_P}\chi^*_{\mathbb Q}(\mathbf{S_P})\qdr{J}
=\left\{ \gamma \in \chi^*_{\mathbb Q}(\mathbf{S_P}) \left|\right.\ \gamma | \mathbf{S_G} =
\nu_P | \mathbf{S_G} \right\} \end{equation*}
is the subset of all characters whose restriction to $\mathbf{S_G}$ agrees with that of
$\nu_P.$  Taking $J = \phi$ gives
\[ \chi^*(\mathbf{S_P})\qdr{\phi} = \left\{ \gamma \in \chi^*_{\mathbb Q}(\mathbf{S_P})\left|
\right.\ \gamma |\mathbf{S_G} = \nu_P|\mathbf{S_G}\text{ and }\langle \gamma, t_{\alpha}
\rangle \ge \langle \nu_P,t_{\alpha} \rangle \text{ for
all } \alpha \in \Delta_P \right\} \]
which was denoted $\chi^*_{\mathbb Q}(\mathbf{S_P})_{+}$ in \cite{GHM} and was denoted
$\chi^*_{\mathbb Q}(\mathbf{S_P})_{\ge \nu_P}$ in \cite{GKM}.  It is the translate by $\nu_P$
of the positive cone $\left\{ \sum_{\alpha \in \Delta_P}m_{\alpha}\alpha \right\}$ with
$m_{\alpha} \ge 0.$  More generally, for $J \subset\Delta_P$  define
\[ \chi^*_{\mathbb Q}(\mathbf{S_P})_{\ge \nu_P(J)} = 
\left\{ \gamma \in \chi^*_{\mathbb Q}(\mathbf{S_P})\left| \right. \
\gamma |\mathbf{S_G} = \nu_P|\mathbf{S_G}\text{ and }
\langle \gamma, t_{\alpha} \rangle \ge \langle \nu_P, t_{\alpha} \rangle \text{ for all }
\alpha \in \overline{J} \right\}. \]
(The subset $\chi^*_{\mathbb Q}(\mathbf{S_P})_{\ge \nu_P(J)}$ was denoted $\chi^*_{\mathbb
Q}(\mathbf{S_P})_{+(J)}$ in \cite{GHM}.  One should also be aware that, in this paper and in
\cite{GKM} each $J$ denotes a subset of $\Delta_P$ while in \cite{GHM} the symbol $J$ was used
to denote a collection of standard maximal parabolic subgroups containing $\mathbf P.$  With
appropriate labelings, one is the complement of the other.) Then
\begin{equation}\label{eqn-difference}
\chi^*_{\mathbb Q}(\mathbf{S_P})\qdr{J} = \chi^*_{\mathbb Q}(\mathbf{S_P})_{\ge \nu_P(J)} - 
\bigcup_{K \subsetneq J} \chi^*_{\mathbb Q}(\mathbf{S_P})_{\ge \nu_P(K)}. \end{equation}
Equation (\ref{eqn-difference}) remains valid if we replace the union on the right hand side
by the union over those $K \subset J$ such that $|K| = |J|-1.$
 
If $H$ is an $\mathbf{S_P}$ module such that $\mathbf{S_G}$ acts on $H$ through the character
$\nu_P|\mathbf{S_G}$ then one may define $H\qdr{J}$ (resp. $H_{\ge \nu_P}$, resp. $H_{\ge
\nu_P(J)}$) to be the sum of those weight spaces $H_{\gamma}$ for which $\gamma \in
\chi^*_{\mathbb Q}(\mathbf{S_P})\qdr{J}$ (resp. $\gamma \in \chi^*_{\mathbb
Q}(\mathbf{S_P})_{\ge \nu_P}$, resp. $\gamma \in \chi^*_{\mathbb Q}(\mathbf{S_P})_{\ge
\nu_P(J)}$).

\subsection{Weighted cohomology}\label{subsec-WC}  
As in \S \ref{sec-notation}, let $D$ denote the symmetric space associated to $G$,  $K' =
A_GK(x_0)$ denote the stabilizer in $G$ of a fixed basepoint $x_0\in D$, $\Gamma \subset
\mathbf G(\mathbb Q)$ be a neat arithmetic group and $X = \Gamma \backslash D.$ 
Let $G \to GL(E)$ be a finite dimensional irreducible representation of $\mathbf G$ on
some complex vectorspace $E.$  It gives rise to a  local system 
$\mathbf E = (G/K') \times_{\Gamma} E$ on $X=\Gamma\backslash G/K.$  Let $\mathbf{P_0}$ be
the standard minimal rational parabolic subgroup with $\mathbf{S_0} = \mathbf{S_{P_0}}.$
Fix $\nu\in \chi^*_{\mathbb Q}(\mathbf{S_0})$ so that $\nu|\mathbf{S_G}$ coincides with
the character by which $\mathbf{S_G}$ acts on $E$.  Then $\nu$ defines a {\it weight
profile} in the sense of \cite{GHM}:  if $\mathbf Q \supseteq \mathbf{P_0}$ is a standard
rational parabolic subgroup then set $\nu_Q = \nu|\mathbf{S_Q}$ and
\begin{equation*}
\chi^*_{\mathbb Q}(\mathbf{S_Q})_{+} = \chi^*_{\mathbb Q}(\mathbf{S_Q})_{\ge \nu_Q}
= \chi^*_{\mathbb Q}(\mathbf{S_Q})_{\left\lceil\nu_Q,{\phi}\right\rfloor}
\end{equation*}  
These definitions may be extended to arbitrary rational parabolic subgroups by
conjugation.  We obtain from \cite{GHM} a complex of fine sheaves,
$\mathbf{W^{\nu}C^{\bullet}}(\mathbf E)$ on the reductive Borel-Serre compactification
$\overline{X}$ of $X$, whose (hyper)-cohomology groups $W^{\nu}H^*(\overline{X},
\mathbf E)$ are the weighted cohomology groups.  Let $i:X \to \overline{X}$ denote the
inclusion.  Recall from \cite{GHM} \S 13 that a choice of basepoint induces an isomorphism
\[ \mathbf H^j_{x}(Ri_*(\mathbf E)) \cong H^j(\mathfrak N_Q,E)\] 
between the stalk cohomology at a point $x \in X_Q$ of the complex of sheaves $Ri_*(\mathbf
E)$ and the Lie algebra cohomology of $\mathfrak N_Q.$  The weighted cohomology complex is
obtained by applying a weight truncation to the complex $Ri_*(\mathbf E)$ with the result that
its stalk cohomology becomes
\begin{equation}\label{eqn-stalkcohomology}
\mathbf H^j_{x}(\mathbf{W^{\nu}C}^{\bullet}(\mathbf E)) = H^j(\mathfrak N_Q,E)_{\ge \nu_Q}.
\end{equation}

\subsection{Remarks on sheaf theory}\label{subsec-sheafremarks}
In the next few sections we will need to use the formalism of the derived category of sheaves,
and some relations between the standard functors, for which we refer to \cite{IHII},
\cite{LFP}, \cite{BorelIH}, \cite{Iverson}, \cite{K&S}, \cite{GeM}.  Specifically, if $X$ is a
subanalytic set we denote by $D^b(X)$ the bounded (cohomologically-) constructible derived
category of sheaves of complex vectorspaces on $X$.  An element $\mathbf S^{\bullet} \in
D^b(X)$ is a complex of sheaves, bounded from below, whose cohomology sheaves $\mathbf
H^i(\mathbf S^{\bullet})$ are finite dimensional and are locally constant on each stratum of
some subanalytic stratification of $X$.  The hypercohomology of $\mathbf S^{\bullet}$ will be
denoted $H^*(\mathbf S^{\bullet})$ and the stalk cohomology at a point $x\in X$ will be
denoted $H^*_x(\mathbf S^{\bullet}).$  Denote by $\mathbf S[n]^{\bullet}$ the shifted sheaf,
$\mathbf S[n]^p = \mathbf S^{n+p}.$  The derived category $D^b(X)$ supports the standard
operations of $R\mathbf{Hom}$, $\otimes$, $Rf_*$, $Rf_!$, $f^*$, and $f^!$.  There are many
relations between these functors, of which we mention a few that we will use:

If $f:Y \to X$ is a normally nonsingular embedding (\cite{IHII} \S 5.4) then there
is a canonical isomorphism
\begin{equation}\label{eqn-normal}
f^!(\mathbf S^{\bullet}) \cong f^*(\mathbf S^{\bullet}) \otimes \mathcal O_{X/Y}[-d]
\end{equation}
where $\mathcal O_{X/Y}$ denotes the orientation bundle (or top exterior power) of the normal
bundle of $Y$ in $X$, and where $d$ denotes the codimension of $Y$ in $X$.  If $f:X \to
\text{pt}$ is the map to a point then $\mathbb D_X =f^!(\mathbb C)$ is the
dualizing complex.  If $X$ is an n-dimensional manifold (or even a rational homology manifold)
then $\mathbb D_X \cong \mathcal O_X[n]$ where $\mathcal O_X$ denotes the orientation bundle.

\subsection{Cohomology with supports}\label{subsec-supports}
Let $\overline{X}$ be a compact subanalytic set and let $\mathbf{S}^{\bullet}$ be
a (cohomologically) constructible complex of sheaves on $\overline{X}.$  Suppose
$Y \subset W \subset\overline{X}$ are locally closed subsets with inclusions
\begin{equation*}
\begin{CD}Y @>>{h_Y}> W @>>{j_W}> \overline{X}
\end{CD}
\end{equation*}
Define the {\it restriction of $\mathbf{S}^{\bullet}$ to $Y$ with compact supports in
$W$} to be the complex of sheaves  
\begin{equation}\label{eqn-Bd}
\mathbf{B^{\bullet}} = h_Y^{!}j_W^*\mathbf{S^{\bullet}}.\end{equation}
If $Y = \{y\}$ is 
a single point, then the cohomology of this complex is the relative cohomology group
\begin{equation}\label{eqn-relativecohomology}
H^m(\mathbf{B^{\bullet}}) = H^m(B_{\epsilon}\cap \overline{X}, \partial B_{\epsilon}\cap W;
\Sd),\end{equation}
where $B_{\epsilon}$ is a sufficiently small ball around $y$ (with respect to some
subanalytic embedding in Euclidean space) and $\partial B_{\epsilon}$ is its boundary.  

Now suppose $\overline{X}$ is the reductive Borel-Serre compactification of a locally
symmetric space $X = \Gamma \backslash G/K$ as in \S \ref{subsec-WC}, and that $\Sd =
\mathbf{W^{\nu}C^{\bullet}}( \mathbf{E})$ is the weighted cohomology sheaf constructed
with respect to some weight profile $\nu$ and local system $\mathbf E$ as in \S
\ref{subsec-WC}.
Let $ Y = {X}_P$ be some stratum and let $W = \overline{X}_Q$
be the closure of a larger stratum, corresponding to some rational parabolic subgroup
$\mathbf Q \supset \mathbf P.$  Form $\mathbf B^{\bullet} = 
h_Y^!j_W^*\mathbf{W^{\nu}C^{\bullet}(E)}$ as above.  Write $\Delta_P = i(\Delta_Q) \amalg I$
as in (\ref{eqn-decomposition}).

\begin{thm}\label{thm-supports}
The cohomology sheaf $\mathbf H^m(\mathbf B^{\bullet})$ is isomorphic to the local system on
$X_P$ which is associated to the following $L_P$-submodule of the $\mathfrak N_P$-cohomology,
\begin{equation}\label{eqn-Nquadrant}
 H^{m - |I|}(\mathfrak N_P,\mathbf E)\qdr{I}
\end{equation}\end{thm}

\subsection{Proof of theorem \ref{thm-supports}}\label{subsec-proof}  The proof follows
closely the computation \cite{GHM} \S 18 of the weighted cohomology of the link $\mathcal
L_y.$  First let us recall some generalities.  Each stratum $X_P$ of the reductive Borel-Serre
compactification $\overline{X}$ is a rational homology manifold.  If $\Gamma$ is neat, then
each stratum is a smooth manifold.  Suppose $\mathbf S^{\bullet}$ is a complex of
sheaves whose cohomology sheaves are locally constant on each stratum of $\overline{X}.$  Let
$Y=X_P \subset W = \overline{X}_Q$ as above.  The choice of basepoint $x_0\in D$ determines a
basepoint $y\in Y.$  Let $N_y \subset \overline{X}$
be a {\it normal slice} (cf. \cite{SMT} \S 5.4) to the stratum $Y$ at the
point $y.$  Let $k:N_y \cap W \to \overline{X}$ denote the inclusion, and let $i_y$ and $a_y$
denote the inclusions of $y$ into $Y$ and $N_y\cap W$ respectively.

Then the stalk cohomology of $\mathbf B^{\bullet} = h_Y^!j_W^*\mathbf S^{\bullet}$ is given by
\begin{equation}\label{eqn-supports2}
H^m_y(\mathbf{B^{\bullet}}) =
H^m(i_y^*h_Y^!j_W^*\Sd) \cong H^m(a_y^!k^*\Sd)
\end{equation}
which in turn may be identified with the relative cohomology group
\begin{equation}\label{eqn-pairsupports}
H^m(B_{\epsilon}\cap N_{y}, \partial B_{\epsilon} \cap N_y \cap
\overline{X}_Q; \mathbf S^{\bullet})
\end{equation}
(where $B_{\epsilon}$ is a sufficiently small ball around $y$, chosen with respect to
some locally defined subanalytic embedding of $\overline{X}$ into some Euclidean space).

These isomorphisms are deduced from the following fiber squares
\begin{equation}\label{eqn-supports}
\begin{CD}
y @>>{a_y}> N_y\cap W @>>{k_W}> N_y\\
@V{i_y}VV @VVV @VV{k_N}V \\
Y @>>{h_Y}> W @>>{j_W}> \overline{X}
\end{CD}\end{equation}
where $k = k_Nk_W.$
In the case that $\mathbf S^{\bullet} = \mathbf{W^{\nu}C^{\bullet}}$ we will compute
(\ref{eqn-pairsupports}) using the long exact cohomology sequence for the pair.
\subsubsection*{Step 1} Construct an isomorphism of $L_P$-modules,
\begin{equation}\label{eqn-sum2terms}
H^c(\partial B_{\epsilon}\cap N_y\cap \overline{X}_Q; \mathbf{W^{\nu}C^{\bullet}(E)}) \cong
H^{c-|I|+1}(\mathfrak N_P,\mathbf E)\qdr{I} \oplus
H^c(\mathfrak N_P,\mathbf E)_{\ge \nu_P} 
\end{equation}

In order to simplify notation, let us choose a labeling $\{\alpha_1,\alpha_2,\ldots,
\alpha_s\} = \Delta_P$ of the simple roots.  As in \cite{GHM} \S 8.8,
the link $\mathcal L_y = \partial B_{\epsilon} \cap N_y $ comes with a
natural mapping $\delta : \mathcal L_y \to \triangleright^{s-1}$ to the $s-1$ dimensional
simplex,
\begin{equation*}
\triangleright^{s-1} = \{ (x_1,x_2,\ldots,x_s)\in \mathbb R^s |\ 0 \le x_i \le 1
\text{ and } \sum x_i = 1 \}
\end{equation*}
For any subset $J\subset \{1,2,\ldots,s\}$ let $\overline{J}$ denote its complement.
Associated  to $J$ there is a (closed) face of dimension $|J|-1$,
\[ \triangleright_J = \{ x\in\triangleright^{s-1} |\ x_j = 0 \text{ for all } j\in 
\overline{J} \} \]
whose interior we denote by $\triangleright_J^o.$
Each $\triangleright _{\{j\}}$ is a vertex of $\triangleright^{s-1}$; the face
$\triangleright_J$ is spanned by the vertices $\triangleright_{\{j\}}$ such that 
$j\in J.$ Let $U_{\{j\}}= St(\triangleright_{\{j\}})$ be the open star of the vertex
$\triangleright_{\{j\}}.$  These form a covering of $\triangleright^{s-1}$  whose
multi-intersections we denote by 
\[ U_{J} = \bigcap_{j\in J} U_{\{j\}}. \]
Then 
\[U_J = St(\triangleright_J^o)=\bigcup \{ F^o |\ F \text{ is a face of $\triangleright^{s-1}$
and } F \supseteq \triangleright_J\} \]
is the open star of the interior of the face $\triangleright_J.$ 

If we think of stratifying the simplex $\triangleright^{s-1}$ by the interiors of its faces,
then
the mapping $\delta: \mathcal L_y \to \triangleright^{s-1}$ is a stratified mapping: for any
$J\subset \{1,2,\ldots,s\}$ it maps $\mathcal L_y \cap X_{P(J)}$ to the interior 
$\triangleright_J^o$ of the face $\triangleright_J,$ and in particular
\begin{equation*}
\mathcal L_y \cap \overline{X}_Q= \delta^{-1} (\triangleright_{I}).
\end{equation*}  The fiber over any interior point $s \in
\triangleright_I^o$ is the nilmanifold $(\Gamma \cap \mathcal U_Q )\backslash 
\mathcal U_Q.$  As in \cite{GHM} \S 18.5, the (weighted) cohomology of $\mathcal L_y \cap
\overline{X}_Q$ can be computed using the Mayer-Vietoris spectral sequence for the covering by
open stars (for $i\in I$),
\begin{equation*}
V_{\{i\}}= \delta^{-1}(U_{\{i\}}\cap \triangleright_I)
\end{equation*} of the vertices of $\triangleright_I.$ Set $V_J = \delta^{-1}(U_J
\cap\triangleright_I).$   The groups $E_1^{a,b}$ are cohomology groups of multi-intersections
of open sets in this covering, and were computed in \cite{GHM} lemma 18.5,
\begin{align*}
E_1^{a,b} &= \bigoplus_{\begin{subarray}{}|J| = a+1 \\ J\subset I    \end{subarray}}
W^{\nu}H^b(\bigcap_{j\in J} V_{\{j\}};E) 
= \bigoplus_{\begin{subarray}{}|J| = a+1 \\ J\subset I    \end{subarray}}
W^{\nu}H^b(V_{J};E) \\
&=  \bigoplus_{\begin{subarray}{}|J| = a+1 \\ J\subset I    \end{subarray}}
H^b(\mathfrak N_P,E)_{\ge\nu_{P(J)}}.
\end{align*} 
  
The $E_1$ differential is given (up to sign) by inclusion, so the argument of \cite{GHM} \S
18.7 applies here as well: the spectral sequence collapses at $E_2$, which has only two
possibly nonzero columns: $E_2^{0,b} = H^b(\mathfrak N_P,E)_{\ge\nu_P}$ and, using
(\ref{eqn-difference}), 
\begin{equation*}
E_2^{|I| -1,b} = \frac{H^b(\mathfrak N_P,E)_{\ge\nu_{P(I)}}}
{\displaystyle\sum_{\begin{subarray}{}
|K| = |I| - 1 \\ K \subset I\end{subarray}}
H^b(\mathfrak N_P,E)_{\ge\nu_{P(K)}}} = H^b(\mathfrak N_P,E)\qdr{I}
\end{equation*}
which contributes to $W^{\nu}H^*(\delta^{-1}St(\triangleright_{ I}),E)$ in degree $|I|-1+b.$
So we obtain a split short exact sequence (with $c=|I|-1+b$),
\begin{equation*}
0 \to H^{c-|I| +1}(\mathfrak N_P,E)\qdr{I} \to
W^{\nu}H^c(\mathcal L_y \cap \overline{X}_Q;E) \to H^c(\mathfrak N_P,E)_{\ge\nu_P} \to 0
\end{equation*}
which completes the proof of (\ref{eqn-sum2terms}).

\subsubsection*{Step 2}
As in \cite{GHM} \S 18.11, the long exact sequence for the pair (\ref{eqn-pairsupports})
splits into split short exact sequences,
\begin{equation*}
0 \to H^c(B_{\epsilon}\cap N_y) \to H^c(\partial B_{\epsilon}\cap N_y\cap \overline{X}_Q) \to
H^{c+1}(B_{\epsilon}\cap N_y, \partial B_{\epsilon} \cap N_y\cap \overline{X}_Q )\to 0.
\end{equation*}
But $H^*(B_{\epsilon} \cap N_y) = H^*(B_{\epsilon}) = H^*(\mathfrak N_P,E)_{\ge \nu_P}$ is the
stalk cohomology at $y$ of the weighted cohomology sheaf.  This kills the second summand in
(\ref{eqn-sum2terms}), leaving
\begin{equation}\label{eqn-localWH}
\mathbf H^m_y(\mathbf B^{\bullet})\cong
H^m(B_{\epsilon}\cap N_{y}, \partial B_{\epsilon} \cap N_y \cap
\overline{X}_Q; \mathbf{W^{\nu}C^{\bullet}(E)}) \cong H^{m - |I|}(\mathfrak N_P,\mathbf
E)\qdr{I}
\end{equation}

\subsubsection*{Step 3}  We briefly indicate why the isomorphism (\ref{eqn-localWH}) extends
to an isomorphism of flat vectorbundles on $X_P$,
\begin{equation*}
\mathbf H^{m}(\mathbf B^{\bullet}) \cong H^{m - |I|}(\mathfrak N_P ,\mathbf E)\qdr{I}
\times_{\Gamma_{L(P)}}L_P/K_PA_P
\end{equation*}
(where $\Gamma_{L(P)} = \nu_P(\Gamma\cap P)$ is the projection of $\Gamma \cap P$ to the Levi
quotient $L_P$ and where it acts on $H^*(\mathfrak N_P,\mathbf E)$ by conjugation).  Let $i:X
\to \overline{X}$ denote the inclusion.  In \cite{GHM} \S 17, special differential forms are
used in order to identify the restriction $\mathbf H^m(Ri_*\mathbf E)|X_P$ with the flat
vectorbundle 
\begin{equation*}
H^m(\mathfrak N_P,E)\qdr{\phi} \times_{\Gamma_{L(P)}}L_P/K_PA_P.
\end{equation*}
But each of the cohomology groups appearing in Step 2 (above) is an $L_P$-submodule of
$H^*(\mathfrak N_P,E)$ and the corresponding bundle on $X_P$ is a sub-bundle of $\mathbf
H^*Ri_*(\mathbf E)|X_P$ (while the shift by $|I|$ corresponds to tensoring with a trivial
vectorbundle on $X_P$).  So it suffices to verify that the stalk cohomology modules agree at
the basepoint, which we have done.  \qed

\subsection{Kostant's theorem}\label{subsec-Kostant}
In this section we will use
Kostant's theorem \cite{Kostant} to explicitly evaluate the cohomology group
(\ref{eqn-Nquadrant}).  Let $\mathbf B \subset \mathbf G$ be a Borel subgroup (over $\mathbb
C$), chosen so that $\mathbf B(\mathbb C) \subset \mathbf{P_0}(\mathbb C) \subset \mathbf
P(\mathbb C).$  Choose a maximal torus $\mathbf T$ (over $\mathbb C$) of $\mathbf G$ so that
\begin{equation}\label{eqn-tori}
 \mathbf{S_P}(\mathbb C) \subset \mathbf{S_0}(\mathbb C) \subset \mathbf{T}(\mathbb C)
\subset \mathbf{B_L}(\mathbb C) 
\end{equation}  
where $\mathbf{B_L} = \mathbf B \cap \mathbf{L_P}$ is the corresponding Borel subgroup of
$\mathbf{L_P}.$  This gives rise to root systems $\Phi_G = \Phi(\mathbf{G}(\mathbb C),
\mathbf{T}(\mathbb C))$
and $\Phi_L = \Phi(\mathbf{L_P}(\mathbb C), \mathbf{T}(\mathbb C))$ with positive roots
$\Phi^{+}_G = \Phi(\mathcal U_{\mathbf B}(\mathbb C), \mathbf{T}(\mathbb C))$ and
$\Phi^{+}_L = \Phi_L \cap \Phi^{+}_G$ (determined by the Borel subgroups $\mathbf B
\subset \mathbf G$ and $\mathbf {B_L} \subset \mathbf {L_P}$ respectively.)  
 Let $\rho_B = \frac{1}{2}\sum_{\alpha \in \Phi_G^{+}}\alpha.$ 

Let $W_G = W(\mathbf G(\mathbb C), \mathbf T(\mathbb C))$ denote the Weyl group of
$\mathbf G(\mathbb C)$ and let $W_P = W(\mathbf{L_P}(\mathbb C), \mathbf T(\mathbb C))$ denote
the Weyl group of $\mathbf{L_P}(\mathbb C).$  The choice of $\mathbf B$ determines a
length function $\ell$ on $W_G.$  Let $W^1_P\subset W_G$ denote the set of Kostant
representatives:  the unique elements of minimal length from each of the cosets $W_Px \in
W_P \backslash W_G.$  As in \cite{Springer} \S 10.2 or \cite{Vogan} \S 3.2.1, it may also be
described as the set 
\begin{equation*}
W^1_P = \{ w\in W|\ w^{-1}(\Phi^+_L) \subset \Phi^+_G \}
\end{equation*}
(and depends on the choice of $\mathbf{B_L} \subset \mathbf{L_P}$).

If $\beta \in \chi^*(\mathbf T(\mathbb C))$ is $\mathbf{B_L}$-dominant, let us write
$V^L_{\beta}$ for the irreducible $\mathbf{L_P}$-representation with highest weight
$\beta.$  Let $\lambda_B \in \chi^*(\mathbf T(\mathbb C))$ be the highest weight of the
irreducible representation $E$ of $\mathbf{G}.$  Kostant's theorem states that for all $w
\in W^1_P$, the weight $w(\lambda_B + \rho_B) - \rho_B$ is $\mathbf{B_L}$-dominant, and
that as an $\mathbf{L_P}$-module, the cohomology group $H^i(\mathfrak N_P,E)$ is
isomorphic to 
\[\bigoplus \{ V^L_{w(\lambda_B + \rho_B)-\rho_B} |\ w \in W^1_P \text{ and } \ell(w) = i
\}.\]
If $w\in W_G$ then the character $w(\lambda_B +\rho_B) - \rho_B -\nu$ is trivial on
$\mathbf{S_G}$ so we may define
\begin{equation}\label{eqn-IWeyl}
 I_{\nu}(w) = \{ \alpha \in \Delta_P |\
\langle (w(\lambda_B + \rho_B)-\rho_B -\nu)|\mathbf{S'_P}, t_{\alpha} \rangle < 0\} 
\end{equation}
where $\{t_{\alpha}\}$ form the basis of the cocharacter group $\chi_*^{\mathbb
Q}(\mathbf{S'_P})$ which is dual to the basis $\Delta_P$ of simple roots, cf.
(\ref{eqn-Inu}).  So in the notation of (\ref{eqn-Inu}),
\[I_{\nu}(w)  = I_{\nu}(\gamma) \text{ where } 
\gamma= (w(\lambda_B + \rho_B) - \rho_B) |\ \mathbf{S_P}.\]
To summarize we have,
\begin{prop}\label{prop-Kostant}
 Kostant's theorem determines an isomorphism of graded $\mathbf{L_P}$-modules, 
\begin{equation}
H^*(\mathfrak N_P,\mathbf E)\qdr{I} \cong \bigoplus_{\substack{w\in W^1_P \\
I_{\nu}(w)=I}} V^L_{w(\lambda_B + \rho_B) - \rho_B}[-\ell(w)]
\end{equation}
where the sum is taken over all $w\in W^1_P$ such that $ I_{\nu}(w) = I $, and where
$V^L_{\beta}[-m]$ means that the irreducible $\mathbf{L_P}$-module $V^L_{\beta}$ appears
in degree $m$.
\end{prop}

%%%%%%%%%%%%%%%%%%%%%%%%%%%%%%%%%%%%%%%%%

%%%%%%%%%%%%%%%%%%%%%%%%%%%%%%%%%%%%%%%%%%%%%%%%%%%%%%%%
\section{Lefschetz numbers} \label{sec-Lefschetz}
\subsection{}\label{subsec-Lef1}
In this section we recall the Lefschetz fixed point theorem for hyperbolic correspondences
from \cite{LFP} \S 10.3.

Suppose $\overline{C}$, $\overline{X}$ and $\overline{Y}$ are compact subanalytic spaces and
that $c=(c_1,c_2):\overline{C} \to \overline{X}\times \overline{Y}$
is a subanalytic mapping.  (The bars are used so that the notation here will agree with
that in the rest of the paper.)  Let $\mathbf S^{\bullet}\in D^b(\overline{X})$ be a (bounded
from below) complex of (cohomologically) constructible sheaves on $\overline{X}$ and let
$\mathbf T^{\bullet}\in D^b(\overline{Y})$ be a (bounded from below) complex of
(cohomologically) constructible sheaves on $\overline{Y}.$  Since $c$ is proper we have $c_* =
c_{!}.$  A {\it lift} of the correspondence $\overline{C}$ to the sheaf level (\cite{Verdier,
GI,BorelIH}) is a morphism
\begin{equation}\label{eqn-sheaf}
\Phi: c_2^*\mathbf T^{\bullet} \to c_1^{!}\mathbf S^{\bullet} \end{equation}
Such a morphism induces a homomorphism $H^*(\overline{Y};\mathbf T^{\bullet}) \to
H^*(\overline{X};\mathbf S^{\bullet})$ as follows.  First apply $(c_1)_{!}$ and adjunction to
obtain a morphism
\begin{equation}\label{eqn-c!c*}
 (c_1)_{!}c_2^*\mathbf T^{\bullet} \to (c_1)_{!}c_1^{!} \mathbf S^{\bullet} \to
\mathbf S^{\bullet} \end{equation}
Let $p:\overline{X}\to \text{pt}$ and $q:\overline{Y}\to \text{pt}$ be the map to a point.
Then the diagram 
\begin{equation*}\begin{CD}
\overline{X}\times \overline{Y} @>>{\pi_2}> \overline{Y} \\ @V{\pi_1}VV @VV{q}V \\
\overline{X} @>>{p}> \text{pt}
\end{CD}\end{equation*}
is a fiber square so there is an adjunction natural transformation \cite{LFP} (2.6b),
$q_{!}(\pi_2)_* \to p_*(\pi_1)_{!}.$  Apply $q_!$ to the adjunction morphism $\mathbf
T^{\bullet} \to (c_2)_*c_2^* \mathbf T^{\bullet}$ and use (\ref{eqn-c!c*}) to obtain
\begin{align*}
q_!\mathbf T^{\bullet} &\to q_!(c_2)_*c_2^*\mathbf T^{\bullet} = q_!(\pi_2)_*c_*c_2^*\mathbf
T^{\bullet}\\
&\to p_*(\pi_1)_!c_*c_2^*\mathbf T^{\bullet} = p_*(c_1)_!c_2^*\mathbf T^{\bullet} 
\to p_*\mathbf S^{\bullet} \end{align*}
This morphism induces the desired mapping on cohomology.  (It may also be constructed by
applying $p_!(c_2)_*$ to (\ref{eqn-sheaf}) rather than $q_*(c_1)_!$.) 

In what follows, we suppose $\overline{X}=\overline{Y}$ and $\mathbf S^{\bullet} = \mathbf
T^{\bullet}$, so  $c=(c_1,c_2):\overline{C} \to \overline{X}\times \overline{X}$ is a
correspondence on $\overline{X}$ and $\Phi:c_2^*\mathbf S^{\bullet}
\to c_1^!\mathbf S^{\bullet}$ is a lift to the sheaf level.  The Lefschetz fixed point theorem
states that the resulting Lefschetz number
\begin{equation*}
L(\mathbf S^{\bullet}, \overline{C}) = \sum_{i\ge 0} \text{Tr
}(\Phi^*:H^i(\overline{X};\mathbf S^{\bullet}) \to
H^i(\overline{X};\mathbf S^{\bullet})) = \sum_F L(\mathbf S^{\bullet}, \overline{C},
F)\end{equation*}
is a sum of locally defined contributions $L(\mathbf S^{\bullet}, \overline{C}, F)$, one for
each connected component $F\subset \overline{C}$ of the fixed point set of the correspondence
$\overline{C}$.

Let $F\subset \overline{C}$ be a connected component of the fixed point set
and suppose that the correspondence $\overline{C}$ is weakly hyperbolic (\S
\ref{subsec-hypercorres}) near $F'=c_1(F)=c_2(F)$ with indicator mapping
$t:W\to \mathbb R_{\ge 0}\times \mathbb R_{\ge 0}$.  (This means that $W\subset \overline{X}$
is a neighborhood of $F'$, that $t$ is a proper subanalytic mapping such that
$t^{-1}(0,0)= F'$, and that for all $x\in c_1^{-1}(W)\cap c_2^{-1}(W)$ we have $t_1c_1(x)\le
t_1c_2(x)$ and $t_2c_1(x) \ge t_2c_2(x)$.)  Denote by $h$ and $j$ the inclusions
\begin{equation*}\begin{CD}
F' @>h>> t^{-1}(\mathbb R_{\ge 0}\times \{0\}) @>j>> \overline{X}
\end{CD}\end{equation*}
of $F'$ into the ``expanding set'' or ``unstable manifold'' $F^{-} = t^{-1}(\mathbb R_{\ge
0}\times \{0\}),$ and of $F^{-}$ into $\overline{X}$.

Let $\mathbf A^{\bullet} = h^{!}j^*(\mathbf S^{\bullet})$ as in \S \ref{subsec-supports}.
Then the lift $\Phi$ determines a lift $\Psi':c_2^*\mathbf A^{\bullet} \to c_1^{!}\mathbf
A^{\bullet}$ which, by adjunction, induces an endormorphism  $\Psi:\mathbf A^{\bullet} \to
\mathbf A^{\bullet}$ (which covers the identity mapping on $F'$). In \cite{LFP} we prove:
\begin{thm}\label{thm-topological}  The contribution $L(\mathbf S^{\bullet},\overline{C},F)$
of $F$ to the global Lefschetz number $L(\mathbf S^{\bullet},\overline{C})$  is given by
\begin{equation*}
 L(\mathbf S^{\bullet},\overline{C},F) = \sum_{i \ge 0} (-1)^i \text{Tr} (\Psi^*:
H^i(F';\mathbf{A^{\bullet}}) \to H^i(F';\mathbf{A^{\bullet}}))
\end{equation*}
\end{thm}

Moreover, if $F' = \coprod_{\alpha=1}^m F'_{\alpha}$ is stratified so that the pointwise
Lefschetz number $n(x) = \sum_{i \ge 0}(-1)^i \text{Tr} (\Psi^{*}_x:
H^i_x(\mathbf A^{\bullet})\to H^i_x(\mathbf A^{\bullet}))$ is constant on each stratum, then
the local contribution is the sum over strata,
\begin{equation}\label{eqn-Eulerchar}
 L(\mathbf S^{\bullet},\overline{C},F) = \sum_{\alpha=1}^m \chi_c(F'_{\alpha})n(x_{\alpha})
\end{equation}
(where $x_{\alpha} \in F'_{\alpha}$ and where $\chi_c$ denotes the Euler characteristic with
compact supports).

\subsection{Morphisms and weighted cohomology}\label{subsec-morph}
In this section we show how to lift any morphism to the weighted cohomology sheaf.
As in \S \ref{sec-notation}, $\mathbf G$ denotes a connected linear reductive algebraic group
defined over $\mathbb Q$, $D$ denotes the associated symmetric space,  $K' = A_GK(x_0)$ is the
stabilizer in $G$ of a fixed basepoint $x_0\in D$, $\Gamma \subset \mathbf G(\mathbb Q)$
is a neat arithmetic group and $X = \Gamma \backslash D.$ As in \S \ref{subsec-WC} let $\tau:
G\to GL(E)$ be a finite dimensional irreducible representation on some complex vectorspace.
It gives rise to the local coefficient system (flat homogeneous vectorbundle) $\mathbf E =
(G/K')\times_{\Gamma}E$ which is the quotient of $(G/K')\times E$ under the equivalence
relation $(xK',v) \sim (\gamma xK', \tau(\gamma)v)$ for all $\gamma \in \Gamma.$  Denote by
$[xK',v]\in \mathbf E$ the resulting equivalence class.   Let $\mathbf{P_0}$ be
the standard minimal rational parabolic subgroup with $\mathbf{S_0} = \mathbf{S_{P_0}}.$  Fix
$\nu\in \chi^*_{\mathbb Q}(\mathbf{S_0})$ so that $\nu|\mathbf{S_G}$ coincides with the
character by which $\mathbf{S_G}$ acts on $E$ and let $\mathbf{W^{\nu}C}^{\bullet}
(\overline{X};\mathbf E)$ denote the resulting weighted cohomology complex of sheaves on
$\overline{X}.$

Suppose $\Gamma' \subset \Gamma$ is a subgroup of finite index, set $C = \Gamma' \backslash
G/K$, and let $f:C \to X$ be a morphism, i.e., there exists $g\in \mathbf G(\mathbb Q)$ such
that $g\Gamma'g^{-1} \subset \Gamma$ and $f(\Gamma'xK) = \Gamma gxK$.  Let $\mathbf{E'} \to C$
be the local coefficient system on $C$ which is determined by the representation $\tau:G \to
GL(E)$.  The morphism $f$ is covered by a mapping $\mathbf{E'}\to \mathbf E$ of local systems
given by $[xK,v] \mapsto [gxK, \tau(g)v]$.  This mapping is easily seen to be well defined,
and it induces an isomorphism of local systems $\mathbf{E'} \cong f^{*}(\mathbf E)$ on $C$.
Since $f:C \to X$ is an unramified finite covering, it further induces a canonical
quasi-isomorphism of the sheaves of smooth differential forms with coefficients in this local
system,  $f^*\boldsymbol{\Omega}^{\bullet}(X; \mathbf E) \to\boldsymbol{\Omega}^{\bullet}
(C;\mathbf{E'}).$ 

The morphism $f:C\to X$ admits a unique continuous extension $\overline{f}:\overline{C} \to
\overline{X}$ to the reductive Borel-Serre compactifications (Lemma \ref{lem-extension}).
 If $i_C:C\to \overline{C}$ and $i_X:X \to \overline{X}$ denote the inclusions then the
adjunction mapping \cite{LFP} equation (2.5a),
\begin{equation*}\begin{CD}
\overline{f}^*(i_{X})_*\boldsymbol{\Omega}^{\bullet}(X;\mathbf E) @>{\cong}>>
(i_{C})_*f^*\boldsymbol{\Omega}^{\bullet}(X;\mathbf E) @>{\cong}>>
(i_{C})_*\boldsymbol{\Omega}^{\bullet}(C;\mathbf E)
\end{CD} \end{equation*}
is a quasi-isomorphism.  It is easy to see that this induces a quasi-isomorphism 
\begin{equation}
\overline{f}^*\mathbf{W^{\nu}C^{\bullet}}(\overline{X};\mathbf E)
\to \mathbf{W^{\nu}C^{\bullet}}(\overline{C};\mathbf{E'})
\end{equation}
of weighted cohomology sheaves.  (In fact the whole construction of the weighted cohomology
sheaf on $\overline{X}$ pulls back to the construction of weighted cohomology on
$\overline{C}$.)  

Let $g\in \mathbf G(\mathbb Q).$  Then $g$ gives rise to a Hecke correspondence
$(c_1,c_2):\overline{C} \to \overline{X}$.  Here, $\overline{C}$ is the reductive Borel-Serre
compactification of $C = \Gamma'\backslash G/K$ with $\Gamma' = \Gamma \cap g^{-1}\Gamma g.$
Both mappings $c_1$ and $c_2$ are finite so there are natural isomorphisms of functors $c_i^*
\cong c_i^{!}$ and $(c_i)_* \cong (c_i)_{!}$ (for $i = 1,2$).  From the preceding paragraph we
obtain a canonical lift
\begin{equation}\label{eqn-lift}
\Phi: c_2^*\mathbf{W^{\nu}C^{\bullet}}(\overline{X};\mathbf E) \to 
c_1^{!}\mathbf{W^{\nu}C^{\bullet}}(\overline{X};\mathbf E)
\end{equation}
to the weighted cohomology sheaves, which is given by the composition
\begin{equation*}\begin{CD}
c_2^*\mathbf{W^{\nu}C^{\bullet}}(\overline{X};\mathbf E) @>{\cong}>>
\mathbf{W^{\nu}C^{\bullet}}(\overline{C};\mathbf E) @<{\cong}<<
c_1^*\mathbf{W^{\nu}C^{\bullet}}(\overline{X};\mathbf E)\cong
c_1^{!}\mathbf{W^{\nu}C^{\bullet}}(\overline{X};\mathbf E).
\end{CD}\end{equation*}

\subsection{Computation of the local contribution}  \label{subsec-computation}
For the remainder of \S
\ref{sec-Lefschetz}, fix a Hecke correspondence $\overline{C} \rightrightarrows \overline{X}$
which is determined by some element $g \in \mathbf G(\mathbb Q).$  Fix a regular
$\Gamma$-equivariant parameter $b \in \mathcal B$ which is so \alarge that the resulting
tilings $\{D^P\}$ of $\overline{D}$, $\{X^P \}$ of $\overline{X}$ and $\{ C^P \}$ of
$\overline{C}$ are  narrow (\S \ref{subsec-narrow}) with respect to the Hecke correspondence.
Choose $\mathbf t \in A_{P_0}(\agt 1)$ to be regular, dominant, and sufficiently close to
$\mathbf 1$ as in Proposition \ref{prop-nofixedpt}, with resulting shrink homeomorphism
$Sh(\mathbf t),$ and let $(c'_1,c'_2):\overline{C} \rightrightarrows \overline{X}$ be the
resulting modified correspondence.  It is easy to see that $Sh(\mathbf
t)^*(\mathbf{W^{\nu}C^{\bullet}}) \cong \mathbf{W^{\nu}C^{\bullet}}$ so we may consider
(\ref{eqn-lift}) to be a lift of the modified correspondence as well.  

Suppose the Hecke correspondence preserves some stratum $C_P.$  According to Proposition
\ref{prop-restriction}, locally near $C_P$ the correspondence is isomorphic to the parabolic
Hecke correspondence $\Gamma'_P\backslash D[P] \rightrightarrows \Gamma_P \backslash D[P]$
which is given by some $y \in \mathbf{P}(\mathbb Q)\cap \Gamma g \Gamma$ and to which we may
associate a decomposition $\Delta_P = \Delta^+_P \cup \Delta^-_P \cup \Delta^0_P$ of the
simple roots.  Suppose that $C_P$ contains fixed points and denote by $F_P(e) \subset C_P$ the
set of fixed points with characteristic element $e\in \Gamma_L  \bar y \Gamma_L \subset
\mathbf{L_P}(\mathbb Q).$ 

By Proposition \ref{prop-char} the torus factor $a_e\in A_P$ of $e$ coincides with the torus
factor $a_y$ so the set $\Delta_P^{+}$ (resp. $\Delta_P^{-}$, resp. $\Delta_P^0$) consists of
those simple roots $\alpha \in \Delta_P$ for which $\alpha(a_e) \alt 0$ (resp. $\agt 0$, resp.
$=0$).  Hence we may write $\Delta_P^{+} = \Delta_P^{+}(e)$ (resp. $\Delta_P^{-} =
\Delta_P^{-}(e)$, resp. $\Delta_P^0 = \Delta_P^0(e)$).

As in \S \ref{subsec-Kostant}, choose a Borel pair $\mathbf{T}(\mathbb C) \subset
\mathbf{B}(\mathbb C)$ so that (\ref{eqn-tori}) holds.  Assume the local system $\mathbf E$
arises from an irreducible representation of $\mathbf G$ with highest weight $\lambda_B\in
\chi^*(\mathbf{T}(\mathbb C)).$  Let $\rho_B = \frac{1}{2}\sum_{\alpha \in
\Phi_G^{+}}\alpha\in \chi^*(\mathbf{T}(\mathbb C))$ denote the half-sum of
the positive roots.  Let $r = [\Gamma \cap\mathcal U_P : \Gamma' \cap \mathcal U_P].$

\begin{thm}\label{thm-answer}  The contribution to the Lefschetz number from the fixed point
constituent $F_P(e)$ is:
\begin{equation*}
r\chi_c(F_P(e)) (-1)^{|\Delta_P^+|}\sum_{\substack{w\in W^1_P \\
I_{\nu}(w) = \Delta^+_P(e)}} (-1)^{\ell(w)} \text{Tr}(e^{-1};V^L_{w(\lambda_B+\rho_B)-\rho_B})
\end{equation*}
where $I_{\nu}(w)$ is defined in \parens{\ref{eqn-IWeyl}}.
\end{thm}

The proof will occupy the rest of this section.

\subsection{The nilmanifold correspondence}\label{subsec-nilmanifold}
The Hecke correspondence $C\rightrightarrows X$ extends to a correspondence  on the
Borel-Serre compactification 
\begin{equation}\label{eqn-BScorresp}\widetilde C \rightrightarrows \widetilde X\end{equation}
which is compatible with the projection $\mu:\widetilde X \to \overline{X}$ to the redutive
Borel-Serre compactification.  Let $w'\in F_P(e)$ and set $w = c_1(w')=c_2(w').$  The
restriction of the correspondence to the relevant Borel-Serre stratum is given by
\begin{align}
Y'_P = \Gamma'_P\backslash P/K_PA_P &\rightrightarrows Y_P = \Gamma_P\backslash P/K_PA_P \\
\Gamma'_PxK_PA_P &\mapsto (\Gamma_PxK_PA_P,\Gamma_PyxK_PA_P).\label{eqn-stratum}\end{align}
(Here, $\Gamma'_P = \Gamma_P \cap y^{-1}\Gamma_P y$.)  The fibers $N_P = \mu^{-1}(w) \subset
Y_P$ and $N'_P = (\mu')^{-1}(w')\subset Y'_P$ are nilmanifolds isomorphic to
$\Gamma_{\mathcal U}\backslash \mathcal U_P$ and $\Gamma'_{\mathcal U} \backslash \mathcal U$
respectively, where $\Gamma_{\mathcal U} = \Gamma_P \cap \mathcal U_P$ and $\Gamma'_{\mathcal
U} = \Gamma'_P \cap \mathcal U_P = \Gamma \cap y^{-1}\Gamma_{\mathcal U} y$.  So the
correspondence \ref{eqn-BScorresp} restricts to a correspondence $N'_P \rightrightarrows N_P$
which will be described below.  The following diagram may help in sorting out these spaces.
\begin{equation*}\begin{CD}
N'_P & \rightrightarrows & N_P &&& Y'_P & \rightrightarrows & Y_P &&& \widetilde{C} &
\rightrightarrows & \widetilde{X} \\
@VVV @VVV\text{\quad in\quad } & @VVV @VVV \text{\quad in \quad} & @VVV @VVV\\
w' & \rightrightarrows & w &&& C_P & \rightrightarrows & X_P &&& \overline{C} &
\rightrightarrows & \overline{X}
\end{CD}\end{equation*}

\begin{lem}\label{lem-nilman} Let $\phi:L_P \to GL(H^*(\mathfrak N_P,E))$ denote the adjoint
representation of the Levi quotient $L_P$ on the Lie algebra cohomology of $\mathfrak N_P.$  
Let $w' \in F_P(e)$ be a fixed point in $C_P$ with characteristic element $e \in L_P.$
Then the nilmanifold correspondence $(c_1,c_2):N'_P \rightrightarrows N_P$ induces a mapping
$(c_1)_*c_2^*:H^*(N_P,\mathbf E) \to H^*(N_P,\mathbf E)$ on cohomology which, under the
Nomizu-van Est isomorphism $H^*(N_P,\mathbf E) \cong H^*(\mathfrak N_P,E)$ may be identified
with the homomorphism
\begin{equation*} r\phi(e^{-1}) \end{equation*}
where $r = [\Gamma_{\mathcal U} : \Gamma'_{\mathcal U}].$ \end{lem}

\subsection{Proof}\label{pf-vanEst}  First we find equations for the nilmanifold
correspondence.  Choose a lift $xK_PA_P \in D=  P /K_PA_P$ of the fixed point
$w'=\Gamma'_PxK_PA_P\mathcal U_P \in C_P.$  This determines a parametrization of the
nilmanifold $N_P$ by
\begin{align}
\Gamma_{\mathcal U} \backslash \mathcal U & \longrightarrow N_P
\subset Y_P = \Gamma_P \backslash P / K_PA_P \\
\Gamma_{\mathcal U}z & \mapsto \Gamma_P zxK_PA_P \label{eqn-parametrization}\end{align}
and similarly $\Gamma'_{\mathcal U} \backslash \mathcal U \longrightarrow N'_P$ by
$\Gamma'_{\mathcal U}z \mapsto \Gamma'_P zx K_PA_P.$

Since $w'$ is fixed, we have $\Gamma_Px\mathcal U_PK_PA_P = \Gamma_Pyx\mathcal U_PK_PA_P$
hence there exists $\gamma \in \Gamma_P$ and $u\in \mathcal U_P$ so that $\gamma yuxK_PA_P =
xK_PA_P$, in other words, so that $\gamma yu$ fixes the point $xK_PA_P$ in the Borel-Serre
boundary component $P/K_PA_P.$  Then $e = \nu_P(\gamma y) = \nu_P(\gamma y u)$ is the
characteristic element of the fixed point $w'.$  Define $N'_P \rightrightarrows N_P$ by
\begin{equation}\label{eqn-nilmanifold}  \Gamma'_{\mathcal U} z \mapsto (\Gamma_{\mathcal U}z,
\Gamma_{\mathcal U}(\gamma y)zu^{-1}(\gamma y)^{-1} ). \end{equation}
A simple calculation shows that the following diagram commutes,
\begin{equation*} \begin{CD}
\Gamma'_{\mathcal U}\backslash \mathcal U & \overset{\ (\ref{eqn-nilmanifold})\ }
{\rightrightarrows} & \Gamma_{\mathcal U} \backslash \mathcal U\\
@V{(\ref{eqn-parametrization})}V{\cong}V @V{\cong}V{(\ref{eqn-parametrization})}V \\
N'_P & \rightrightarrows & N_P \\
@VVV @VVV \\
Y'_P & \overset{(\ref{eqn-stratum})}{\rightrightarrows}  & Y_P 
\end{CD}\end{equation*}

Next we will apply the theorem of Nomizu \cite{Nomizu-vanEst} and van Est \cite{vanEst} to
this correspondence.  The local system
$\mathbf E \to X$ which is defined by the representation $\tau:G \to GL(E)$ extends
canonically to a local system on the Borel-Serre compactification $\widetilde{X}$ and its
restriction to the nilmanifold $N_P$ is given by the quotient $\mathbf E | N_P = \mathcal U
\times_{\Gamma_{\mathcal U}}E$ under the relation $(z,v) \sim (\gamma z, \tau(\gamma)v)$
(for $\gamma \in \Gamma_{\mathcal U}$, $z\in \mathcal U_P$, and $v\in E$).  The complex
$\Omega^{\bullet}(N_P, \mathbf{E})$ of smooth $\mathbf E$-valued differential forms on $N_P$
consists of sections of the (flat) vectorbundle
\[ \mathbf{C}^{\bullet}(N_P,\mathbf{E}) = \mathcal U_P \times_{\Gamma_P} C^{\bullet}
(\mathfrak N_P,E)\]
where $C^{\bullet}(\mathfrak N_P,E) = \text{Hom}_{\mathbb C}(\wedge^{\bullet} \mathfrak
N_P,E)$ is the complex of Lie algebra cochains.  Let $\phi$ be the representation of $P$ on
this complex:  if $\wedge^{\bullet} \text{Ad}(p): \wedge^{\bullet}\mathfrak N_P \to
\wedge^{\bullet}\mathfrak N_P$ denotes the adjoint action of $p \in P$ on the exterior algebra
of $\mathfrak N_P$, then
\[ \phi(p)(s) = \tau(p) \circ s \circ \wedge^{\bullet} \text{Ad}(p). \]
Denote by
\begin{equation}\label{eqn-flatforms}
 \Omega^{\bullet}_{\text{inv}}(\mathcal U_P,E) = 
\left\{ \omega: \mathcal U_P \to C^{\bullet}(\mathfrak N_P,E) \left| \right.
\omega(ux) = \phi(u) \omega(x) \text{ for all } u,x \in \mathcal U_P \right\} \end{equation}
the complex of (left) $\mathcal U_P$-invariant $E$-valued differential forms on $\mathcal
U_P.$  Such a differential form is determined by its value $s = \omega(1) \in
C^{\bullet}(\mathfrak N_P,E)$, and it passes to a differential form on $N_P.$  Denote by
$\Omega^{\bullet}_{\text{inv}}(N_P, \mathbf{E})$ the collection of all such ``left''-invariant
differential forms.  The Nomizu-van Est theorem (\cite{Nomizu-vanEst, vanEst})states that the
inclusion $\Omega^{\bullet}_{\text{inv}}(N_P,\mathbf{E}) \hookrightarrow \Omega^{\bullet}
(N_P,\mathbf{E})$ induces an isomorphism on cohomology.  In summary we have a diagram
\begin{equation*}\begin{CD}
C^{\bullet}(\mathfrak N_P,E) @<{\cong}<< \Omega^{\bullet}_{\text{inv}}(\mathcal U_P,E)
@<{\cong}<< \Omega^{\bullet}_{\text{inv}}(N_P,\mathbf{E}) \hookrightarrow
\Omega^{\bullet}(N_P,\mathbf{E})
\end{CD} \end{equation*}
of isomorphisms and quasi-isomorphisms.
Although the group $P$ does not act on the vectorbundle $\mathbf{C}^{\bullet}(N_P,
\mathbf{E})$, it does act on the complex $\Omega^{\bullet}_{\text{inv}}(N_P, \mathbf{E}) \cong
\Omega^{\bullet}_{\text{inv}}(\mathcal U_P,E)$ of invariant sections by
\[ (p \cdot \omega)(x) = \phi(p)^{-1}\omega(pxp^{-1})\]
and the group $\mathcal U_P$ acts on this complex by
\[ (u\cdot \omega)(x) = \omega(xu^{-1}).\]
If $\omega \in \Omega^{\bullet}_{\text{inv}}(N_P,\mathbf{E})$ is given by
(\ref{eqn-flatforms}) then by (\ref{eqn-nilmanifold}) its pullback by $c_2$ is given by
\[ c_2^*(\omega)(z) = \phi(\gamma y)^{-1} \omega((\gamma y) zu^{-1} (\gamma y)^{-1}).\]
Evaluating at $z=1$ and using the fact that $\omega$ is left invariant,
\begin{equation*} c_2^*(\omega)(1)  = \phi(u)^{-1} \phi(\gamma y)^{-1}\omega(1).
\end{equation*}
Let $s=\omega(1) \in C^{\bullet}(\mathfrak N_P,E)$, suppose $ds=0$ and let $[s] \in
H^*(\mathfrak N_P,E)$ be the resulting cohomology class.  Since $\mathcal U_P$ acts trivially
on this cohomology,
\[ c_2^*([s]) = \phi(e)^{-1}[s] \]
where $e = \nu_P(\gamma y)$ is the characteristic element of the fixed point $w$.  Finally,
observe that the pushforward mapping $(c_1)_*:H^*( N'_P,\mathbf E) \to H^*(N_P,\mathbf E)$ is
given by multiplication by $r = [\Gamma_{\mathcal U}:\Gamma'_{\mathcal U}].$  This completes
the proof of lemma \ref{lem-nilman}.  \qed

\subsection {Proof of Theorem \ref{thm-answer}} \label{subsec-outline}
We will apply the Lefschetz fixed point formula to the modified Hecke correspondence. 
By Proposition \ref{prop-nofixedpt}, after modifying the correspondence by composing with
$Sh(\mathbf t)$, the fixed point constituent $F_P(e)$ becomes ``truncated'', that is, it
becomes replaced by the intersection $F_P^0(e) = F_P(e) \cap C_P^0$ of $F_P(e)$
with the central tile in $C_P.$  Denote by $\partial F^0 = F_P(e) \cap \partial C_P^0$ its
intersection with the boundary of the central tile.  Set $F' = c_1(F_P(e)) = c_2(F_P(e)).$ Set
$E = F'\cap X_P^0 = c_i(F_P^0(e))$ and $\partial E = F' \cap \partial X_P^0 = c_i(\partial
F^0).$  (Having used up all the letters some time ago, we temporarily re-use the notation $E$
here, hoping the reader will not confuse it with the local system.)  Note that $E - \partial
E$ is diffeomorphic to $F'.$

By theorem \ref{thm-hyperbolic} the (modified) Hecke correspondence is weakly hyperbolic near
$F_P(e)$ and an indicator mapping (defined in a neighborhood $U\subset \overline{X}$ of
$F'$) is given by
\begin{equation}\label{eqn-indicator}
t(x) = \Bigl(\sum_{\alpha \in \Delta^+_P}r_{\alpha}^P(x) + \ d_E\pi_P(x), \sum_{\alpha \in
\Delta^-_P} r_{\alpha}^P(x) + \sum_{\alpha \in
\Delta^0_P}r_{\alpha}^P(x) \Bigr)
\end{equation}
Let $\mathbf Q \supset \mathbf P$ be the rational parabolic subgroup corresponding to the
subset $I = \Delta_P^+  \subset \Delta_P$ consisting of the simple roots for
which the Hecke correspondence is (strictly) expanding.  Then, in the notation of
(\ref{eqn-decomposition}), $ \Delta_P = i(\Delta_Q) \amalg \Delta_P^+ .$
The partial distance function $r_{\alpha}^P$ vanishes on the stratum $X_Q$ whenever
$\alpha \in \Delta_P^{-}\cup \Delta_P^0$, cf. (\ref{eqn-rootpi2}).  Hence $\overline{X}_Q\cap
U = t^{-1}(\mathbb R_{\ge 0} \times \{0\})$ is the ``expanding set'' of the correspondence.

According to Theorem \ref{thm-topological} we need to compute the stalk cohomology (at points
$w\in E$) of the sheaf 
\[ \mathbf{A^{\bullet}} = h^{!}j^{*}\mathbf{W^{\nu}C^{\bullet}(E)} \]
where 
\begin{equation*}\begin{CD}
E @>>{h}> \overline{X}_Q @>>{j}> \overline{X}.
\end{CD}\end{equation*}
This is best accomplished by decomposing $h$,
\begin{equation*}\begin{CD}
E @>>{h_1}> F'  @>>{h_2}> X_P @>>{h_3}>  \overline{X}_Q @>>{j}> \overline{X}
\end{CD}\end{equation*}
Then $\mathbf B^{\bullet} = h_3^!j^*\mathbf{W^{\nu}C^{\bullet}(E)}$ is the sheaf studied in
Theorem \ref{thm-supports}, where we have taken $I = \Delta_P^+.$  Its stalk cohomology is
locally constant on $X_P$ and was shown to be 
\begin{equation*}
H^i_{w}(\mathbf B^{\bullet}) \cong H^{j - |\Delta_P^+|}(\mathfrak N_P,\mathbf
E)\qdr{\Delta_P^+}
\end{equation*}  
Since $h_2$ is a smooth closed embedding we have a canonical quasi-isomorphism
(\ref{eqn-normal})
\begin{equation*}
\mathbf C^{\bullet} := h_2^!\mathbf B^{\bullet} \cong h_2^*(\mathbf B^{\bullet}) \otimes
\mathcal O[-c]
\end{equation*}
where $c = \dim(X_P) - \dim(F' \cap X_P)$ and where $\mathcal O$ is the orientation bundle
(i.e. the top exterior power) of the normal bundle of $F'\cap X_P$ in $X_P.$  The complex
$\mathbf C^{\bullet}$ is constructible with respect to the stratification of $\overline{X}$,
meaning that its cohomology sheaves are locally constant on $X_P,$ hence also on $E$.  But $E$
is a manifold with boundary, so
\begin{equation}\label{eqn-shriek}
 h_1^!\mathbf C^{\bullet} \cong i_!\mathbf C^{\bullet}|(E- \partial E)\end{equation}
is obtained by first restricting to the interior $E- \partial E$ and then extending by 0.
(Here, $i:E - \partial E \to E$ denotes the inclusion.)  
 Thus the cohomology of $h_1^!\mathbf C^{\bullet}$ is the compactly supported cohomology
$H^i_c(E-\partial E;\mathbf C^{\bullet})\cong H^i_c(F'; \mathbf C^{\bullet}).$

Next we must compute the pointwise Lefschetz number $n(w)$ for $w\in E,$ that is, the
alternating sum of the traces on the stalk cohomology of $\mathbf {A^{\bullet}} = h_1^!
\mathbf C^{\bullet}.$   By (\ref{eqn-shriek}) it is 0 when $w\in \partial E$, so let $w\in E -
\partial E.$  Then
\begin{align}
H^i_w(\mathbf C^{\bullet}) &= H^{i-c}_w(h_2^*\mathbf B^{\bullet} \otimes\mathcal O) \\
&= H^{i-c-|\Delta_P^+|}(\mathfrak N_P,E)\qdr{\Delta_P^+} \otimes \mathcal O_w.
\label{eqn-stalk cohomology}
\end{align}

By \S \ref{subsec-addendum}, the mapping $c_1:F_P(e) \to F'$ is a covering of degree
$d = [ \Gamma_L \cap \bar y^{-1} \Gamma_L \bar y : \nu_P(\Gamma_P \cap y^{-1} \Gamma_P y)].$
Near each fixed point $w' \in c_1^{-1}(w)$ the Hecke correspondence acts on the
$\mathfrak N_P$-cohomology through the homomorphism $r\phi(e^{-1})$ (using Lemma
\ref{lem-nilman}),  and by \S \ref{subsec-orientation} it acts on $\mathcal O_w$ by $(-1)^c$.
Summing these contributions over the $d$ different points in $c_1^{-1}(w)$ gives

\begin{align}\label{eqn-dr}
n(w) &= dr(-1)^{-c-|\Delta_P^+|}(-1)^c\sum_{i\ge 0}(-1)^i \text{Tr}(\phi(e^{-1});
H^i(\mathfrak N_P,\mathbf E)\qdr{\Delta_P^+}) \\
&= dr (-1)^{|\Delta_P^+|} \sum_{\substack{v\in W^1_P \\
I_{\nu}(v) = \Delta^+_P}} (-1)^{\ell(v)} \text{Tr}(e^{-1};V^L_{v(\lambda_B+\rho_B)-\rho_B})
\label{eqn-mess}
\end{align}
by Proposition \ref{prop-Kostant}.  The contribution arising from $F_P(e)$ is this quantity
times $\chi_c(E-\partial E) = \chi_c(F').$  However (by \S \ref{subsec-addendum}),
$\chi_c(F_P(e)) = d \chi_c(F')$ which absorbs the factor of $d$ in (\ref{eqn-mess}) and
therefore completes the proof of  Theorem \ref{thm-answer}.  \qed

\section{Proof of Theorem \ref{thm-intro}} 
 \subsection{}  As in \S \ref{sec-notation}, $\mathbf G$ denotes a connected reductive linear
algebraic group defined over $\mathbb Q$, $D = G/K'$ is its associated symmetric space with
basepoint $x_0 \in D$ and stabilizer $K' = A_GK(x_0).$ Let $\Gamma \subset \mathbf G(\mathbb
Q)$ denote an arithmetic subgroup which we assume to be neat, and $X = \Gamma \backslash
D$.  Throughout this section we fix a Hecke correspondence $(c_1,c_2):\overline{C}
\rightrightarrows \overline{X}$ defined by some element $g\in \mathbf G(\mathbb Q).$  So $C =
\Gamma'\backslash D$ with $\Gamma' = \Gamma \cap g^{-1}\Gamma g.$  We also fix a
$\Gamma$-equivariant tiling of $D$ which is narrow with respect to the Hecke correspondence. 
Choose $\mathbf t \in A_{P_0}(\agt 1)$ in accordance with Proposition \ref{prop-nofixedpt}.

 Let $F \subset \overline{C}$ denote the (full) fixed point set of the Hecke correspondence
$\overline{C} \rightrightarrows \overline{X}$ and let $E$ denote the (full) fixed point set of
the modified Hecke correspondence (\ref{eqn-modified}). Then
\[ F = \coprod_{\{\mathbf P\}}F \cap C_P \text{ and } E = \coprod_{\{\mathbf P\}}F \cap
C_P^0\]
where the union is over the strata of $\overline{C}$, that is, over $\Gamma'$-conjugacy
classes of rational parabolic subgroups $\mathbf P \subseteq \mathbf G.$  Each $F \cap C_P^0$
is a union of connected components of $E$ by Proposition \ref{prop-nofixedpt}.  The Lefschetz
fixed point theorem (Theorem \ref{thm-topological}) may be used to write the Lefschetz number
as a sum over these individual strata.

\subsection{Contribution from a single stratum}
Let $\mathbf P \subseteq \mathbf G$ be a rational parabolic subgroup and suppose that
$c_1(C_P) = c_2(C_P) = X_P.$  By Proposition \ref{prop-restriction}, in a neighborhood of
$C_P$ the correspondence is isomorphic to the parabolic Hecke correspondence determined by
some $y \in \Gamma g \Gamma \cap P$ and moreover (in this neighborhood) the fixed points of
the modified correspondence coincide with those of $E \cap C_P=F\cap C_P^0.$ 

If $F_P(e)$ denotes the set of fixed points in $C_P$ with characteristic element $e\in
\Gamma_L \bar y \Gamma_L$, then by Proposition \ref{prop-char},
\begin{equation}\label{eqn-conjclass} 
F\cap C_P = \coprod_{\{e\}} F_P(e) \text{ and } E \cap C_P = \coprod_{\{e\}} F_P(e) 
\cap C_P^0 \end{equation}
where the union is over $\Gamma_{L}-$conjugacy classes of elements $\{e\} \subset \Gamma_L
\bar y \Gamma_L$ which are elliptic modulo $A_P.$  (Here, $\Gamma_L = \nu_P(\Gamma \cap P)
\subset L_P$ and $\bar y = \nu_P(y).$)
For each such conjugacy class $\{e\}$, the set $F_P(e)$ consists of
finitely many connected components, say, $F_1, F_1,\ldots, F_m.$  The contribution to the
Lefschetz number from the component $F_j$ is given by Theorem \ref{thm-answer}.  By
(\ref{eqn-centralizer}) (see also \S \ref{subsec-addendum}, \S \ref{subsec-Euler} ),
\begin{equation*}
 \sum_{j=1}^m \chi_c(F_j)=\chi_c(F_P(e)) =  \chi_c (\Gamma'_e \backslash
L_e/K'_e).\end{equation*}
So the contribution to the Lefschetz number from the stratum $C_P \rightrightarrows
X_P$  is \begin{equation}\label{eqn-LPy}
L(P,y):=\sum_{\{e\}} \chi_c(\Gamma'_e \backslash L_e /K'_e)  r
(-1)^{|\Delta_P^{+}|}\sum_{\substack{w\in W^1_P \\
I_{\nu}(w) = \Delta^+_P(e)}} (-1)^{\ell(w)} \text{Tr}(e^{-1};V^L_{w(\lambda_B+\rho_B)-\rho_B})
\end{equation}
where the index set for the first sum is the same as that for the union in
(\ref{eqn-conjclass}).
This quantity $L(P,y)$ depends only on the local system $\mathbf E$, the choice of parabolic
subgroup $\mathbf P$ and the element $y \in P.$

\subsection{Sum over strata}
Let $\mathbf P_1, \mathbf P_2, \ldots, \mathbf P_t$ denote a collection of representatives,
one from each $\Gamma$-conjugacy class of rational parabolic subgroups $\mathbf P \subseteq
\mathbf G.$  These index the strata of $\overline{X}.$  For each such $i$ the intersection
$\Gamma g \Gamma \cap P_i$ decomposes:
\[ \Gamma g \Gamma \cap P_i = \coprod_j \Gamma_{P_i} y_{ij} \Gamma_{P_i} .\]
 Lemma \ref{lem-Xi} gives a one-to-one correspondence between this collection $\{ y_{ij}  \}$
and strata $C_{ij}$ of $\overline{C}$ such that $c_1(C_{ij})=c_2(C_{ij}).$  Moreover the
restriction of the Hecke correspondence to a neighborhood of $C_{ij}$ is locally isomorphic to
the parabolic Hecke correspondence defined by $y_{ij}$ so the local contribution to the
Lefschetz number from $C_{ij}$ equals the number $L(P_i,y_{ij})$ given in (\ref{eqn-LPy}).  In
summary, the total Lefschetz number is
\begin{equation}\label{eqn-finalsum}
L(g)= \sum_{i=1}^t \sum_j L(P_i,y_{ij})\end{equation}
as claimed in Theorem \ref{thm-intro}.  \qed   

\subsection{Another formula}
  If a little expansion $Sh(\mathbf t)^{-1}$ is used instead of the shrink, this will convert
neutral directions normal to each stratum into expanding directions, and it will convert the
tangential distance into a contracting direction.  An indicator mapping replacing
(\ref{eqn-indicator}) is
\begin{equation*}
t(x) = \bigl( \sum_{\alpha \in \Delta_P^+ \cup \Delta_P^0}r_{\alpha}^P(x),\ \sum_{\alpha \in
\Delta_P^-} r_{\alpha}^P(x) + d_E\pi_P(x)
\bigr)
\end{equation*}
This changes the nature of the sheaf $\mathbf A^{\bullet}$ with the result that the Euler
characteristic (rather than the Euler characteristic with compact supports) appears in the
formula.  So, in equation (\ref{eqn-finalsum}), the contribution $L(P,y)$ (\ref{eqn-LPy}) from
the stratum $C_P \rightrightarrows X_P$ will be replaced by the quantity
\begin{equation*}
L'(P,y)=\sum_{\{e\}}
r(-1)^{|\Delta_P^+ \cup \Delta_P^0|}\chi(\Gamma'_e \backslash L_e /K'_e)  
\sum_{\substack{w\in W^1_P \\
I_{\nu}(w)=\Delta_P^+\cup\Delta_P^0}} (-1)^{\ell(w)} \text{Tr}\bigl(e^{-1}; 
V^L_{w(\lambda_B + \rho_B) - \rho_B}\bigr)
\end{equation*}
where the summations are over the same index sets as in (\ref{eqn-LPy}), and where $\Delta_P^+
= \Delta_P^+(e)$ and $\Delta_P^0 = \Delta_P^0(e).$
\section{Remarks on the Euler characteristic}\label{sec-Harder}
As in \S \ref{sec-notation}, $\mathbf G$ denotes a connected linear reductive algebraic group
defined over $\mathbb Q$; $D$ denotes the associated symmetric space;  $\mathbf{S_G}$ denotes
the greatest $\mathbb Q$-split torus in the center of $G$; $A_G = \mathbf{S_G}(\mathbb R)^0$
denotes the identity component of its real points; $K' = A_GK$ is the
stabilizer in $G$ of a fixed basepoint $x_0\in D$; $\Gamma \subset \mathbf G(\mathbb Q)$
is an arithmetic group, and $X = \Gamma \backslash D.$  Recall the following classical result
of Harder \cite{Harder} (a more streamlined proof of which may be found in \cite{Leuzinger}).

\begin{thm}\label{prop-Harder}
  The Euler characteristic $\chi(X)$ is given by the integral over $X$ of the
Euler form with respect to any invariant Riemannian metric on $X.$  \qed\end{thm}
  For completeness we also include a proof of the following often-cited fact.

\begin{lem}\label{lem-Harder}
 Suppose the real Lie group $G/A_G$ does not contain a compact maximal torus.  Then the Euler
form vanishes identically on $X$. \end{lem}
\subsection{Proof} By replacing $\mathbf G$ by the algebraic group $\mathbf{{}^0G}$ (and
noting that $X = \Gamma\backslash {}^0G /K$), we may  assume that
$\mathbf{S_G}$ is trivial.  Let $\mathfrak g = \mathfrak k \oplus \mathfrak p$ be the Cartan
decomposition of $\text{Lie}(G)$ corresponding to the choice $K$ of maximal compact subgroup.
 Choose a $K$-invariant inner product on $\mathfrak p.$  This determines a $G$-invariant
Riemannian metric on $D = G/K$ which passes to a Riemannian metric on $X$.  Let $\Omega$ be
the curvature form of the torsion-free Levi-Civita connection which is associated to this
metric.  The resulting Euler
form $Eu$ is defined to be $0$ if $\dim(D)$ is odd.  If $\dim(D)=2k$ then $Eu$ is the
$G$-invariant differential form on $D$ whose value at the basepoint $x_0$ is $Eu(\dot x_1,
\dot y_1, \dot x_2, \dot y_2, \ldots, \dot x_k, \dot y_k) = P(\Omega(\dot x_1,\dot
y_1),\ldots, \Omega(\dot x_k, \dot y_k))$ (for any $\dot x_1, \ldots, \dot y_k \in \mathfrak p
= T_{x_0}D)$, where $P$ is the polarization of the Pfaffian $\text{Pf}: \text{End}(\mathfrak
p)^- \to \mathbb R.$  (Here, $\text{End}(\mathfrak p)^-$ denotes the skew-adjoint
endomorphisms of $\mathfrak p.$)  The form $Eu$ on $D$ passes to a differential form on $X =
\Gamma \backslash D$, which is the Euler form for $X.$

Let $\text{Ad}: K \to \text{GL}(\mathfrak p)$ be the adjoint representation and let
$ad:\mathfrak k \to \text{End}(\mathfrak p)^-$ be its derivative.  We claim that $\det(ad(\dot
k)) = 0$ for any $\dot k \in \mathfrak k.$  Modify $\dot k$ by conjugacy if necessary, so as
to guarantee that $\dot k$ lies in a maximal torus $\mathfrak t \subset \mathfrak g$ which is
stable under the Cartan involution (\cite{Warner} \S 1.2, 1.3).  Then $\mathfrak t = \mathfrak
t_{+} \oplus \mathfrak t_{-}$ with $\mathfrak t_{+} \subset \mathfrak k$ and $\mathfrak t_{-}
\subset \mathfrak p.$  By assumption, $\mathfrak t_{-}$ contains a nonzero vector $\dot t$,
and $ad(\dot k)(\dot t) = [\dot k, \dot t] = 0,$ which proves the claim.

The principal $K$-bundle $G \to D = G/K$ admits a canonical $G$-invariant connection
(\cite{KN} Chapt. II Thm. 11.5).  Its curvature form $\omega \in \mathcal A^2(D, \mathfrak k)$
is the $G$-invariant differential form whose value at the basepoint $x_0$ is given by
$\omega_{0}(\dot p_1, \dot p_2) = -[\dot p_1, \dot p_2] \in \mathfrak k$ for any $\dot p_1,
\dot p_2 \in \mathfrak p.$   By a theorem of Nomizu \cite{Nomizu-invariant}, for any real
representation $\lambda: K \to \text{GL}(E)$, the resulting connection in the associated
$G$-homogeneous vectorbundle $\mathbf E = G \times_K E$ coincides with the torsion-free metric
(Levi-Civita) connection of any $G$-invariant metric on $E$. Its curvature is the
$G$-invariant $\text{End}(\mathbf E)$-valued differential form whose value at the basepoint is
$\Omega_0 = \lambda' \circ \omega_0$ where $\lambda': \mathfrak k \to \text{End}(E)$ is the
differential of $\lambda.$  Taking $\lambda = \text{Ad}:K \to \text{GL}(\mathfrak p)$ as above
gives $\Omega_0(\dot p_1, \dot p_2) = - ad([\dot p_1,\dot p_2])$.  By the above claim,
this has determinant 0 hence its Pfaffian vanishes also. Therefore the Euler form is zero on
$D$, so it is also zero on $X.$
\qed

\begin{lem}\label{lem-cptspt}  Let $X = \Gamma \backslash G /KA_G$ as above.
  Then the Euler characteristic and the Euler characteristic with compact supports
coincide:  $\chi(X) = \chi_c(X).$ \end{lem}
\subsection{Proof}  Let $\widetilde{X}$ denote the Bore-Serre compactification of $X$.
Topologically, it is a manifold with boundary $\partial \widetilde{X} = \widetilde{X} - X.$
Since $H^i_c(X) = H^i(\widetilde{X}, \partial \widetilde{X})$, it suffices to show that
$\chi(\partial \widetilde{X}) = 0.$  The boundary $\partial \widetilde{X}$ is a union of
finitely many boundary strata $Y_P,$ each of which fibers over the corresponding stratum $X_P$
(of the reductive Borel-Serre compactification) with fiber a nilmanifold $N_P$ (cf. \S
\ref{subsec-boundary components}, \ref{subsec-nilmanifold}).  So $\chi(Y_P) =
\chi(N_P)\chi(X_P) = 0.$  It follows from Mayer-Vietoris that $\chi(\partial \widetilde{X}) =
0.$  \qed

The above results together give:
\begin{prop}\label{thm-Harder}  
Suppose $(\mathbf G/\mathbf{S_G})(\mathbb R)$ does not contain a compact maximal torus.  Then
$\chi(X) = \chi_c(X)=0,$ that is, both the Euler characteristic and the Euler
characteristic with compact supports vanish. \qed\end{prop}

\subsection{Euler characteristic of a fixed point component}\label{subsec-Euler}
Now suppose that $X_P \subset \overline{X}$ is a boundary stratum corresponding to a rational
parabolic subgroup $ P = \mathcal U_P L_P.$  Let $F_P(e) \subset X_P$ be the set of fixed
points with some fixed (elliptic) characteristic element $e \in \mathbf{L_P}(\mathbb Q).$ Let
$L_e$ be the centralizer of $e$ in $L_P.$  By (\ref{eqn-centralizer}), $F_P(e) \cong \Gamma'_e
\backslash L_e / K'_e$ where $\Gamma'_e = \Gamma'_L \cap L_e$ and where $K'_e = L_e \cap
( z (K_PA_P) z^{-1})$ (for appropriate $z$).  By \S \ref{subsec-addendum},
\begin{equation}\label{eqn-dEuler}
 \chi_c(F_P(e)) = d \chi_c(\Gamma_e \backslash L_e /K'_e) \end{equation}
where $\Gamma_e = \Gamma_L \cap L_e$ and $d = [\Gamma_L:\Gamma'_L].$  This expression has the
following merit.  The
contribution (\ref{eqn-LPy})  to the Lefschetz number from the stratum $C_P$ depends on the
subgroup $\Gamma'_P \subset \Gamma_P.$ However once this expression (\ref{eqn-dEuler}) has
been substituted into (\ref{eqn-LPy}), the dependency on this subgroup $\Gamma'_P$ occurs only
in the two integers $r$ and $d$.

Let $\mathbf S_e$ be the greatest
$\mathbb Q$-split torus in the center of $L_e$ and let $A_e$
be the identity component of its group of real points.  As explained in \cite{GKM} \S 7.11,
the group $K'_e$ does not necessarily contain $\mathbf{S_e}(\mathbb R)$, so although
$F_P(e)$ is not necessarily a ``locally symmetric space'' in the sense of \S
\ref{sec-notation}, it fibers over the locally symmetric space $\Gamma'_e \backslash L_e /
K'_eA_e$ with fiber $A_e/A_P$ which is diffeomorphic to a Euclidean space.  Therefore
\[ \chi_c(F_P(e)) = (-1)^{\dim(A_e/A_P)} \chi_c(\Gamma'_e \backslash L_e /K_eA_e)
= (-1)^{\dim(A_e/A_P)} d \chi_c(\Gamma_e \backslash L_e / K_e A_e)\]
(where $K_e = L_e \cap (zK_Pz^{-1})$).

Now suppose that $L_P/A_P$ does not contain a compact maximal torus.  Then the same is true of
$L_e/A_P$.  If $A_e=A_P$ then it follows from Proposition \ref{thm-Harder} that $\chi(F_P(e))
= \chi_c(F_P(e)) = 0.$  However it is not the case in general that $A_e=A_P$.  It is possible
that $L_e/A_e$ does contain a compact maximal torus and that $\chi(F_P(e)) \ne 0$ (in which
case $F_P(e)$ is necessarily non-compact since it is fibered by $A_e/A_P$ as described above).
See Example \ref{subsec-notzero} in which $L_e = A_e$ and $F_P(e) \cong A_e/A_P$ is the orbit
of a split torus.  In such cases it is possible to re-attribute the contribution
(\ref{eqn-LPy}) (to the Lefschetz number) from the stratum $C_P \rightrightarrows X_P$ to
smaller strata in the correspondence.  This procedure is carried out in \cite{GKM} p. 531,
resulting in a Lefschetz formula in which the only nonzero contributions come from strata $C_P
\rightrightarrows X_P$ such that $L_P$ has a compact maximal torus.

In the adelic setting, the Euler characteristic with compact support $\chi_c(F_P(e))$ can be
expressed in terms of orbital integrals (cf 
\cite{GKM} \S 7.11 and \S 7.14).

\section{Examples and special cases}
\subsection{Reducible fixed point components} \label{subsec-reducible}
For $G = SL(3, \mathbb R)$, $D = G/K$, and $\Gamma \subset SL(3, \mathbb Z)$ a neat principal
congruence subgroup, the reductive Borel-Serre compactification $\overline{X}$ contains a
singular 0-dimensional stratum $X_B$ corresponding to the standard Borel subgroup $\mathbf B.$
This stratum is contained in the closures of the strata $X_{P_1}$ and $X_{P_2}$ corresponding
to the standard maximal parabolic subgroups $\mathbf{P_1}, \mathbf{P_2}$ containing $\mathbf
B.$  Let $g$ be a generic element of $\mathbf{G}(\mathbb Q) \cap \mathcal U(P_1) \cap
\mathcal U(P_2)$ which is not in $\Gamma$, for example, 
\[ g = \left( \begin{matrix} 1 & 0 & \frac{1}{2} \\
0 & 1 & 0 \\ 0 & 0 & 1 \end{matrix} \right).\]
Let $\overline{C} \rightrightarrows
\overline{X}$ be the resulting Hecke correspondence.  Then these three strata $X_{P_1},
X_{P_2},$ and $X_B$ are fixed by the correspondence.  However points in $X$ which are
sufficiently close to these strata are not fixed.

\subsection{Middle weight for $\mathbf{Sp}_4$}
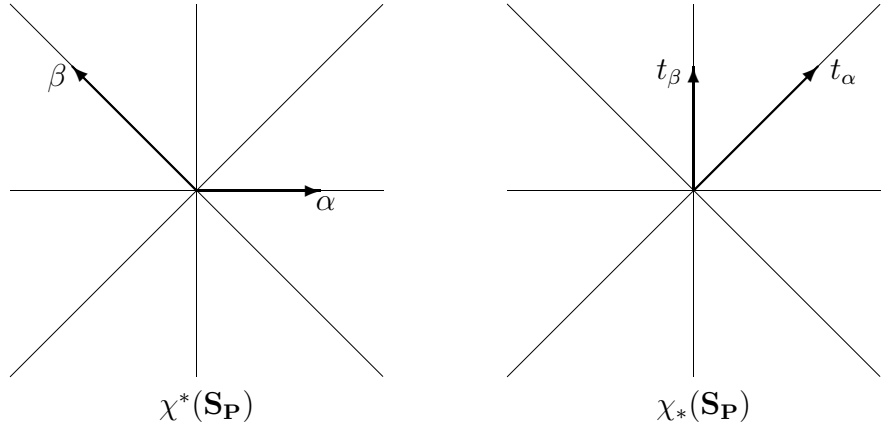
\begin{figure}[h!]
\begin{picture}(350,180)(-50,-20)  %% (450,220)(0,0)
\setlength{\unitlength}{.08em}   %% {.1em} is the normal size
\allinethickness{0.1pt}

\put(0,75){\line(1,0){150}} %% first weyl diagram
\put(75,0){\line(0,1){150}}
\put(0,0){\line(1,1){150}}
\put(0,150){\line(1,-1){150}}

\put(200,75){\line(1,0){150}}  %% second weyl diagram
\put(275,0){\line(0,1){150}}
\put(200,0){\line(1,1){150}}
\put(200,150){\line(1,-1){150}}

\allinethickness{0.5pt}
\put(75,75){\vector(1,0){50}} \put(123,67){$\alpha$}
\put(275,75){\vector(0,1){50}} \put(260,120){$t_{\beta}$}
\allinethickness{1.0pt}
\put(75,75){\vector(-1,1){50}} \put(15,117){$\beta$}
\put(275,75){\vector(1,1){50}} \put(330,120){$t_{\alpha}$}

\put(60,-15){$\chi^*(\mathbf{S_P})$}
\put(260,-15){$\chi_*(\mathbf{S_P})$}

\end{picture}\caption{Simple roots and dual basis}
\label{fig-roots}\end{figure}

Let $\mathbf G = \mathbf{Sp}_4$, fix a neat arithmetic subgroup $\Gamma \subset \mathbf
G(\mathbb Q),$  and choose a Hecke correspondence $\overline{C} \rightrightarrows
\overline{X}$ which is determined by some $g \in \mathbf G(\mathbb Q).$  If $\mathbf P$ is a
minimal parabolic subgroup of $\mathbf G$ then its Levi quotient $\mathbf{L_P} = \mathbf{S_P}$
is a maximal split torus and the boundary stratum $C_P$ consists of a single point.  Suppose
this point is an isolated fixed point of the Hecke correspondence.  Let $e$ be its
characteristic element.  The vectorspace $\chi_{\mathbb Q}^*(\mathbf{S_P})$ has a basis
consisting of the simple roots $\Delta_P = \{ \alpha,\beta \}$.  Let $\{ t_{\alpha}, t_{\beta}
\}$ be the dual basis of $\chi_*^{\mathbb Q}(\mathbf{S_P})$, so that $\langle \beta, t_{\beta}
\rangle = 1$, $\langle \beta, t_{\alpha} \rangle = 0$ and the same with $\alpha$ and $\beta$
interchanged.  See Figure \ref{fig-roots}.

\begin{figure}[h!]

\begin{picture}(350,180)(-50,-20)  %% (450,220)(0,0)
\setlength{\unitlength}{.08em}   %% {.1em} is the normal size
%% first weyl diagram

\allinethickness{1.5pt}
\put(0,75){\line(1,0){150}}
\put(0,150){\line(1,-1){150}}
\allinethickness{0.1pt}
\put(75,0){\line(0,1){150}}
\put(0,0){\line(1,1){150}}

%%% second weyl diagram %%%

\put(200,75){\line(1,0){150}}  
\put(200,150){\line(1,-1){150}}
\allinethickness{1.5pt}
\put(275,0){\line(0,1){150}}
\put(200,0){\line(1,1){150}}
\allinethickness{.1pt}

\put(75,75){\vector(1,0){50}}   %% \put(123,67){$\alpha$}
\put(275,75){\vector(0,1){50}}  %% \put(263,120){$t_{\beta}$}

\put(75,75){\vector(-1,1){50}}  %% \put(15,117){$\beta$}
\put(275,75){\vector(1,1){50}}  %% \put(330,120){$t_{\alpha}$}

%%%  first weyl diagram:  dots and labels  %%%
\newcommand{\ci}{\circle*{2.5}}
\put(50,0){\ci}
\put(100,0){\ci}
\put(150,50){\ci}
\put(150,100){\ci}
\put(100,150){\ci}
\put(50,150){\ci}
\put(0,100){\ci}
\put(0,50){\ci}

\put(50,135){$\phi$}
\put(95,135){$\phi$}
\put(130,97){$\phi$}
\put(120,48){$\{\beta\}$}
\put(85,10){$\{\alpha,\beta\}$}
\put(30,10){$\{\alpha,\beta\}$}
\put(10,48){$\{\alpha,\beta\}$}
\put(10,97){$\{\alpha\}$}

%%% second weyl diagram:  labels  %%%

\put(297,135){$\phi$}
\put(320,48){$\{\beta\}$}
\put(320,97){$\{\beta\}$}
\put(290,10){$\{\beta\}$}
\put(230,10){$\{\alpha,\beta\}$}
\put(210,48){$\{\alpha\}$}
\put(210,97){$\{\alpha\}$}
\put(242,135){$\{\alpha\}$}

\put(60,-20){$\chi^*(\mathbf{S_P})$}
\put(260,-20){$\chi_*(\mathbf{S_P})$}

\end{picture}
\caption{Diagram of $I_{\nu}(w)$ and of $\Delta_P^{+}(e)$}\label{fig-weights}
\end{figure}
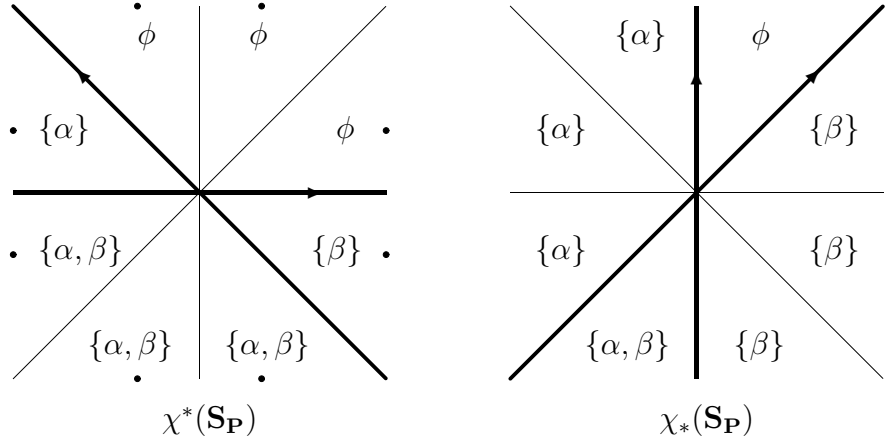

Let us take $\mathbf E$ to be the trivial local system, and the weight profile $\nu = -\rho_B$
to be the middle weight (where $\rho_B$ is the half-sum of the positive roots).  The
cohomology $H^*(\mathfrak N_P, \mathbb C)$ decomposes into a sum of 1-dimensional weight
spaces,
\[ V_{w\rho_B - \rho_B} \subset H^{\ell(w)}(\mathfrak N_P, \mathbb C) \]
as $w\in W$ varies over the elements of the full Weyl group.
These weights are the dots in the left hand part of Figure \ref{fig-weights}, in which the
origin is at $-\rho_B.$  For each weight space indexed by a given $w \in W$ we have indicated
the corresponding set 
\[ I_{\nu}(w) = \{ \theta \in \Delta_P \left| \right. \langle w\rho_B, t_{\theta} \rangle < 0
\} \]
of simple roots.  The cohomology $H^*(\mathfrak N_P)$ is divided into four ``quadrants''
according to the value of $I_{\nu}(w).$  
%%The cohomology degree of each weight space is indicated in parentheses.%%

If necessary, project the characteristic element $e$ to the identity component $A_P$ of the
torus $\mathbf{S_P}(\mathbb R)$ and let $t \in \mathfrak A_P = \text{Lie}(A_P)$ denote its
log.  The right hand half of Figure \ref{fig-weights} may be identified with the Lie algebra
$\mathfrak A_P.$  The chamber containing $t$ determines the expanding-contracting nature of
the Hecke correspondence near this fixed point.  In each chamber we have indicated the set
of expanding roots,
\[ \Delta^{+}_P(e) = \{ \theta \in \Delta_P \left| \right. \theta(t) \alt 0\} \]
(where now $\theta \in \Delta_P$ has been identified with a homomorphism $\mathfrak A_P \to
\mathbb R$).  The Lie algebra $\mathfrak A_P$ is divided into four ``quadrants'' according to
the value of $\Delta^{+}_P(e)$ (although we have not indicated which quadrant contains a given
``wall'').

Theorem \ref{thm-answer} states that the portion of $H^*(\mathfrak N_P)$ which contributes to
the Lefschetz number at this fixed point depends on the quadrant in which $t = \log(e)$ lies:
if $\Delta^{+}_P(e) = J \subset \Delta_P$ then only the portion of $H^*(\mathfrak N_P)$ which
lies in the quadrant indexed by $J$ contributes to the Lefschetz number.  A further degree
shift by $|J|$ occurs when this portion  $H^*(\mathfrak N_P)\qdr{J}$ is identified (in Theorem
\ref{thm-supports}) with the local weighted cohomology with supports. 

It is a remarkable
fact that, globally in the Hecke correspondence, the fixed points occur in Weyl group orbits.
Assuming $t$ is regular (does not lie on a wall) then, after summing over all the fixed
points, each chamber will appear the same number of times.  It is the sum of these local
contributions over a $W$-orbit of fixed points (\cite{GKM} p. 529, last paragraph) which gives
rise to the combinatorial formula for the averaged discrete series characters as described in
\cite{GKM}.

\subsection{Very positive and very negative weights }
Let $i:X \to \overline{X}$ denote the inclusion.  Suppose the weight profile $\nu = -\infty$
(or is very negative).  Then the weight truncation does nothing,
and the weighted cohomology sheaf $\mathbf{W^{\nu}C}^{\bullet}(\mathbf E) \cong Ri_*(\mathbf
E)$ becomes the ``full'' direct image of $\mathbf E.$  For any stratum $X_Q$,  $I_{\nu}(w) =
\phi$ for any $w \in W^1_Q.$  Theorem \ref{thm-answer} then says that a fixed point stratum $F
\cap C_Q$ (with characteristic element $e$) makes a contribution to the Lefschetz number only
if $\Delta_Q^{+}(e) = \phi$, which is to say, only if the Hecke correspondence is either
contracting or neutral in every direction normal to the stratum $X_Q.$  

In this case the local contribution to the Lefschetz number may be expressed in terms of the
character of the finite dimensional representation $G \to \text{GL}(E).$  We briefly recall 
the argument in \cite{GKM} \S 7.18.  The quantity $\sum_i(-1)^i \text{Tr}(e^{-1};
H^i(\mathfrak N_P,E))$ is equal to 
\[\text{Tr}(e^{-1};E)\] 
times the following quantity:
\begin{align*} 
 \sum_i(-1)^i \text{Tr}(e^{-1}; \wedge^i(\mathfrak N_P^*)
&=  \det(1-Ad(e); \mathfrak N_P(\mathbb C)) \\
&=  \prod_{\alpha \in \Phi^+_L}(1-\alpha^{-1}(e)) \prod_{\alpha \in
\Phi^+_L}\alpha(e)(-1)^{\dim \mathfrak N_P}\\
&=  \Delta_P(e) \det(e; \mathfrak N_P) (-1)^{\dim \mathfrak N_P}
\end{align*}
where $\Delta_P(e) = \prod_{\alpha \in \Phi^+_L} (1 - \alpha^{-1}(e))$ denotes the (partial)
Weyl denominator.  (These quantities may be further expressed in terms of $|D^G_L(e)|$,
$\delta_P(e)$, and $\chi_G(e)$ using \cite{GKM} (7.16.11), (7.18.3) and \cite{GKM} p. 497.)

Similarly, suppose the weight profile $\nu = +\infty$ or is a very large
positive weight.  The stalk cohomology (at a point $x \in X_Q$ in some boundary stratum $X_Q$)
of the weighted cohomology sheaf $\mathbf{W^{\nu}C}^{\bullet}(\mathbf E)$ vanishes because the
weight truncation (\ref{eqn-stalkcohomology}) kills everything.  In this case, the weighted
cohomology sheaf is quasi-isomorphic to the sheaf $Ri_!(\mathbf E)$ which is obtained as the
extension by 0 of the local system $E$.  Its cohomology is the compact support cohomology
$H^*_c(X;E)$ of the locally symmetric space.  For any stratum $X_Q$, According to
(\ref{eqn-IWeyl}), $I_{\nu}(w)=\Delta_Q$ for any $w \in W^1_Q.$  Theorem \ref{thm-answer} then
says that a fixed point stratum $F \cap C_Q$ (with characteristic element $e \in L_Q$) makes a
contribution to the Lefschetz number only if $\Delta_Q^{+}(e) = \Delta_Q$, that is, only if
the Hecke correspondence is strictly expanding in all directions normal to the stratum $X_Q.$  
Then the same quantity $\sum_i(-1)^i \text{Tr}(e^{-1}; H^i(\mathfrak N_P,E))$
occurs in the formula, but with a (possibly) different sign.

In these cases (of $\nu = \pm \infty$) the Lefschetz formula of Franke \cite{Franke} can be
recovered, cf. \cite{GKM} \S 7.17, 7.18.

\subsection{Hyperboic 3-space}\label{subsec-notzero}
For $\mathbf{G}(\mathbb R) = \text{SL}_2(\mathbb C)$ the symmetric space $D = G/K$ may be
identified with hyperbolic 3-space.  If $\Gamma$ is a torsion-free arithmetic group, then $X =
\Gamma \backslash D$ is a hyperbolic 3-manifold.  The group $G$ does not contain a compact
maximal torus.  Consequently, $\chi(X) = 0$ (cf \S \ref{thm-Harder}).  However, when $X$ is
not compact, there exists a Hecke correspondence on $X$ whose fixed point set consists of a
smooth curve which passes from one cusp to another cusp.  The Euler characteristic of this
fixed point set is not zero, although the Euler form vanishes identically\footnote{It does
not violate Theorem \ref{prop-Harder}:  it is not a ``locally symmetric space'' in the
sense of \S \ref{subsec-LSS} because it contains (and in fact consists of) a Euclidean
factor.}. It is possible to find particular weight profiles such that
the (global) Lefschetz number of this correspondence on the weighted cohomology is nonzero.
However, the formula \cite{GKM} (thm. 7.14.B) would attribute the contribution from this fixed
curve to the cusps, rather than to the interior stratum.  This re-attribution is a result of
equation (7.14.2) of \cite{GKM}.

\subsection{Nielsen fixed point theory}
Suppose $X$ is a compact manifold with fundamental group $\Gamma = \pi_1(X,x_0).$  Let $f:X
\to X$ be a self-map.  A choice of path from the basepoint $x_0$ to its image $f(x_0)$
determines a homomorphism $\phi:\Gamma \to \Gamma.$  Two elements $\gamma_1,\gamma_2 \in
\Gamma$ are said to be {\it $\phi$-conjugate} if there exists $\gamma \in \Gamma$ so that
$\gamma_2 = \gamma \gamma_1 \phi(\gamma)^{-1}.$  Let $(\Gamma)_{\phi}$ denote the set of
$\phi$-conjugacy classes in $\Gamma$ and let $\mathbb R (\Gamma)_{\phi}$ be the vectorspace of
finite formal linear combinations of such classes.  For each connected component $F$ of the
fixed point set of $f$, let $L(F)\in \mathbb R$ denote the contribution of $F$ to the
Lefschetz number $L$, that is,
\[ L=\sum_i(-1)^i\ \text{Tr }(f^*:H^i(X) \to H^i(X)) = \sum_F L(F).\]
The Nielsen theory (see \cite{Nicas}) assigns \begin{itemize}
\item a $\phi$-conjugacy class $\{ F \} \in (\Gamma)_{\phi}$ to each connected component $F$
of the fixed point set, and
\item a (cohomologically defined) Nielsen number $N(\{ \gamma \},f)$ to each $\phi$-conjugacy
class $\{ \gamma \}$
\end{itemize}
such that \begin{equation}\label{eqn-Nielsen}
 \sum_{\{\gamma\}\in (\Gamma)_{\phi}} N(\{\gamma\},f) \{ \gamma \} =
\sum_{F} L(F) \{F\} \in \mathbb R(\Gamma)_{\phi}\end{equation}
thereby ``refining'' the Lefschetz fixed point formula.  (The sum on the left is over
$\phi$-conjugacy classes and the sum on the right is over connected components of the fixed
point set.)

Now suppose that $X = \Gamma \backslash D$ is a {\em compact} locally symmetric space.  Fix $g
\in\mathbf G(\mathbb Q)$ and let $C \rightrightarrows X$ be the resulting Hecke
correspondence. Let $(\Gamma g \Gamma)_1$ be the set of $\Gamma$-conjugacy classes of elements
$e \in \Gamma g \Gamma.$  Let $\mathbf E$ to be the local system corresponding to a
representation $\tau: G \to \text{GL(E)}.$  Theorem \ref{thm-intro} then says that the
Lefschetz number of this correspondence is:
\begin{equation}\label{eqn-shortanswer} L = \sum_{\{e\}}
\chi(F(e))\text{tr}(\tau(e)^{-1};E).\end{equation}
Here, the sum is taken over all conjugacy classes $\{e\} \in (\Gamma g \Gamma)_1$, and $F(e)$
denotes the set of fixed points which have characteristic element equal to $e$.  This set is
empty unless $e$ is elliptic (modulo $A_G$).  If $F(e)$ is not empty, then it is compact.

It turns out that if the local system $E$ is trivial, and if the correspondence $C
\rightrightarrows X$ is actually a self-map $f:X \to X$ then the terms in
(\ref{eqn-shortanswer}) are exactly the terms in the Nielsen formula (\ref{eqn-Nielsen}). 
The group $\Gamma$ may be identified with the fundamental group $\pi_1(X,x_0).$
The Hecke correspondence is actually a self-map iff the element $g$ normalizes $\Gamma.$  In
this case, the automorphism $\phi:\Gamma \to \Gamma$ is given by conjugation:  $\phi(\gamma) =
g \gamma g^{-1}.$ Finally, the association $a \mapsto ag$ (for $a \in \Gamma$) determines a
one to one correspondence
\[ (\Gamma)_{\phi} \to (\Gamma g \Gamma)_1.\]
There is a slightly more general Nielsen formula for correspondences, also with coefficients
in a local system.  The terms in this formula again coincide with the terms in the sum
(\ref{eqn-shortanswer}).

\end{document}